\documentclass[11pt]{amsart}
\usepackage{amssymb}
\usepackage{amscd}
\usepackage{curves}
\usepackage{epsfig}
\usepackage[pass]{geometry}
\usepackage{graphicx}
\usepackage{pstricks}
\usepackage{pstricks-add}
\usepackage{pst-grad}
\usepackage{pst-plot}
\usepackage{verbatim}
\usepackage{latexsym,eucal,amsfonts,amssymb,amsmath,graphicx}
\usepackage{auto-pst-pdf}
\numberwithin{equation}{section}
\newtheorem{theorem}{Theorem}[section]

\newtheorem{lemma}[theorem]{Lemma}
\newtheorem{prop}[theorem]{Proposition}

\def \bpf {\begin{proof}}
\def \epf {\end{proof}}
\def \beq {\begin{eqnarray*}}
\def \eeq {\end{eqnarray*}}
\def \bsp{\begin{split}}
\def \esp{\end{split}}
\def \beqq {\begin{eqnarray}}
\def \eeqq {\end{eqnarray}}
\def \mca {{\mathcal A}}
\def \mcb {{\mathcal B}}

\def \mcd {{\mathcal D}}
\def \mce {{\mathcal E}}
\def \mcf {{\mathcal F}}
\def \mcg {{\mathcal G}}
\def \mch {{\mathcal H}}
\def \mci {{\mathcal I}}
\def \mcj {{\mathcal J}}

\def \mcl {{\mathcal L}}
\def \mcm {{\mathcal M}}

\def \mco {{\mathcal O}}
\def \mcp {{\mathcal P}}

\def \mcr {{\mathcal R}}

\def \mct {{\mathcal T}}
\def \mcu {{\mathcal U}}
\def \mcv {{\mathcal V}}
\def \mcw {{\mathcal W}}
\def \mcx {{\mathcal X}}
\def \mcy {{\mathcal Y}}

\def \mbn {{\mathbb N}}
\def \mbr {{\mathbb R}}

\def \id {\operatorname{Id}}
\def \comp {\operatorname{comp}}

\def \loc {\operatorname{loc}}
\def \div {\operatorname{div}}
\def \det {\operatorname{det}}

\def \tr {\operatorname{Tr}}
\def \ric {\textrm{Ric}}

\def \supp {\text{supp }}

\def\Id {\operatorname{Id}}
\def \eps {\epsilon}   
\def \la {\lambda}   
\def \La {\Lambda}

\def \p {\partial}

\def \eps {\epsilon}
\def \det {\text{det}}

\def \vol {\text{vol}}

\def \ha {\frac{1}{2}}

\def \fnf {\frac{n}{4}}
\def \WF {\operatorname{WF}}
\def \loc {\operatorname{loc}}
\def \ein{\operatorname{Ein}}

\def \bfq {\textbf{Q}}
\def \bfb {\textbf{B}}
\def \diag {\operatorname{Diag}}
\def \arv {\vec{v}}
\def \singsupp{\operatorname{singsupp}}
\def \div {\hbox{div}}
\def\bra{\langle}
\def\cet{\rangle}
\def\R{\mbr}

\begin{document}
\title{Determination of vacuum space-times from the Einstein-Maxwell equations}
\author{Matti Lassas}
\address{Matti Lassas
\newline
\indent Department of Mathematics, University of Helsinki}
\email{matti.lassas@helsinki.fi}

\author{Gunther Uhlmann} 
\address{Gunther Uhlmann
\newline
\indent Department of Mathematics, University of Washington,
\indent 
\newline  
\indent Institute for Advanced Study, the Hong Kong University of Science and Technology 
\newline
\indent and Department of Mathematics, University of Helsinki}
\email{gunther@math.washington.edu}

\author{Yiran Wang}
\address{Yiran Wang
\newline
\indent Department of Mathematics, University of Washington 
\newline
\indent and Institute for Advanced Study, the Hong Kong University of Science and Technology}
\email{wangy257@math.washington.edu}
\begin{abstract} 
We study inverse problems for the Einstein-Maxwell equations. We prove that it is possible to generate gravitational waves from the nonlinear interactions of electromagnetic waves. By sending electromagnetic waves from a neighborhood of a freely falling observer and taking measurements of the gravitational  perturbations in the same neighborhood, one can determine the vacuum space-time structure up to diffeomorphisms in the largest region where these waves can travel to from the observer and return. 
\end{abstract}
 
\maketitle

%==================================%
\section{Introduction}
This paper continues the study of inverse problems with sources for non-linear hyperbolic equations in 4 dimensional space-times  initiated in \cite{KLU1}, where Einstein equations are coupled with matter fields. It is shown in \cite{KLU1} that one can recover the topology, differentiable structure and conformal class of the metric in a larger set than  where the observations are made. The article \cite{LUW} studies a general class of semilinear wave equations extending the previous work \cite{KLU} that considered quadratic non-linearities. See also \cite{WZ} for equations with quadratic derivative nonlinear terms. In this paper we consider Einstein field equations in vacuum coupled with Maxwell equations. We show that just using electromagnetic sources we can generate gravitational waves and use them to determine the metric, topology and differentiable structure in a larger set than where the observations are made. The method of proof  uses asymptotic expansions and the precise study of the singularities of each term in the expansion. We use as in \cite{KLU1, KLU} four plane waves interacting to create  new point singularities and use these to determine the earliest light observation set. Then we appeal to the geometric result in \cite{KLU} that states that from the earliest light observation set we can determine the topology, differentiable structure and conformal class of the metric in a larger set than where the observations are made. The asymptotics needed in the current paper are more subtle than the previously mentioned works since the highest order term vanishes and one needs to look at the lower order terms carefully. As in the previously mentioned papers we remark that the linearized inverse problem is not known to be solvable. We use the non-linear interaction of waves in a significant way to create new point sources.

%In this work, we study the inverse problem of determining vacuum space-time structures by sending electromagnetic waves. 
In the following, we first introduce the mathematical model, then we formulate the inverse problem and state the main result. 
%%%%%%%%%%%
\subsection{The Einstein-Maxwell equations}
Let $M$ be a $4$ dimensional smooth manifold and $g$ be a Lorentzian metric on $M$ satisfying the vacuum Einstein equations    
\beq
\ein(g) = 0.
\eeq
Here $\ein(g)$ denotes the Einstein tensor given by 
\beq
\ein(g) = \ric(g) - \ha R(g)g
\eeq
where $\ric(g)$ denotes the Ricci curvature tensor and $R(g)$ the scalar curvature. In particular, $(M, g)$ is called a vacuum space-time. The propagation of electromagnetic waves on $(M, g)$ are governed by the Maxwell equations. In the covariant formulation, the electromagnetic field can be described as a two form $F$ and the (four) electric current $J$ is a vector field on $M$. The Maxwell equations for $F$ on $(M, g)$ with source $J$ are given by 
\beqq\label{maxwell}
dF = 0, \ \ \delta_g F = J^\flat,
\eeqq
where the codifferential $\delta_g$ is the dual of the exterior differential $d$ with respect to  $g$ and $J^\flat$ denotes the one form obtained from $J$ by lowering the index using $g$. In simply connected domains, one can find a one form $\phi$ by Poincar\'e lemma such that $F = d\phi$, and $\phi$ is called the electromagnetic potential. The Maxwell equations can be written as
\beq
\delta_g d\phi = J^\flat,
\eeq
while the first set of Maxwell equations is automatically satisfied as $d^2 = 0$. The source $J$ is subject to the conservation law $\div_g J = 0$, where $\div_g$ denotes the divergence operator. 

To understand the gravitational effects of incidenting electromagnetic waves on $(M, g)$ with electric current $J$ as a source, we need to couple the Einstein and Maxwell equations in a physically meaningful way. There are many models for the Einstein-Maxwell equations in the literature, for example equations with no electric current \cite[Section 6.10]{Cb} and \cite{MTW}, and equations for charged dust  \cite{Cb} and \cite[Chap.\ 18]{Tay3}. In this work, we derive the equations with sources in Section \ref{lagfor} by ignoring the contribution from the mass of the dust in the charged dust model. Roughly speaking, the total Lagrangian consists of the Einstein-Hilbert Lagrangian, the electromagnetic Lagrangian and the interaction term. The Einstein-Maxwell equations are the Euler-Lagrange equations of the total Lagrangian and they can be written as
\beq 
\begin{gathered}
\ein(g) = T_{sour},\\
\delta_g d\phi = J^\flat,
\end{gathered}
\eeq 
where $T_{sour}$ is the source stress-energy tensor which can be written as
\beqq\label{eqtsour}
\begin{gathered}
T_{sour} = T_{em} + T_{inter},\\
T_{em, \alpha\beta} = F_\alpha^\la F_{\beta\la} - \frac{1}{4}g_{\alpha\beta}F^{\la\mu}F_{\la\mu}, \ \
T_{inter} = -\ha (J^\mu \phi_\mu)g.
\end{gathered}
\eeqq
Here $T_{em}, T_{inter}$ are the stress-energy tensor of the electromagnetic field $F$ and the interaction terms respectively. 
Moreover,  $T_{sour}$ satisfies the conservation law $\hbox{div}_gT_{sour} = \nabla_\alpha T^{\alpha\beta} = 0$ if the Maxwell equations are satisfied. This puts the requirement that the four current $J$ should satisfy the conservation law $\div_g J = 0$. Now we must be careful formulating the source problem, because it is not possible to solve the Einstein-Maxwell equations for every four current $J$. This can be seen easily from the Maxwell equations on Minkowski space-time. In this case, the component of the four current are $J^0 = q$ and $J^i =  \jmath^i, i = 1, 2, 3,$ where $q$ is the electric charge and $\jmath$ is the electric current. The conservation law is just the conservation of charges 
\beq
\frac{dq}{dt} + \sum_{i = 1}^3 \frac{\p \jmath^i}{\p x^i} = 0, 
\eeq
where $(x^i), i = 0, 1,2, 3$ denotes the coordinates of $\mbr^4$, see for example, \cite[Section II.3.2]{Cb} and \cite[Chap.\ 4]{MTW}. So once $\jmath$ is prescribed,  $q$ is determined from the above equation. Moreover, even if $\jmath$ is smooth and compactly supported, $q$ may not be compactly supported. Therefore, one should think of the conservation law as a constitute equation for the Einstein-Maxwell equations with sources, and construct certain component of the source $J$ given the others. We remark that a similar problem for the Einstein equations with matter fields is addressed by Girbau and Bruna \cite{GB} for the Cauchy problem.

Now we formulate the local problem of Einstein-Maxwell equations with sources to be considered in this work, starting with some notions of Lorentzian geometry. For $p, q$ on a Lorentzian manifold $(M, g)$, we denote by $p\ll q$ ($p < q$) if $p\neq q$ and there is a future pointing time-like (causal) curve from $p$ to $q$. We denote by $p\leq q$ if $p = q$ or $p<q$. The chronological (causal) future of $p\in M$ is the set $I_g^+(p) = \{q\in M: p\ll q\}$ ($J_g^+(p) = \{q\in M: q\leq p\}$). Similarly, we can define the chronological past and causal past, which are denoted by $I_g^-(p)$ and $J_g^-(p)$ respectively. For any set $A\subset M$, we denote the causal future by $J_g^\pm(A) = \bigcup_{p\in A}J_g^\pm(p)$. Also, we denote $J_g(p, q) = J_g^+(p)\cap J_g^-(q)$ and $I_g(p, q) = I_g^+(p)\cap I_g^-(q)$.

We consider globally hyperbolic manifold $(M, g)$, which can be identified with the product manifold $\mbr\times \mcm$ with $\mcm$ a $3$-dimensional manifold and metric $g = -\beta(t, y) dt^2 + \kappa(t, y)$ where $\beta > 0$ is a smooth function and $\kappa$ is a family of Riemmanian metrics on $\mcm$ smooth in $t$, see Section \ref{loren}. For $T_0\in \mbr$, we denote $M(T_0) = (-\infty,T_0)\times \mcm$. For a vector field $J$ on $M$, we shall write $J = (J^0, \bar J)$ where $\bar J \in C^\infty(\mbr_t; T\mcm)$. In particular, $\bar J$ is a family of smooth vector fields on $\mcm$ and smooth in $t\in \mbr$. 

We consider $(M, \widehat g)$ as a background space-time where $\widehat g$ satisfies the vacuum Einstein equations. Let $\widehat \mu(t) \subset M$ be a time-like geodesic where $t\in [-1, 1]$. In the language of general relativity, $\widehat \mu$ represents a freely falling observer. Let $V\subset M$ be an open relatively compact neighborhood of $\widehat \mu([s_-, s_+])$ where $-1<s_-<s_+<1$. Take $T_0 >0$ such that $V\subset M(T_0)$. For $\bar J$ compactly supported in $V$, the Einstein-Maxwell equations with sources we consider are following equations for the field $(g, \phi, J^0)$ 
\beqq\label{einmax}
\begin{gathered}
\left\{\begin{array}{c}
\ein(g) = T_{sour}(g, \phi, J) \\ 
\delta_g d\phi = J^\flat\\
\div_g J = 0
\end{array}\right. \text{ in } M(T_0),\\ 
g = \widehat g, \ \ \phi = 0, \ \ J^0 = 0,  \text{ in } M(T_0)\backslash J_g^+(\supp (\bar J)). 
\end{gathered}
\eeqq
As is well-known, the solution of the Einstein equation is unique only up to diffeomorphisms and the electric potential is unique up to a gauge choice. Hence one should think of $g, \phi$ and $J$ in the above equations as representatives of their equivalence classes.  

In Section \ref{secwell}, we study the well-posedness of \eqref{einmax} in wave and Lorentz gauge. Moreover, in Section \ref{gwave}, we consider the gravitational perturbations in the weak field approximation and show that small electric current $\bar J$ will generate non-trivial gravitational waves.  We remark that the coupling of progressive electromagnetic waves and gravitational waves is known for the electro-vacuum space-time (with no source) using the WKB method, see for example Choquet-Bruhat \cite[Chap.\ XI]{Cb}. Also, we'd like to point to the interesting work by F\"uzfa \cite{Fu}, who demonstrated numerically that strong electromagnetic fields (current loops and solenoid) can generate gravitational field curving the space-time. 

\subsection{The inverse problem and the main results}
\begin{figure}[htbp]
\centering
% Generated with LaTeXDraw 2.0.1
% Sat Oct 31 18:06:19 CST 2015
% \usepackage[usenames,dvipsnames]{pstricks}
% \usepackage{epsfig}
% \usepackage{pst-grad} % For gradients
% \usepackage{pst-plot} % For axes
\scalebox{0.7} % Change this value to rescale the drawing.
{
\begin{pspicture}(0,-4.68)(6.8418946,4.68)
\psellipse[linewidth=0.04,dimen=outer](3.11,-0.13)(0.55,3.69)
\psbezier[linewidth=0.04](3.0557144,-4.66)(3.26,-3.04)(3.2428572,-2.26)(3.1114285,-0.18)(2.98,1.9)(3.0557144,3.58)(3.2414286,4.66)
\psellipse[linewidth=0.04,linestyle=dashed,dash=0.16cm 0.16cm,dimen=outer](3.1,-0.23)(3.1,0.61)
\psdots[dotsize=0.12](3.08,2.84)
\psdots[dotsize=0.12](3.18,-3.18)
\psline[linewidth=0.04cm,linestyle=dashed,dash=0.16cm 0.16cm](3.12,2.84)(6.1,-0.08)
\psline[linewidth=0.04cm,linestyle=dashed,dash=0.16cm 0.16cm](3.04,2.84)(0.08,-0.1)
\psline[linewidth=0.04cm,linestyle=dashed,dash=0.16cm 0.16cm](6.16,-0.3)(3.18,-3.2)
\psline[linewidth=0.04cm,linestyle=dashed,dash=0.16cm 0.16cm](0.04,-0.3)(3.16,-3.2)
\usefont{T1}{ptm}{m}{n}
\rput(4.5,2.965){$p_+ = \widehat\mu(s_+)$}
\usefont{T1}{ptm}{m}{n}
\rput(4.5,-3.255){$p_- = \widehat\mu(s_-)$}
\usefont{T1}{ptm}{m}{n}
\rput(1.5614551,-0.19){$I_{\widehat g}(p_-, p_+)$}
\usefont{T1}{ptm}{m}{n}
\rput(2.95,-2.235){$V$}
\usefont{T1}{ptm}{m}{n}
\rput(3.5071455,4.385){$\widehat\mu$}
\end{pspicture} 
}
\caption{$\widehat\mu$ is a time-like geodesic (the observer). The set $V$ is an open neighborhood of $\widehat\mu([-1, 1])$. The source perturbation $J$ is supported in $V$ and produces gravitational and electromagnetic perturbations. The measurements are taken in $V$. We study the inverse problem of determining the space-time in $I_{\widehat g}(p_-, p_+)$ bounded by the dashed curves.}
\label{figinv}
\end{figure}
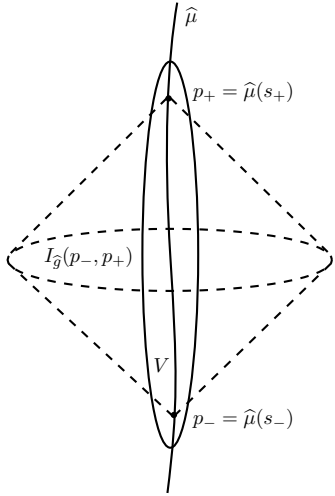

Assume that we are given $V$ a neighborhood of an observer as a differentiable manifold. Suppose we can control  the electric current $\bar J$ which are compactly supported in $V$ and we take the measurements of the gravitational perturbations $g$ in $V$.  Can we determine the background metric $\widehat g$ on $I_{\widehat g}(p_-, p_+)$ using these information? See Fig.\ \ref{figinv}. This problem becomes particularly interesting thanks to the breakthrough discovery by the LIGO project \cite{Ab, Ab1}, i.e.\ direct detection of  gravitational waves. We remark that the region $I_{\widehat g}(p_-, p_+)$ is the maximal region where the electromagentic waves can travel to from $V$ and from where the waves can travel back to $V$. Therefore, it is the maximal region that we can expect to determine the space-time structure in this setup. 

To formulate the inverse problem precisely, one need to keep in mind that the solutions to the Einstein-Maxwell equations are unique up to diffeomorphisms of certain regularity. Here we follow \cite{KLU1} to formulate the observations in Fermi coordinates associated with freely falling observers $\widehat \mu([-1, 1])$. Assume that $p = \widehat \mu(-1)\in \{0\}\times \mcm$. We shall take some basis $X_j(p) \in T_p M, j = 0, 1,2, 3$ such that $X_0(p) = \dot{\widehat \mu}(-1)$ is time-like. We let $Z_j, j = 0, 1, 2, 3$ be a frame along $\widehat \mu$ obtained from parallel translating $X_j(p)$. Then we consider the following Fermi coordinates:
\beqq\label{Fermi coordinates}
\Phi_{g, p}(z^0, z^1, z^2, z^3) = \exp_{g, \widehat\mu(z^0)}(\sum_{j = 1}^3 z^j Z_j), \ \ \Phi_{g, p}: \widetilde V_p  \rightarrow V^g \subset V,
\eeqq 
where $\widetilde V_p = (-1, 1)\times B_p$ is a subset of $\mbr^4$, $B_p$ is a small neighborhood of $0\in \mbr^3$ depending on $p$ (or $\widehat \mu$) so that the exponential maps are well-defined and $V^g = \Phi_{g,p}(\tilde V_p)$. In particular,  $(\widetilde V_p, \Phi_{g, p})$ gives a coordinate system of $V^g$, the Fermi coordinates. We formulate the problem in this coordinate. 

For the inverse problem, it is natural to think of the four current $J$ as the source and define the source-to-solution map of the Einstein-Maxwell system as $J \rightarrow (g, \phi)$. Moreover, it is convenient to take the graph of the source-to-solution map as the data set. On the manifold $M$, we introduce a complete Riemannian metric $\widehat g^+$ so we can define seminorms in $C^s(M)$.  For $\delta>0$ small, we define the data set as
\beqq\label{dataset}
\begin{gathered}
 \mcd(\delta) \doteq \{(\Phi_{g, p}^*g|_{\widetilde V_p}, \Phi_{g, p}^*\phi|_{\widetilde V_p}, \Phi_{g, p}^*J|_{\widetilde V_p}):  J = (J^0, \bar J) \in C^{4}(M), 
 \| \bar J\|_{C^{4}(M)}< \delta, \\
 \supp (\bar J) \subset V^g  \text{ and } (g, \phi, J^0) \in C^4(M) \text{ satisfy \eqref{einmax}}\}.
 \end{gathered}
\eeqq
Roughly speaking, the data set $\mathcal{D}(\delta)$ corresponds to the following measurement settings: We control the electric current $\overline J$ that is supported in the neighborhood $V^g$ of the geodesic $\hat \mu$. This gives rise to the moving electric charges $J^0$. Together $J=(J^0,\overline J)$ cause a perturbation in  the electric field $\phi$ and the gravitational field $g$ that start to propagate in space time. This perturbed fields can be considered as a non-linear wave that scatters from non-homogenous background metric $\widehat g$. Some of the scattered waves return back to $V^g$ where we observe the fields $(g,\phi)$. In the data set $\mcd(\delta)$ the controlled source $\overline J$ and the other fields $J^0,g,\phi$ are given in the Fermi coordinates $\Phi_{g,p}$ associated to the freely falling observer (time-like geodesic) $\mu_g$. Our main result is the unique determination of two vacuum space-times up to isometries. 
\begin{theorem}\label{main1}
Let $M$ be a $4$ dimensional simply connected smooth manifold and $\widehat g^{(i)}, i = 1, 2$ be two globally hyperbolic Lorentzian metrics on $M$  satisfying the vacuum Einstein equations.  Let $\widehat \mu^{(i)}([-1, 1])$ be time like geodesics with respect to $\widehat g^{(i)}$. Assume that $p^{(i)} = \widehat\mu^{(i)}(-1)\in \{0\}\times \mcm$ and take $p^{(i)}_\pm = \widehat \mu^{(i)}(s_\pm)$ with $-1<s_-<s_+ < 1$.  Suppose that $V^{(i)}$ are neighborhoods of $\widehat \mu([s_-, s_+])$  and $V^{(i)}\subset M(T_0)\backslash M(0)$ for $T_0>0$. Consider the Einstein-Maxwell systems  
\beq
\begin{gathered}
\left\{\begin{array}{c}
\ein(g^{(i)}) = T_{sour}(g^{(i)}, \phi^{(i)}, J^{(i)}) \\[3pt]
\delta_{g^{(i)}} d\phi^{(i)} =  J^{(i), \flat}
\end{array}\right. \text{ in } M(T_0),\\[3pt]
g^{(i)} = \widehat g^{(i)}, \ \ \phi^{(i)} = 0,  \text{ in } M(T_0)\backslash J_{g^{(i)}}^+(\supp (\bar J^{(i)})),
\end{gathered}
\eeq
which is well-posed on data set $\mcd^{(i)}(\delta), i = 1, 2$ defined as \eqref{dataset} for some small $\delta$. If 
\beq
\mcd^{(1)}(\delta) = \mcd^{(2)}(\delta),
\eeq
then there is a diffeomorphism $\Psi: I_{\widehat g^{(1)}}(p^{(1)}_-, p^{(1)}_+)\rightarrow I_{\widehat g^{(2)}}(p^{(2)}_-, p^{(2)}_+)$ such that $\Psi^*\widehat g^{(2)} = \widehat g^{(1)}$ in $I_{\widehat g^{(1)}}(p^{(1)}_-, p^{(1)}_+)$. 
\end{theorem}

%This type of inverse problem was introduced by Kurylev, Lassas and Uhlmann in \cite{KLU} for semilinear wave equations with a quadratic nonlinear term and the extended preprint \cite{KLU1} for the Einstein-scalar field equations. Roughly speaking, their work showed that the topological, differentiable structure as well as the conformal class of the metric in $I_{\widehat g}(p_-, p_+)$ can be determined from the so-called source-to-solution map of associated equations. If the manifolds are Ricci flat, the metric can be determined up to isometries. The results are extended and improved to equations with general semilinear terms with no derivatives in Lassas-Uhlmann-Wang \cite{LUW} and with quadratic derivative terms in Wang-Zhou \cite{WZ}. 

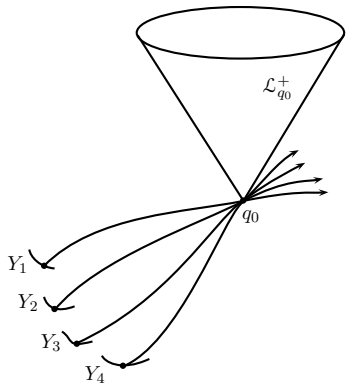
\begin{figure}[htbp]
% Generated with LaTeXDraw 2.0.1
% Thu Dec 03 16:32:23 CST 2015
% \usepackage[usenames,dvipsnames]{pstricks}
% \usepackage{epsfig}
% \usepackage{pst-grad} % For gradients
% \usepackage{pst-plot} % For axes
\scalebox{.7} % Change this value to rescale the drawing.
{
\begin{pspicture}(0,-3.5188477)(6.0210156,3.4788477)
\psellipse[linewidth=0.04,dimen=outer](4.0310154,2.9688478)(1.99,0.51)
\psline[linewidth=0.04cm](2.0810156,2.9188476)(4.0610156,-0.18115234)
\psline[linewidth=0.04cm](6.0010157,2.9388475)(4.1010156,-0.20115234)
\psdots[dotsize=0.12](4.0810156,-0.22115235)
\psbezier[linewidth=0.04,arrowsize=0.05291667cm 2.0,arrowlength=1.4,arrowinset=0.4]{<-}(5.6010156,0.19884765)(5.321016,0.13884765)(4.9010158,0.19884765)(4.1010156,-0.22115235)(3.3010156,-0.6411523)(1.3010156,-1.3211523)(0.5010156,-2.3011522)
\psbezier[linewidth=0.04,arrowsize=0.05291667cm 2.0,arrowlength=1.4,arrowinset=0.4]{<-}(5.2610154,0.49884766)(5.0410156,0.37884766)(4.619152,0.18169029)(4.0579753,-0.23587316)(3.4967985,-0.6534366)(2.6610155,-2.1411524)(0.90101564,-2.9411523)
\psbezier[linewidth=0.04,arrowsize=0.05291667cm 2.0,arrowlength=1.4,arrowinset=0.4]{<-}(5.7010155,-0.061152343)(5.341016,-0.08115234)(5.095066,-0.021152344)(4.066549,-0.23026125)(3.0380318,-0.43937016)(1.2678777,-0.47876698)(0.30101562,-1.4611523)
\psbezier[linewidth=0.04,arrowsize=0.05291667cm 2.0,arrowlength=1.4,arrowinset=0.4]{<-}(5.1410155,0.7388477)(4.7610154,0.5388477)(4.4010158,0.078847654)(4.0410156,-0.28115234)(3.6810157,-0.6411523)(2.8410156,-2.9611523)(1.7810156,-3.3811524)
\psdots[dotsize=0.12](0.30101562,-1.4611523)
\psdots[dotsize=0.12](0.5010156,-2.2811522)
\psdots[dotsize=0.12](0.9210156,-2.9411523)
\psdots[dotsize=0.12](1.8010156,-3.3611524)
\psbezier[linewidth=0.04](0.021015625,-1.15)(0.041015625,-1.29)(0.22101563,-1.49)(0.5010156,-1.51)
\psbezier[linewidth=0.04](0.30101562,-2.05)(0.36101562,-2.27)(0.5610156,-2.35)(0.8810156,-2.23)
\psbezier[linewidth=0.04](0.64101565,-2.71)(0.87301564,-2.8605883)(0.6990156,-3.03)(1.2210156,-2.8982353)
\psbezier[linewidth=0.04](1.4010156,-3.1670732)(1.4410156,-3.27)(1.7410157,-3.49)(2.3010156,-3.23)
\usefont{T1}{ptm}{m}{n}
\rput(4.2324707,-0.51615233){$q_0$}
\usefont{T1}{ptm}{m}{n}
\rput(-0.2,-1.4361523){$Y_1$}
\usefont{T1}{ptm}{m}{n}
\rput(0,-2.1561522){$Y_2$}
\usefont{T1}{ptm}{m}{n}
\rput(0.4324707,-2.9161522){$Y_3$}
\usefont{T1}{ptm}{m}{n}
\rput(1.2724707,-3.52961524){$Y_4$}
\usefont{T1}{ptm}{m}{n}
\rput(4.7424707,1.9238477){$\mcl^+_{q_0}$}
\end{pspicture} 
}
\caption{Interaction of singularities. Four electric current sources are placed on submanifolds $Y_i$, which produce four electromagnetic waves (distorted plane waves) propagating along geodesics from $Y_i$. Their nonlinear interactions at $q_0$ produce new gravitational waves propagating on the light-cone $\mcl^+_{q_0} = \exp_{q_0}(L^+_{q_0}M)$.}%\cup \{q_0\}$.}
\label{fourinter}
\end{figure}

A key ingredient in the proof of Theorem \ref{main1}  is to produce artificial point sources in $I_{\widehat g}(p_-, p_+)$ from the nonlinear interaction of linear waves. For the Einstein-Maxwell equations we study in this paper, the linear waves are electromagnetic waves caused by small perturbations of  the electromagnetic sources. In Section \ref{singu}, we show that the nonlinear interaction of four electromagnetic waves could produce a point source of gravitational waves and this helps us to solve the inverse problem, see Fig.\ \ref{fourinter}.

We remark that the propagation and interaction of singularities for nonlinear hyperbolic equations was actively studied in the 80's and 90's, mainly for $2+1$ dimension by Bony \cite{Bo}, Melrose-Ritter \cite{MR}, Rauch-Reed \cite{RR}, etc, see also Beals \cite{Bea} for an overview. The Einstein-Maxwell equations reduces to a quasilinear hyperbolic system in $3+1$ dimensions in wave and Lorentz gauge, for which the methods in the mentioned literatures do not apply directly. We use the method developed in \cite{KLU} and \cite{LUW} to analyze the singularities in the asymptotic expansion of solutions with sources depending on some small parameters. We point out that we only treat weak conormal type of singularities. For the propagation and interaction  of strong gravitational waves, for example impulse waves, see \cite{LR0, LR}. 

The paper is organized as following. In Section \ref{secwell}, we formulate the source problem for the Einstein-Maxwell equations using Lagrangians. Then we discuss the well-posedness of the Einstein-Maxwell equations, first for the reduced equations in wave and Lorentz gauge then the full equations. In Section \ref{einlin}, we consider the linearization of the Einstein-Maxwell equations and construct distorted plane waves. Before solving the inverse problem, we show in Section \ref{gwave} the generation of gravitational waves using electromagnetic sources in wave gauge. In Section \ref{singu}, we study the leading singularities generated from nonlinear interaction of linear waves and prove that they are not always vanishing. Finally, we solve the inverse problem in Section \ref{secinv}.

\section*{Acknowledgments} 
The authors would like to thank Peter Hintz for valuable comments on a preliminary version of the paper and Judith Arms for helpful discussions. Matti Lassas has been supported by Academy of Finland. Gunther Uhlmann is partially supported by NSF, a Si-Yuan Professorship of HKUST and FiDiPro Professorship of Academy of Finland.

%==================================%
\section{Locall well-posedness of the Einstein-Maxwell equations}\label{secwell}

%%%%%%%%%%%
\subsection{Notations}\label{loren}
We review some facts from Lorentzian geometry and set up notations. Some of them are briefly discussed in the introduction and we collect them here for convenience. Our main references are \cite{BEE, Cb, KLU1}. 

Assume that $(M, g)$ is an $4$ dimensional Lorentzian manifold which is time oriented and globally hyperbolic. From \cite{BS0}, we know that $(M, g)$ is globally hyperbolic if there is no closed causal paths in $M$ and for any $p, q\in M$ and $p<q$, the set $J_g(p, q)$ is compact. We take the signature of the metric $g$ as $(-, +, +, +)$. It is proved in \cite{BS} that $(M, g)$ is isometric to a product manifold $\mbr\times \mcm$ with $g =-\beta(t, y)dt^2 + \kappa(t, y),$ where $\mcm$ is a $3$ dimensional manifold, $\beta: \mbr\times \mcm \rightarrow \mbr_+$ is smooth and $\kappa$ is a Riemannian metric on $\mcm$ and smooth in $t$. In this work, we identify $(M, g)$ with this isometric image. We shall use $x = (t, y) = (x^0, x^1, x^2, x^3)$ as the local coordinates on $M$. For $T_0\in\mbr$, we set $M(T_0) = (-\infty, T_0)\times \mcm.$ It is worth noting that for each $t\in \mbr$, the submanifold $\{t\}\times \mcm$ is a Cauchy surface, i.e.\ every inextendible causal curve intersects the submanifold only once, see for example \cite[Page 65]{BEE}. 

For $p\in M$, we denote the collection of light-like vectors at $p$ by $L_pM = \{\theta \in T_pM\backslash \{0\}: g(\theta, \theta) = 0\}$ and the bundle by $LM = \bigcup_{p\in M}L_pM$. The future (past) light-like vectors are denoted by $L^+_pM$ ($L^-_pM$), and the bundle $L^\pm M = \bigcup_{p\in M}L_p^\pm M$. Also, the set of light-like covectors at $p$ is denoted by $L^*_pM$ and the bundles $L^{*}M, L^{*,\pm}(M)$ are defined similarly. Since the metric $g$ is nondegenerate, there is a natural isomorphism $i_p: T_pM\rightarrow T_p^*M$. With this isomorphism, we sometimes use vectors and co-vectors interchangeably. Let $\exp_p: T_pM  \rightarrow M$ be the exponential map. The geodesic from $p$ with initial direction $\theta$ is denoted by $\gamma_{p, \theta}(t) = \exp_p (t\theta), t\geq 0$. The forward light-cone at $p\in M$ 
\beq
\mcl^+_p = \{\gamma_{p, \theta}(t): \theta \in L^+_pM, t> 0\}
\eeq 
is a submanifold of $M$ and we notice that $p\notin \mcl^+_pM$. 
 
We recall some useful formulas. We use the standard Einstein summation notation, i.e.\ summation is over repeated indices. In local coordinates, the Christoffel symbols are given by
\beq
\Gamma_{\alpha\beta}^\mu = \ha g^{\mu\la}( \p_{\beta} g_{\la\alpha}  +  \p_\alpha g_{\la\beta} -  \p_\la g_{\alpha\beta}), \ \  \p_i = \frac{\p }{\p x^i}, i = 0, 1, 2, 3.
\eeq
The contracted Christoffel symbol is 
\beq
\Gamma^\mu = g^{\alpha\beta}\Gamma_{\alpha\beta}^\mu = -|\det g|^{-\ha} \p_\la(|\det g|^\ha g^{\la\mu}),
\eeq
see for example \cite[Section 3]{FM} and \cite{KLU1}. The Laplace-Beltrami operator defined via $\square_g = \delta_g d + d\delta_g$ can be written as 
\beqq\label{eqbeltra}
\square_g = -|\det g|^{-\ha}  \p_\beta(|\det g|^\ha g^{\beta\alpha}  \p_\alpha) = -g^{\alpha\beta}  \p_\alpha\p_\beta + \Gamma^\alpha  \p_\alpha.
\eeqq
Let $J$ be a smooth vector field. The divergence of $J$ is given in local coordinates by 
\beqq\label{divJ}
\div_g J = \sum_{i = 0}^3 \frac{\p}{\p x^i} ((-\det g)^\ha J^i),
\eeqq
see for example \cite{Tay1} and \cite[p.\ 222]{MTW}. Let $\ric(g)$ denote the Ricci curvature tensor. For any covariant $2$-tensor $T$, recall that the trace of $T$ with respect to $g$ is defined as $\tr_g(T) = g^{\alpha\beta}T_{\alpha\beta}$. In particular, the scalar curvature $R(g)=\tr_g(\ric(g))$.  

Besides the Lorentzian metric, we will take a complete Riemannian metric $\widehat g^+$ on $M$, whose existence is guaranteed by \cite{NO}. With this metric, we can introduce distances on $M, TM$, and Sobolev spaces on $M$.  

%%%%%%%%%%%
\subsection{The Lagrangian formulation}\label{lagfor}
We derive the coupled Einstein-Maxwell equations from the Lagragian point of view. One can compare this with the equations for the charged dust, see Taylor \cite[Chap.\ 18]{Tay3}. Recall that the (vacuum) Einstein equations are the Euler-Lagrange equations of the Einstein-Hilbert Lagrangian:
\beq
\mcl_{grav}(g) = \int_M R(g) dg,
\eeq
see for example \cite[Theorem 7.1, Chap.\ III]{Cb}. Here $dg$ denotes the volume element and in local coordinates, $dg = |\det g|^\ha dx$.  Next, the electromagnetic field can be described as a two form $F$ on $(M, g)$ i.e.\ $F = F_{\alpha\beta} dx^\alpha\wedge dx^\beta.$ On a simply connected domain, we can write $F = d\phi \text{ or equivalently } F_{\alpha\beta} = \p_{\alpha}\phi_\beta - \p_{\beta} \phi_{\alpha}$ 
 with a one form $\phi$, the electromagnetic potential. The Lagrangian of the electromagnetic field is given by
\beq
\mcl_{em} =  \frac{1}{4}\int_M F^{\alpha\beta}F_{\alpha\beta} dg.
\eeq
In case there are electric current $J$ a vector field, we should add the interaction Lagrangian 
\beq
\mcl_{inter}  = - \int_M J^\mu \phi_\mu dg, 
\eeq
so the source Lagrangian is $\mcl_{sour} = \mcl_{em}+\mcl_{inter}$. Indeed, the Maxwell equations \eqref{maxwell} with source $J$ are the Euler-Lagrange equations of this Lagrangian. 

Now we consider the coupled Lagrangian of gravitational field and electromagnetic field:
\beq
\mcl_{total} = \mcl_{grav} + \mcl_{sour}.
\eeq
According to \cite[Theorem 7.3, Chap.\ III]{Cb}, we can write the first variation of $\mcl_{sour}$ as
\beq
\delta \mcl_{sour} = \int_M \Psi(g, \phi) \delta \phi dg - \int_M T_{sour} \delta g dg.
\eeq
(Note here $\delta$ denotes the variation and is different from the codifferential $\delta_g$.) Then we get the field equations i.e.\ Maxwell equations $\delta_g F = J^\flat,$ from $\Psi(g, \phi) = 0$ and  the Einstein equations with source stress-energy tensor $T_{sour}$ i.e.\ $
\ein(g) = T_{sour}$ from $\delta \mcl_{grav} = \int_M T_{sour} \delta g dg$. 
More importantly, since the Lagrangian $\mcl_{sour}$ is invariant under diffeomorphisms, Theorem 7.3 of \cite[Chap.\ III]{Cb} guarantees that 
the stress-energy tensor $T_{sour}$ satisfies the conservation law $\div_g T_{sour} = 0$ if the Maxwell equations are satisfied. It also follows from the Maxwell equations that the electric current satisfies the conservation law $\div_g J = 0$.  A more physical elaboration of electrodynamics on curved space-time, especially the validation of the conservation laws, can be found in \cite[Section 22.4]{MTW}. 

We find $T_{sour}$ explicitly. We can write $T_{sour} = T_{em} + T_{inter}$, where 
\beqq\label{stau}
T_{em, \alpha\beta} = F_\alpha^\la F_{\beta\la} - \frac{1}{4}g_{\alpha\beta}F^{\la\mu}F_{\la\mu}
\eeqq
is the stress-energy tensor of the electromagnetic field $F$. To find the stress-energy tensor $T_{inter}$ of $\mcl_{inter}$, we consider the first variation 
\beq
\delta \mcl_{inter} = -\int_M J^\mu \phi_\mu  \delta(dg)= \int_M  \ha  J^\mu \phi_\mu   g (\delta g) dg,
\eeq
since $\delta(dg) = -\ha g_{\alpha\beta}\delta g^{\alpha\beta}dg$, see formula (7.4) of \cite[Chap.\ III]{Cb}. So the stress-energy tensor $T_{inter} = -\ha (J^\mu \phi_\mu)g$. Finally, we can write the coupled Einstein-Maxwell equations with sources as
\beqq\label{einmaxsour}
\begin{gathered}
\ein(g) =   T_{em}  -\ha (J^\mu \phi_\mu)g, \\
\delta_g d\phi = J^\flat.
\end{gathered}
\eeqq
We remark that in absence of the electric current i.e.\ $J = 0$, these equations  reduce to the sourceless Einstein-Maxwell equations in literatures, see e.g.\ \cite{Cb, MTW}. 

For our analysis, it is convenient to use an equivalent form of the Einstein equations. We take the trace of the Einstein equation in \eqref{einmaxsour} to obtain
\beq
\tr_g(\ein(g)) = -R(g) =  -2 J^\mu \phi_\mu,
\eeq
because $T_{em}$ is trace-free in dimension $4$. Hence the Einstein equations in \eqref{einmaxsour} is equivalent to 
\beqq\label{ein1}
\ric(g) = T_{em} + \ha(J^\mu\phi_\mu)g.
\eeqq 

%%%%%%%%%%%%%%%%%%
\subsection{The reduced Einstein-Maxwell equations}\label{secwmap}
The Einstein equations form a second order system of PDEs. In harmonic gauge i.e.\ local coordinates such that $\Gamma^{\alpha} = 0, \alpha = 0, 1, 2, 3$, the Einstein equation \eqref{ein1} becomes a second order quasi-linear hyperbolic system for $g$. Actually, in local coordinates, the Ricci tensor can be written as
\beqq\label{richar}
\begin{gathered}
\ric_{\mu\nu}(g) = \ric^{(h)}_{\mu\nu}(g) + \ha (g_{\mu q}\p_\nu \Gamma^q  + g_{\nu q} \p_\mu \Gamma^q), \ \ \mu, \nu = 0, 1, 2, 3,
\end{gathered}
\eeqq
 where  $\ric^{(h)}_{\mu\nu}(g) = -\ha g^{pq} \p_p \p_q g_{\mu\nu} + P_{\mu\nu}$ is the harmonic Ricci tensor and
\beq
 P_{\mu\nu } = g^{ab}g_{ps}\Gamma^p_{\mu b}\Gamma^s_{\nu a} + \ha (\p_a g_{\mu\nu} \Gamma^a + g_{\nu l} \Gamma^l_{ab} g^{aq}g^{bd}  \p_\mu g_{qd} + g_{\mu l} \Gamma^l_{ab} g^{aq} g^{bd}  \p_\nu g_{qd}),
\eeq
see for example \cite[Appendix A]{KLU1} and \cite[Section VI.7]{Cb}. Note that $P_{\mu\nu}$ is a polynomial of $g_{pq}, g^{pq}$ and the first derivatives of $g_{pq}$. When $\Gamma^{\alpha} = 0$, $\ric(g)$ is reduced to $\ric^{(h)}(g)$. 

This reduction can be made coordinate invariant using wave map. Let $\widehat g$ be a smooth Lorentzian metric and $g'$ be a Lorentzian metric close to $\widehat g$ in $C^m(M)$ with $m$ to be specified later. Consider the wave map $f: (M, g')\rightarrow (M, \widehat g)$ satisfying the following equation 
\beqq\label{wavemap}
\begin{gathered}
\square_{g', \widehat g}f = 0 \text{ in } M(T_1),\\
f = \id \text{ in } M(0), 
\end{gathered}
\eeqq 
where $T_1>T_0$ and $\square_{g', \widehat g}$ is the wave map operator  defined independent of local coordinates, see \cite[Section VI.7]{Cb} and \cite[Appendix A]{KLU1} for more details. In local coordinates $x = (x^0, x^1, x^2, x^3)$ for $(M, g')$ and $y = (y^0, y^1, y^2, y^3)$ for $(M, \widehat g)$ so that $y^A = f^A(x), A = 0, 1,2, 3$, we have  
\beq
(\square_{g', \widehat g} f)^A(x) = (g')^{jk}(x)(\frac{\p^2 }{\p x^j \p x^k} f^A(x) - \Gamma'^n_{jk}(x) \frac{\p }{\p x^n}f^A(x) + \widehat \Gamma^A_{BC}(f(x)) \frac{\p}{\p x^j}f^B(x)\frac{\p}{\p x^k} f^C(x)),
\eeq
where $\Gamma', \widehat \Gamma$ denote the Christoffel symbols on $(M, g')$ and $(M, \widehat g)$ respectively. In particular, \eqref{wavemap} is a second order semilinear hyperbolic system for $f$. The well-posedness of \eqref{wavemap} has been studied in  \cite[Appendix A.3]{KLU1}, see also \cite[Appendix III]{Cb}. Actually, for $g'$ sufficiently close to $\widehat g$ in $C^m(M(T_1)), m\geq 5$, the wave map equation \eqref{wavemap} has a unique solution 
\beq
f\in C^0([0, T_1]; H^{m-1}(\mcm))\cap C^1([0, T_1]; H^{m-2}(\mcm)). 
\eeq
Here we take $T_1>T_0$ such that $M(T_0)\subset f(M(T_1))$ and $g = f_*g'$ is defined on $M(T_0)$. Let 
\beqq\label{eqspace}
E^{s}(M(T_0)) =  \bigcap_{p = 0}^{s} C^p([0, T_0]; H^{s-p}(\mcm)).
\eeqq
We also abbreviate the notation as $E^s$ when the manifold is clear. If $m$ is even, the wave map $f\in E^{m-1}$ and $f$ depends continuously on $g' \in C^m(M(T_1))$. In particular, for $s\geq 4$, we notice that $E^s \subset C^p([0, T_0]\times \mcm)$ for $0\leq p < s-2.$ Therefore $f\in C^{m-3}(M(T_0))$ if $m\geq 4$ is even. We remark that these regularity results may not be optimal but are sufficient for our analysis. 

Suppose $f$ is a wave map with respect to $(g', \widehat g)$ and set $g = f_* g'$, then the identity map $\id$ is a wave map with respect to $(g, \widehat g)$ and the wave map equation for $\id$ is equivalent to the harmonicity condition $
\Gamma^n = \widehat \Gamma^n, n = 0, 1, 2, 3,$ where  $\widehat \Gamma^n = g^{\alpha\beta}\widehat \Gamma^n_{\alpha\beta}$. Let $\widehat F^n = \Gamma^n - \widehat \Gamma^n$, the Ricci curvature tensor can be written as 
\beq
\ric_{pq}(g) = (\ric_{\widehat g}(g))_{pq} + \ha (g_{pn} \widehat \nabla_q \widehat F^n + g_{qn} \widehat \nabla_p \widehat F^n), \ \ p, q = 0, 1, 2, 3,
\eeq
where $\ric_{\widehat g}(g)$ is called the reduced Ricci curvature tensor. Therefore, in wave gauge where $\widehat F^n$ vanishes, the Ricci tensor is the reduced Ricci tensor. Using \eqref{richar} and the harmonicity condition, we obtain that
\beqq\label{ricre}
\begin{gathered}
(\ric_{\widehat g}(g))_{\mu\nu}  = -\ha g^{pq} \p_p\p_q g_{\mu\nu}  + Q_{\mu\nu}, \ \ \mu, \nu = 0, 1, 2, 3, \\
\text{ where } Q_{\mu\nu } = g^{ab}g_{ps}\Gamma^p_{\mu b}\Gamma^s_{\nu a} + \ha ( \p_a g_{\mu\nu} \widehat \Gamma^a + g_{\nu l} \Gamma^l_{ab} g^{aq}g^{bd}  \p_\mu g_{qd}   
+ g_{\mu l} \Gamma^l_{ab} g^{aq} g^{bd}  \p_\nu g_{qd}) \\
+ \ha (g_{\mu q} \p_\nu \widehat\Gamma^q  + g_{\nu q} \p_\mu \widehat \Gamma^q).
 \end{gathered}
\eeqq 
Suppose that $g'$ is a solution to the Einstein equation $\ric(g') = T'_{em} + \ha(J^\mu \phi_\mu)g'$ in $M(T_0)$ where $T'_{em}$ is defined using $g'$. Then $g = f^*g'$ satisfies the reduced Einstein equations
\beqq\label{ein2}
\begin{gathered}
\ric_{\widehat g}(g) =  T_{em} +\ha (J^\mu \phi_\mu) g
\end{gathered}
\eeqq
in wave gauge, where  $T_{em}$ is defined with respect to $g.$

Next we consider the Maxwell equations $\delta_g d\phi = J^\flat$ for electromagnetic potential $\phi$. We impose the Lorentz gauge condition $\delta_g \phi = 0$ to fix the gauge choice of $\phi$, see for example Section 10.2 of \cite[Chap. VI]{Cb}. Then we have $\square_g \phi = \delta_g d \phi + d\delta_g\phi  = J^\flat$. So we can write the Maxwell equation as
\beqq\label{maxeq}
\square_g \phi_\beta = -g^{\alpha\la} \p_\la\p_\alpha \phi_\beta + \widehat\Gamma^\mu \p_\mu \phi_\beta  = J_\beta, \ \ \beta = 0, 1, 2, 3,
\eeqq
where we have used the expression \eqref{eqbeltra} and the harmonicity condition $\Gamma^\mu = \widehat \Gamma^\mu$.

Now we formulate the reduced Einstein-Maxwell equations in wave and Lorentz gauge. We consider small perturbations near a background field. Suppose that $(\widehat g, \widehat \phi)$ is the background field where $\widehat \phi = 0$  and $\widehat g$ satisfies $\ric(\widehat g) = 0$. Let $\bar J \in C^\infty(\mbr; T\mcm)$ with $\supp(\bar J) \subset V$. In wave and Lorentz gauge, the Einstein-Maxwell equations \eqref{einmax} with unknowns $(g, \phi, J^0)$ are reduced to
\beqq\label{einmax1}
\begin{gathered}
\left\{\begin{array}{c}
\ric_{\widehat g}(g) = T_{em}(\phi) + \ha(J^\mu\phi_\mu) g\\[3pt]
\square_g \phi =  J^\flat \\
\div_g J = 0
%\p_t ((-\det g)^\ha J^0) + \sum_{i = 1}^3 \p_i((-\det g)^\ha \bar J^i) = 0
\end{array}\right. \text{ in } M(T_0),\\[2pt]
g = \widehat g, \ \ \phi  = 0, \ \ J^0 = 0, \text{ in } M(T_0)\backslash J_g^+(\supp(\bar J)).
\end{gathered}
\eeqq
Here we denote by $T_{em}(\phi)$ the stress-energy tensor associated with $F = d\phi$. Note that $T_{em}(\phi)$ is a polynomial involving $g_{ij}, g^{ij}$ and first derivatives of $\phi$. To establish the well-posedness for this system, it is convenient to consider the perturbed fields
\beqq\label{perfield}
\vec w \doteq (u, \phi) = (g, \phi) - (\widehat g, 0), 
\eeqq
where $\doteq$ means by definition. Using \eqref{divJ}, we can express $J^0$ in terms of $\bar J$ and $g$ i.e.\
\beqq\label{eqj0}
J^0 = F(x, g, \bar J, \p g, \p \bar J) \doteq  -(-\det g)^{-\ha} \sum_{i = 1}^3 \int_0^t \p_i((-\det g)^\ha  \bar J^i )ds.
\eeqq
In particular, $F$ is a smooth function of its arguments when $g$ is non-degenerate. Note that here we can replace the
vector field $\partial_t$ along which the divergence equation $\div_gJ=0$ is solved by any time-like vector field, but for clarity of notations we use just $\p_t$. Using \eqref{ricre} and the expression of $\square_g \phi$, we see that equations \eqref{einmax1} in terms of $\vec w$ can be reduced to 
\beqq\label{pernon}
\begin{gathered}
\left\{\begin{array}{c}
-g^{pq} \p_p \p_q u_{\mu\nu}   + \mcb^{\mu\nu}(x, \vec w, \p \vec w) =  G(x, \vec w, \bar J, \p \bar J),\\[3pt]
-g^{pq} \p_p \p_q \phi_{\beta}  + \mcb^\beta(x, \vec w, \p \vec w) = \sum_{\alpha = 1}^3g_{\beta\alpha}\bar J^\alpha + g_{\beta 0} F(x, g, \bar J, \p g, \p \bar J),
\end{array}\right. \text{ in } M(T_0),\\[3pt]
\vec w = 0, \text{ in } M(T_0)\backslash J_g^+(\supp(\bar J)),
\end{gathered}
\eeqq
where $\mu, \nu, \beta = 0, 1, 2, 3$, $\mcb^\bullet$ are smooth functions of $\vec w, \p \vec w$ and $G(\bullet)$ is a smooth function in its arguments. In particular, $G(\bullet)$ comes from the term $\ha(J^\mu\phi_\mu) g$.  We shall establish first the well-posedness of system \eqref{pernon} and then the coupled system \eqref{einmax1}. The system \eqref{pernon} is a quasilinear hyperbolic system of PDEs for $\vec w$. In Appendix B of \cite{KLU1}, the authors studied the local solvability and stability for certain nonlinear hyperbolic system of PDEs arising from the Einstein-scalar field equations. The proof follows similar fixed point iteration arguments as for the Cauchy problem, see for example \cite{HKM} and \cite[Appendix III]{Cb}. Although our nonlinear function $G(x, \vec w, \bar J, \p \bar J)$ is non-local due to \eqref{eqj0}, the proofs in \cite{KLU1} can be slightly modified to apply to  our system \eqref{pernon} since we consider local problem for $t$ finite. We briefly recall the results in \cite{KLU1}.

In \eqref{perfield}, $u$ is a symmetric $4\times 4$ matrix and $\phi$ is a $4$-vector. By renumbering, we can regard $\vec w = (u, \phi)$ as a $14$-vector. For integer $\kappa$, we let $\bfb^\kappa$ be a vector bundle on $M$ such that the fiber $\bfb_x^\kappa$ is a $\kappa$ dimensional vector space. In the following, we use section-valued Sobolev space $H^m(M; \bfb^\kappa)$. We denote
\beq
E^m(M(T_0); \bfb^\kappa) = \bigcap_{j = 0}^m C^j([0, T_0]; H^{m-j}(\mcm; \bfb^\kappa)), \ \ m\in \mbn.
\eeq
Also, we use the abbreviation $E^m$ and denote $E_0^m$ the compactly supported functions in $E^m.$  Now we state the local well-posedness for \eqref{pernon}, hence the reduced Einstein-Maxwell equations \eqref{einmax1}. As we already mentioned, the proof can be found in Appendix B of \cite{KLU1}. 
\begin{prop}\label{stabest}
Let $m_0\geq 4$ be an even integer and $T_0>0$.  If $\bar J$ is compactly supported and $\|\bar J\|_{E^{m_0}}< c_0$ is small enough, there exists a unique solution $\vec w$ satisfying the equation \eqref{pernon} on $M(T_0)$ and 
\beqq\label{eqstab}
\|\vec w\|_{E^{m_0}} \leq C \|\bar J\|_{E^{m_0+1}},
\eeqq
where $C$ denotes a generic constant.
\end{prop}
We remark that these regularity requirements are not optimal but sufficient for our work. Also, the space for $\bar J$ is $E^{m_0 + 1}$ because the functions $F$ and $G$ contains derivative of $\bar J$. We see that for $\|\bar J\|_{E^{m_0+1}}$ sufficiently small so that $g = \widehat g + u$ is non-degenerate,  we can solve $J^0$ from \eqref{eqj0}. Moreover, we have $J^0 \in E^{m_0}$ and $\| J^0\|_{E^{m_0}}\leq C \|\bar J\|_{E^{m_0+1}}$ for $m\geq 4$ even. This proves the well-posedness of the system \eqref{einmax1}.  

%%%%%%%%%%%%%%%%%%
\subsection{The full Einstein-Maxwell equations}\label{secfull}
We've found solutions to the reduced Einstein-Maxwell equations in wave and Lorentz gauge. To make sure that they give solutions to the full Einstein-Maxwell equations, we need to check that the wave gauge condition $\widehat\Gamma^n = \Gamma^n, n = 0, 1, 2, 3$ and the Lorentz gauge condition $\delta_g \phi = 0$ are satisfied. For the sourceless Einstein-Maxwell equations, this can be found in e.g.\ Lemma 10.1 of \cite[Section 10.3]{Cb}. The same idea applies here. The key is that the stress-energy tensor $T_{sour}$ satisfies the conservation law when $J$ does. 
 
\begin{lemma}\label{gasa}
If $(g, \phi)$ and $J$ satisfies \eqref{einmax1} in wave and Lorentz gauge, the harmonic gauge condition $\Gamma^n = \widehat \Gamma^n$ and Lorentz gauge condition $\delta_g \phi = 0$ are satisfied on $M(T_0)$. 
\end{lemma}
\bpf
First consider the function $\mcg = \delta_g\phi.$  
%It follows from \eqref{eqpoF} and \eqref{eqpoA} that 
%\beqq\label{eqlg}
%\nabla^\alpha F_{\alpha\beta} = \max(\phi) + \p_{\beta} \mcg.
%\eeqq
We know that 
\beqq\label{eqlg}
\square_g\phi = \delta_g F + d \mcg, \ \ F = d\phi.
\eeqq
Also, we know from \eqref{einmax1} that $\square_g \phi  = J^\flat$ and $J$ satisfies the conservation law for Maxwell equations which implies $\delta_g J^\flat = 0$. Now we take the codifferential of equation \eqref{eqlg} and use $\delta_g^2 = 0$ %and use the identity $\nabla^\beta\nabla^\alpha F_{\alpha\beta} = 0$ 
to obtain $\delta_g(d\mcg) = 0,$ which is  $\square_g \mcg = 0$ as $\mcg$ is a zero form (function). Therefore, $\mcg$ is a solution to the linear hyperbolic system
\beq
\begin{gathered}
\square_g \mcg = 0 \text{ in } M(T_0),\\[2pt]
 \mcg = 0 \text{ in } M(T_0)\backslash J_g^+(\supp(\bar J)).
\end{gathered}
\eeq
The system has a unique solution. Hence we conclude that $\mcg = 0$ i.e.\ $\delta_g \phi = 0$ in $M(T_0)$. 

Next, consider $\widehat F = (\widehat F^n)_{n = 0}^3$ with $\widehat F^n = \Gamma^n - \widehat \Gamma^n$. The proof of $\widehat F = 0$ follows from standard arguments. We briefly repeat the proof for completeness. If we write the reduced Einstein equation \eqref{einmax1} as 
\beq
\ein_{\widehat g}(g) \doteq \ric_{\widehat g} (g) - \ha \tr_g(\ric_{\widehat g}(g)) g = T_{sour},
\eeq
it follows from the Maxwell equation that $T_{sour}$ satisfies the conservation law $\div_g T_{sour} = 0$. From standard calculation for example \cite[Appendix A.4]{KLU1}, we know  that 
\beq
\ein_{jk}(g) - T_{jk} = \ha(g_{jn}\widehat \nabla_k \widehat F^n + g_{kn}\widehat\nabla_j \widehat F^n - g_{jk}\widehat\nabla_n\widehat F^n), \ \ j, k = 0, 1, 2, 3.
\eeq
Taking divergence and using Bianchi identity, we arrive at
\beq
0 = g^{nm}\nabla_n\widehat \nabla_m \widehat F^q + W^q(\widehat F), \ \ q = 0, 1, 2, 3, 
\eeq
where $W$ is a first order linear differential operator whose coefficients are polynomials of $\widehat g, g$ and their derivatives. Thus, $\widehat F$ satisfies a second order linear hyperbolic system
\beq
\begin{gathered}
g^{nm}\nabla_n \widehat \nabla_m \widehat F^q + W^q(\widehat F) = 0,\ \ \text{ in } M(T_0),\\[2pt]
\widehat F = 0, \text{ in } M(T_0)\backslash J_g^+(\supp(\bar J)).
\end{gathered}
\eeq
Since this system is uniquely solvable, we conclude that $\widehat F = 0$ i.e. the harmonic gauge condition is satisfied. This finishes the proof.
\epf

Finally, we discuss the well-posedness the full Einstein-Maxwell equations and we'd like to explain the regularity requirements in the set $\mcd(\delta)$. In wave gauge, we take $J = (J^0, \bar J)$ a one form with $J^0\in E^{m_0}, \bar J\in E^{m_0 + 1}, m_0\geq 4$ even and by Prop.\  \ref{stabest}, we know that the solution $g = \widehat g + u$ and $\phi$ are in $E^{m_0}$ as well. Since $E^{m_0}$ is an algebra, we can take the physical field $J$ as a vector field in $E^{m_0}$ in wave gauge.  For the wave map $f$,  we know that  $f\in C^{p-3}$ if $g'\in C^{p}$ where $g'$ is the solution to the full Einstein-Maxwell equations with electric current $J' = f^*J$. For solving the inverse problem, we'd like to take into account all $J$ in $E^{12}\subset C^8$ in wave gauge. Then $g\in C^{8}$ and $f\in C^4$. Thus we should include $J'\in C^4$ and $g' \in C^4$ in the data set.

%=============================================%
\section{Linearization of the Einstein-Maxwell equations}\label{einlin}
We've found the set  of $J$ for which the Einstein-Maxwell equations are well-posed. In the rest of the paper, we shall work with the following system in wave and Lorentz gauge
\beqq\label{einmax2}
\begin{gathered}
\left\{\begin{array}{c}
\ric_{\widehat g}(g) = T_{em}(\phi) + \ha(J^\mu\phi_\mu) g\\[3pt]
\square_g \phi =  J^\flat
\end{array}\right. \text{ in } M(T_0),\\[2pt]
g = \widehat g, \ \ \phi  = 0,  \text{ in } M(T_0)\backslash J_g^+(\supp(\bar J)).
\end{gathered}
\eeqq
The corresponding system for the perturbed fields $\vec w = (u, \phi) = (g, \phi) - (\widehat g, 0)$ are
\beqq\label{pernon1}
\begin{gathered}
\left\{\begin{array}{c}
-g^{pq} \p_p \p_q u_{\mu\nu}   + \mcb_{\mu\nu}(x, \vec w, \p \vec w) =  (J^\mu \phi_\mu)g,\\[3pt]
-g^{pq} \p_p \p_q \phi_{\beta}  + \mcb_\beta(x, \vec w, \p \vec w) = J_\beta,
\end{array}\right. \text{ in } M(T_0),\\[3pt]
\vec w = 0, \text{ in } M(T_0)\backslash J_g^+(\supp(\bar J)).
\end{gathered}
\eeqq
This is a second order quasilinear hyperbolic system for $\vec w$. We consider its linearization in this section and construct the so-called distorted plane wave solutions. The theory for linear hyperbolic PDEs, for example wave equations, is well-developed, see e.g.\ \cite{Fr, Ho3, Bar}. We review the causal inverse for the linearized Einstein equation as a paired Lagrangian distribution \cite{MU}. This is particularly convenient for the analysis of propagation and interaction of singularities.  

%%%%%%%%%%%
\subsection{The causal inverse}
Consider the equation \eqref{pernon1}  with source $J = \eps \mcj$ depending on a small parameter $\eps >0$.  For the moment, we can assume that $\mcj \in C^m(M; \bfb^4)$ is compactly supported. We will return to the issue of constructing $\mcj$ belonging to the data set in Section \ref{secmuls}. By the stability result Prop.\ \ref{stabest}, we can write the solution $(g, \phi) = (\widehat g, 0) + \eps (\dot g, \dot \phi)$ modulo terms of order $\eps^2$. Then $(\dot g, \dot \phi)$ satisfy the linearized Einstein-Maxwell equations (for $\mu, \nu, \beta = 0, 1, 2, 3$)
\beqq\label{perlin}
\begin{gathered}
\left\{\begin{array}{c}
 -\widehat g^{pq} \p_p \p_q \dot g_{\mu\nu} + \mca_{\mu\nu} (\dot g) = 0,\\[3pt]
-\widehat g^{pq} \p_p \p_q \dot \phi_{\beta}  + \mca_{\beta} (\dot \phi) = \widehat g_{\beta\alpha} \mcj^{\alpha},
\end{array}\right. \text{ in } M(T_0),\\[3pt]
(\dot g, \dot \phi) = 0, \text{ in } M(T_0)\backslash J_{\widehat g}^+(\supp(\mcj)),
\end{gathered}
\eeqq
where $\mca_\bullet$ are first order linear differential operators. The form of $\mca_\bullet$ follows from direct calculations using \eqref{ricre}, \eqref{stau} and \eqref{maxeq}. We emphasis that the system \eqref{perlin} is actually decoupled and we  have $\dot g = 0$. 
%\HOX{It seems that $\dot g=0$. Should we point this out?}
This facilitates our analysis and is  the reason why we consider the vacuum background.

It is convenient to write the system in matrix form. Let $\arv = (\dot g, \dot \phi)$ be section valued in $\bfb^{14}$. We let $\vec{\mcj} = (\bold{0}, (\widehat g_{\beta\alpha} \mcj^{\alpha})_{\beta=0}^3)$ where $\bold{0}$ is the zero section in $\bfb^{10}$. Then equation \eqref{perlin} can be written as
\beqq\label{lineq}
\bold{P}\arv \doteq (-\widehat g^{pq}  \p_p \p_q )\bold{Id} \arv + V(x, \p) \arv = \vec{\mcj},
\eeqq
where $\bold{Id}$ denotes the $14\times 14$ identity matrix and $V$ is a $14\times 14$ block-diagonalized matrix whose elements are first order differential operators.  From \eqref{perlin}, we notice that $\bold{P}$ is hyperbolic with principal term the wave operator $\square_{\widehat g}\bold{Id}$ for which we know there exists a causal inverse if $(M, \widehat g)$ is globally hyperbolic, see for example \cite{Bar}. Moreover, the Schwartz kernel of the causal inverse can be described as a paired Lagrangian distribution as shown in \cite{MU}, see also \cite{KLU1}. We review the construction in the following.

For two Lagrangians $\La_0, \La_1 \subset T^*X$ intersecting cleanly at a codimension $k$ submanifold i.e. 
\beq
T_q\La_0\cap T_q\La_1  = T_q(\La_0\cap \La_1),\ \ \forall q\in \La_0\cap \La_1,
\eeq
the paired Lagrangian distribution associated with $(\La_0, \La_1)$ is denoted by $I^{p, l}(\La_0, \La_1)$, see \cite{DUV, MU, GrU93, GrU90} for details. Let $\mcp(x, \xi) = |\xi|^2_{\widehat g^*}$ be the symbol of $-\widehat g^{pq}  \p_p \p_q$ with $\widehat g^* = \widehat g^{-1}$ the dual metric. Let $\Sigma_{\widehat g} = \{(x, \xi)\in T^*M: \mcp(x, \xi) = 0\}$ be the characteristic set which consists of light-like co-vectors. The Hamilton vector field of $\mcp$ is denoted by $H_\mcp$ and in local coordinates
\beq
H_\mcp = \sum_{i = 0}^3( \frac{\p \mcp}{\p \xi_i}\frac{\p }{\p x_i} - \frac{\p \mcp}{\p x_i}\frac{\p }{\p \xi_i}).
\eeq
The integral curves of $H_\mcp$ in $\Sigma_{\widehat g}$ are called null bicharacteristics and there projections to $M$ are geodesics. Let $\diag = \{(z, z')\in M\times M: z = z'\}$ be the diagonal and 
\beq
N^*\diag = \{(z, \zeta, z', \zeta')\in T^*(M\times M)\backslash 0: z = z', \zeta' = -\zeta\}
\eeq 
be the conormal bundle of $\diag$ minus the zero section. We let $\La_{\widehat g}$ be the Lagrangian obtained by flowing out $N^*\diag\cap \Sigma_{\widehat g}$ under $H_\mcp$. Here we regard $\Sigma_{\widehat g}, H_\mcp$ as objects on  product manifold $T^*M\times T^*M$ by lifting from the left factor. 
It is proved in Lemma 3.1 of \cite{KLU1} that for the linear differential operator $\bold{P}$, there exists a causal inverse $\bfq \in I^{-\frac{3}{2}, -\ha}(N^*\diag, \La_{\widehat g}; \bold{B}^{14})$ such that $\bold{P}\bfq = \bold{Id}$ on $\mce'(M; \bfb^{14})$.  Later, we also use $\bfq_{\widehat g}$ to emphasis the dependence on the metric. Moreover, from \cite[Prop.\ 5.6 ]{DUV} or \cite[Theorem 3.3]{GrU90}, we know that $\bfq: H_{\comp}^{m}(M; \bold{B}^{14})\rightarrow H^{m+1}_{\loc}(M; \bold{B}^{14})$ is continuous for $m\in \mbr$.

%%%%%%%%%%%
\subsection{Distorted plane waves}\label{sectdist}
We shall consider source $\mcj$ with singularities in the normal directions of a submanifold. For a submanifold $Y\subset X$ of codimension $k$, the conormal bundle $N^*Y$ is a Lagrangian submanifold. We review Lagrangian distributions in the scalar case, see \cite{Ho3, Ho4}. 
 
Let $X$ be a $n$ dimensional smooth manifold and $\La$ be a smooth conic Lagrangian submanifold of $T^*X\backslash 0$. Following the standard notation, we denote by $I^\mu(\La)$ the Lagrangian distribution of order $\mu$ associated with $\La$. In particular, for $U$ open in $X$, let $\phi(x, \xi): U\times \mbr^N \rightarrow \mbr$ be a smooth non-degenerate phase function that locally parametrizes $\La$ i.e.\ $\{(x, d_x\phi): x\in U, d_\xi \phi = 0\} \subset \La.$ Then $u\in I^{\mu}(\La)$ can be locally written as a finite sum of oscillatory integrals
\beq
\int_{\mbr^N} e^{i\phi(x, \xi)} a(x, \xi) d\xi, \ \ a\in S^{\mu + \fnf - \frac{N}{2}}(U\times \mbr^N),
\eeq
where $S^\bullet(\bullet)$ denotes the standard symbol class, see \cite[Section 18.1]{Ho3}.  For $u\in I^\mu(\La)$, we know that the wave front set $\WF(u)\subset \La$ and $u\in H^s(X)$ for any $s< -\mu-\fnf$. The principal symbol of $u$ is well-defined as a half-density bundle tensored with the Maslov bundle on $\La$, see \cite[Section 25.1]{Ho4}.

For a submanifold $Y\subset M$, the conormal distributions to $Y$ are denoted by $I^\mu(N^*Y)$. Occasionally, we also use the notation $I^{m}(Y) = I^{m + \frac{d}{2} - \frac n4}(N^*Y)$ with $d$ the codimension of $Y$ so that $m$ is the order of the symbols. Consider  local representations of such distributions. We can find local coordinates $x = (x', x''), x'\in \mbr^k, x''\in \mbr^{n-k}$ such that $Y = \{x' = 0\}$. Let $\xi = (\xi', \xi'')$ be the dual variable, then $N^*Y = \{x' = 0, \xi'' = 0\}$. We can write $u\in I^\mu(N^*Y)$ as
\beq
u = \int_{\mbr^k} e^{ix'\xi'} a(x'', \xi')d\xi', \ \ a\in S^{\mu + \fnf - \frac{k}{2}}(\mbr^{n-k}_{x''}; \mbr^k_{\xi'}).
\eeq
In this case, the principal symbol is 
\beq
\sigma(u) = (2\pi)^{\fnf - \frac{k}{2}} a_0(x'', \xi')|dx''|^\ha |d\xi'|^\ha,
\eeq
where  $a_0 \in S^{\mu + \fnf - \frac{k}{2}}(\mbr^{n-k}_{x''}; \mbr^k_{\xi'})$ is such that $a - a_0 \in  S^{\mu + \fnf - \frac{k}{2}- 1}(\mbr^{n-k}_{x''}; \mbr^k_{\xi'})$, see \cite[Section 18.2]{Ho3}. Later, we also use the notation $\sigma_{N^*Y}(u)$ to emphasis where the symbol is defined. For a Lorentzian manifold $(M, g)$, there is a natural choice of the density bundle, the volume element $d\vol_g$. Hence we can trivialize the distributional half densities and regard $u$ as distributions. {\em We emphasis that in our notation for principal symbols, we do not specify the order but refer to the distribution space for the orders.}

\begin{figure}[htbp]
\centering
% Generated with LaTeXDraw 2.0.1
% Sun May 08 15:43:20 PDT 2016
% \usepackage[usenames,dvipsnames]{pstricks}
% \usepackage{epsfig}
% \usepackage{pst-grad} % For gradients
% \usepackage{pst-plot} % For axes
\scalebox{.7} % Change this value to rescale the drawing.
{
\begin{pspicture}(0,-3.6725378)(10.202974,3.6525378)
\definecolor{color130b}{rgb}{0.8,0.8,0.8}
\rput{28.0817}(1.5960418,-3.7945466){\psellipse[linewidth=0.04,dimen=outer,fillstyle=solid,fillcolor=color130b](8.3845215,1.2937199)(0.69,1.8)}
\rput{28.0817}(0.6909437,-2.8008814){\psellipse[linewidth=0.04,dimen=outer,fillstyle=solid,fillcolor=color130b](5.945321,-0.019025208)(0.4,0.99)}
\psline[linewidth=0.04cm,dotsize=0.07055555cm 2.0,arrowsize=0.05291667cm 2.0,arrowlength=1.4,arrowinset=0.4]{*->}(0.31991574,-2.9864004)(10.175545,2.249302)
\psdots[dotsize=0.12,dotangle=28.0817](2.8879344,-1.6389287)
\psline[linewidth=0.04cm](2.8432288,-1.6401124)(8.34944,3.6325378)
\psline[linewidth=0.04cm](2.8255832,-1.649527)(10.182975,-0.058932148)
\psdots[dotsize=0.12,dotangle=28.0817](5.931199,0.0074430085)
\psdots[dotsize=0.12,dotangle=28.0817](8.436867,1.3443165)
\usefont{T1}{ptm}{m}{n}
\rput(0.4124707,-3.4498425){$x_0$}
\usefont{T1}{ptm}{m}{n}
\rput(3.3824706,-1.8098425){$x'$}
\usefont{T1}{ptm}{m}{n}
\rput(6.392471,-1.3898425){$Y$}
\psline[linewidth=0.04cm,arrowsize=0.133cm 2.0,arrowlength=1.4,arrowinset=0.4]{->}(2.8410156,-1.6548425)(4.2210155,-0.9148425)
\psline[linewidth=0.04cm,arrowsize=0.153cm 2.0,arrowlength=1.4,arrowinset=0.4]{->}(0.34101564,-2.9548426)(1.2010156,-2.4948425)
\usefont{T1}{ptm}{m}{n}
\rput(1.3824707,-2.9498425){$\theta_0$}
\usefont{T1}{ptm}{m}{n}
\rput(4.7524705,-1.0298425){$\theta'$}
\psline[linewidth=0.04cm,arrowsize=0.05291667cm 2.0,arrowlength=1.4,arrowinset=0.4]{->}(2.9010155,-1.5948424)(3.8810155,-0.6548425)
\psline[linewidth=0.04cm,arrowsize=0.05291667cm 2.0,arrowlength=1.4,arrowinset=0.4]{->}(2.8410156,-1.6548425)(4.341016,-1.3148425)
\end{pspicture} 
}
\label{figdpw}
\caption{Illustration of distorted plane waves.}
\end{figure}
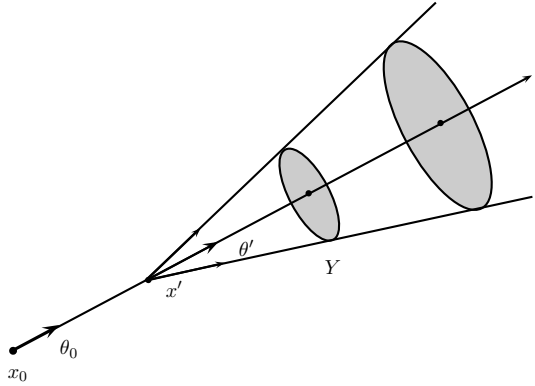

We follow \cite{KLU1} to construct distorted plane waves. For $(x_0, \theta_0)\in L^+M$, recall that $\gamma_{x_0, \theta_0}(t), t\geq 0$ is the geodesic from $x_0$ with direction $\theta_0$. For $s_0>0$ a small parameter, we let 
\beq
K(x_0, \theta_0; t_0, s_0) = \{\gamma_{x', \theta'}(t)\in M(T_0); \theta \in \mco(s_0), t\in (0, \infty)\},
\eeq
where $(x', \theta') = (\gamma_{x_0, \theta_0}(t_0), \dot \gamma_{x_0, \theta_0}(t_0))$ and $\mco(s_0) \subset L^+_{x'}M$ is an open neighborhood of $\theta'$ consisting of $\zeta\in L_{x'}^+M$ such that $\|\zeta - \theta'\|_{\widehat g^+}< s_0$. We observe that as $s_0\rightarrow 0$, $K(x_0, \theta_0; t_0, s_0)$ tends to the geodesic $\gamma_{x_0, \theta_0}$.  Next, let 
\beqq\label{ymani}
Y(x_0, \theta_0; t_0, s_0) = K(x_0, \theta_0; t_0, s_0)  \cap \{t = 2t_0\}
\eeqq
be a $2$-dimensional surface, which intersects the geodesic at $\gamma_{x_0, \theta_0}(2t_0)$, see Fig.\ \ref{figdpw}.  We let $\La(x_0, \theta_0; t_0, s_0)$ be the Lagrangian submanifold obtained from flowing out $N^*K(x_0, \theta_0; t_0, s_0)\cap N^*Y(x_0, \theta_0; t_0, s_0)$ under the Hamilton vector field of $\mcp$ in $\Sigma_{\widehat g}$. More precisely, for a smooth vector field $V$ on $M$ and $(x, \xi)\in TM$, we denote the integral curve of $V$ from $(x, \xi)$ by $\exp (tV)(x, \xi), t\in \mbr$. Then we have 
\beqq\label{lamani}
 \La(x_0, \theta_0; t_0, s_0) = \bigcup_{t\geq 0} \exp \big( t H_\mcp\big) (N^*K(x_0, \theta_0; t_0, s_0)\cap N^*Y(x_0, \theta_0; t_0, s_0) \cap \Sigma_{\widehat g}).
\eeqq
It is convenient to introduce a notation for the flow out under $\La_{\widehat g}$. For any $\Gamma\subset T^*M$, we denote the flow out of $\Gamma$ by $\Gamma^{\widehat g} = \La'_{\widehat g} \circ (\Gamma\cap \Sigma_{\widehat g})$, where as usual in microlocal analysis
\beq
\La_{\widehat g}' = \{(x, \xi, y, \eta)\in T^*M\times T^*M : (x, \xi, y, -\eta)\in \La_{\widehat g}\},
\eeq
and $\circ$ denotes the composition of sets as relations. 
The proposition below is essentially proved in \cite{KLU1} and is stated for the current setting. 
\begin{prop}\label{distor}
Suppose $Y(x_0, \theta_0; t_0, s_0)$ and $\La(x_0, \theta_0; t_0, s_0)$ are defined as in \eqref{ymani} and \eqref{lamani} respectively. For $\vec{\mcj}\in I^{\mu+1}(N^*Y; \bold{B}^{14})$ with $\mu$ an integer, the solution to the linearized Einstein-Maxwell equations $\bold{P}\arv =\vec{\mcj}$ is 
\beq
\arv = \bfq(\vec{\mcj}) \in I^{\mu - \ha}(\La; \bold{B}^{14}) \text{ in } M(T_0)\backslash Y,
\eeq
and $\arv$ is called distorted plane waves. Moreover, for $(p, \xi)\in N^*Y\cap \Sigma_{\widehat g}$ and $(q, \eta)\in \La(x_0, \theta_0; t_0, s_0)$ which lie on the same bicharactersitics, the principal symbol of $\arv$ is given by
\beq
\sigma(\arv)(q, \eta) = \sigma(\bfq)(q, \eta; p, \xi) \sigma(\vec{\mcj})(p, \xi),
\eeq
where $\sigma(\bfq)$ is a $14\times 14$ invertible matrix.
\end{prop} 

We remark that in general $\arv$ is a well-defined Lagrangian distribution and before the first conjugate point, it is a conormal distribution. As discussed in \cite{KLU}, it is complicated to analyze the interaction of singularities past the conjugate points and this difficulty will be overcome using another argument. Let ${\tau_0}>0$ be such that $\gamma_{x_0, \theta_0}({\tau_0})$ is the first conjugate point of $x_0$ along $\gamma_{x_0, \theta_0}$. Then the exponential map $\exp_{x_0}$ is a local diffeomorphism from a neighborhood of $t\theta_0\in T_{x_0}M$ to a neighborhood of $\gamma_{x_0, \theta_0}(t)$ for $t< {\tau_0}$. Therefore, $K(x_0, \theta_0; t_0, s_0)$ is a co-dimension $1$ submanifold near $\gamma_{x_0, \theta_0}(t)$ and 
\beq
{\La (x_0, \theta_0; t_0, s_0) = N^*K(x_0, \theta_0; t_0, s_0) }
\quad \text{ near } \gamma_{x_0, \theta_0}(t) \text{ for }  t< {\tau_0}.
\eeq
In particular, before the first conjugate point of $x_0$ along $\gamma_{x_0, \theta_0}$, $\vec v$ is a conormal distribution to $K(x_0, \theta_0; t_0, s_0)$.  We mention that  in \cite{KLU}, the null cut points was considered, see \cite[Section 2.1]{KLU}. On a globally hyperbolic Lorentzian manifold, the first null cut point $\widehat x$ of $x_0$ along the geodesic $\gamma_{x_0, \theta_0}$ is either the first conjugate point or there are at least two light-like geodesics joining $x_0$ and $\widehat x$. So the first cut point appears on or before the first conjugate point.

%%%%%%%%%%%%%%%%%%%%%%%
\subsection{Microlocal linearization conditions}\label{secmuls}
Now we address the issue that for the source $J$ to be in the data set, it must satisfy the conservation law. Suppose $J = J_\eps$ depending on some small parameter $\eps$ and we let $\mcj = \p_\eps J_\eps|_{\eps = 0}$ be its linearization. From the conservation law $\div_g J = 0$, we derive its linearization  
\beqq\label{eqlincons}
\p_i ((-\det \widehat g)^\ha \mcj^i) = 0. 
\eeqq
Notice that this is one equation and the Einstein summation is over $i = 0, 1, 2, 3.$ For $\mcj\in I^{\mu}(N^*Y; \bfb^4)$ where $Y$ is defined in \eqref{ymani}, the principal symbol at $(x, \xi)\in N^*Y$ should satisfy the microlocal linearized conservation law
\beqq\label{eqmlincons}
 \xi_i\sigma(\mcj^i)(x, \xi) = 0.
\eeqq
We shall denote such conormal distributions by $I^\mu_C(N^*Y; \bfb^4)$. We emphasis that the space $\mcx_{x, \xi}$ of principal symbols of $\vec \mcj = (\bold{0}, \mcj)$ at $(x, \xi)$ satisfying \eqref{eqmlincons} is a $3$-dimensional subspace of a $14$-dimensional vector space. 

Since we derived the linearized Einstein-Maxwell equations from the reduced equations in wave and Lorentz gauge, the linearized solutions must also satisfy the linearized gauge conditions. Their microlocalization imposes conditions on the symbols of the distorted plane waves. We start with the harmonicity condition $g^{nm}\Gamma^j_{nm} = g^{nm}\widehat \Gamma^j_{jm}, j = 0, 1, 2, 3.$ By linearization we obtain 
\beq
-\widehat g^{an}\p_a \dot g_{nj} + \ha \widehat g^{pq} \p_j \dot g_{pq} = m_{j}^{pq} \dot g_{pg}, \ \ j = 0, 1, 2, 3, 
\eeq
where $m_j$ depends on $\widehat g$ and its derivative, see \cite[Sec.\ 3.2.3]{KLU1}. When the background metric $\widehat g$ is Minkowski, we have $m_j = 0$. The Lorentz gauge condition is $g^{\alpha \la}\p_{\la}\phi_{\alpha} = \widehat \Gamma^\mu \phi_{\mu}$. The linearization gives $
\widehat g^{\alpha\la}\p_{\la}\dot \phi_{\alpha} = \widehat \Gamma^{\mu} \dot \phi_{\mu}. $ Also, when $\widehat g$ is Minkowski, the right hand side vanishes. Let $\dot g, \dot \phi$ be conormal distributions in $I^{\mu}(\La)$ where $\La$ is defined in \eqref{lamani}. For any $(y, \eta)\in \La$, we obtain that the principal symbols of $\dot g, \dot \phi$ should satisfy the microlocal linearized gauge conditions i.e.\
\beqq
\label{lingauge1}
& &
-\widehat g^{an}\xi_a \sigma(\dot g_{nj}) + \ha \widehat g^{pq} \xi_j \sigma(\dot g_{pq}) = 0, \ \ j = 0, 1, 2, 3, \\
\label{lingauge2}
& &\widehat g^{\alpha\la}\xi_{\la}\sigma(\dot \phi_{\alpha}) =0. 
\eeqq
We shall denote conormal distributions satisfying these conditions by $I^{\mu}_G(\La)$. In particular, for our consideration, the linearized metric component $\dot g$ vanishes thus the microlocal gauge condition is reduced to \eqref{lingauge2}. Then we notice that the space $\mcy_{y, \eta}$ of principal symbols of $(\dot g, \dot \phi)$ at $(y, \eta)$ is also a $3$-dimensional vector space. Now if we consider the linear map $R=\sigma(\bfq)(y, \eta, x, \xi)$ in Proposition \ref{distor}, we see that 
\beqq\label{eqbij}
\text{$R$ is bijective from $\mcx_{x, \xi}$ to $\mcy_{y, \eta}$.} 
\eeqq

Now we make a connection with the so called microlocal linearization stability condition, introduced in \cite[Assumption $\mu$-LS]{KLU1} for the Einstein-scalar field equations. The condition plays an important role because it allows one to construct distorted plane wave solutions by prescribing the principal symbols. Let's consider the setup of the Einstein-Maxwell equations \eqref{einmax} in the introduction, see also Theorem \ref{main1}. Let $Y$ be defined as in \eqref{ymani} and $Y\subset V$.  For any $(x, \xi)\in N^*Y\cap L^*M$ and vector $A = (A^i)_{i = 0}^3\in \mbr^4$ satisfying $\xi_i A^i = 0$, one can find $\bar \mcj^i \in I^{\mu+1}(N^*Y), i = 1, 2, 3$ compactly supported in $V$ such that $\sigma(\bar \mcj^i)(x, \xi) = A^i, i = 1, 2, 3.$ Then by the well-posedness of the Einstein-Maxwell equations with sources in wave and Lorentz gauge, we know that there exists a family of solutions $(g_\eps, \phi_\eps, J^0_\eps)$ for $\bar J_\eps = \eps \bar \mcj, \eps > 0$ such that $J_\eps = (J^0_\eps, \bar J_\eps)$ satisfies the conservation law. Let $\mcj^0$ be the linearization of $J^0_\eps$. From \eqref{eqj0}, we obtain that
\beqq\label{eqmcj0}
\mcj^0 = -(-\det \widehat g)^{-\ha} \sum_{i = 1}^3 \int_0^t \p_i((-\det \widehat g)^\ha   \bar\mcj^i )ds.
\eeqq 

\begin{lemma}\label{lmconormal}
Let $Y$ be defined as in \eqref{ymani} and suppose that $\bar \mcj^i \in I^{\mu+1}(N^*Y), i = 1, 2, 3$. Then we have  $\mcj^0\in I^{\mu+ 2}(N^*Y)$. Moreover, near $L^*M\cap N^*Y$, we have $\mcj^0\in I^{\mu+1}(N^*Y)$. 
\end{lemma}
\bpf
We recall the equivalent definition of conormal distributions e.g.\ \cite[Def.\ 18.2.6]{Ho3}. So $u\in I^\mu(N^*Y)$ if and only if $V_1\cdots V_k u\in {}^\infty H^{-\mu-\frac n4}_{\loc}(M), \forall k \geq0$ where the $V_j$ are smooth vector fields on $M$ tangent to $Y$ and ${}^\infty H^s_{\loc}$ denotes Besov spaces, see e.g.\ \cite[Appendix B]{Ho3}. By our construction, we know that $Y \subset \{t = 2t_0\}$ so that $TY\subset T\mcm$. In particular, the vector fields tangent to $Y$ are spanned by $\frac{\p}{\p y^i}, i = 1, 2, 3$ with $(y^i)$ the local coordinates for $\mcm$. We observe that 
\beq
V_1\cdots V_k \mcj^0 = - \sum_{i = 1}^3 V_1\cdots V_k \int_0^t  (-\det \widehat g)^{-\ha}\p_i((-\det \widehat g)^\ha   \bar\mcj^i )ds= \sum_{i = 1}^N \int_0^t f_i(s, y) ds,
\eeq
where $f_i\in {}^\infty H^{-\mu + 2 -\frac n4}_{\loc}(M)$. Using for example \cite[Corollary B.1.6]{Ho3}, we conclude that $V_1\cdots V_k \mcj^0 \in {}^\infty H^{-\mu+ 2-\frac n4}_{\loc}(M)$ so that $\mcj^0\in I^{\mu+ 2}(N^*Y)$. Now we use that the principal symbols of $J$ should satisfy the microlocal linearized conservation law. So for $(x, \xi)\in N^*Y$ and $\xi$ close to $L_x^*M$ such that $\xi_0\neq 0$, we have $\xi_0 \sigma(\mcj^0)(x, \xi) = 0$ which implies $\sigma(\mcj^0) = 0$. Thus the principal symbol vanishes and $\mcj^0\in I^{\mu+1}(N^*Y)$ microlocally near $L^*M\cap N^*Y$. 
\epf
Notice that $\mcj^0$ and $\bar\mcj$ satisfy the linearized conservation law, therefore, their principal symbols satisfy the microlocal linearized conservation law. Continue our argument, we have $\sigma(\mcj^0)(x, \xi) = A^0$. To conclude, given the symbols $A$ at $(x, \xi)\in N^*Y\cap L^*M$, we can find a one parameter family of sources $J_\eps$ such that $\mcj = \p_\eps J_\eps|_{\eps = 0}\in I^\mu(N^*Y)$ with $\sigma(\mcj)(x, \xi) = A$. Moreover, for $-\mu$ large enough, there exists a unique solution $(g_\eps, \phi_\eps)\in C^4(M)$ to the Einstein-Maxwell equations \eqref{einmax} with the source $J_\eps$. Notice that we can choose $\bar J_\eps = \eps \bar \mcj$  compactly supported in $V$ and from the formula \eqref{eqj0} of $J^0_\eps$, we can shrink the support of $\bar \mcj$ so that $J^0_\eps$ is supported in $V$ (but not compactly supported). This is actually the corresponding version of the microlocal linearization stability condition for the Einstein-Maxwell equations (compare with \cite[Assumption $\mu$-LS]{KLU1}). Finally, combining with the observation \eqref{eqbij}, we conclude that for any given symbols of $(\dot g, \dot \phi)$ in $\mcy_{y, \eta}$, one can find the corresponding source $J_\eps$.

%==================================%
\section{Linearized electromagnetic and gravitational waves }\label{gwave} %\HOX{The name of the section is changed.}
In this section, we demonstrate that due to the nonlinearity and coupling of the Einstein-Maxwell equations, it is possible to generate gravitational waves using electromagnetic sources. For simplicity, {we consider in this section the linearization of Einstein-Maxwell equations, in wave gauge, when the background manifold is Minkowski.} 

Let  $(M, \widehat g)$ be the Minkowski space-time $(\mbr^4, h)$ with $h = -(dx^0)^2 + \sum_{i = 1}^3 (dx^i)^2.$ In this case, the Einstein-Maxwell equations are simple because $\widehat \Gamma_{ij}^k$ and the derivatives of $\widehat g_{ij}$ all vanish. 
We start with the expression of the nonlinear equations \eqref{pernon1}. Recall that $(u, \phi) = (g - \widehat g, \phi - 0)$ is the perturbed fields. 
From \eqref{ricre}, we get the reduced Ricci tensor ($\mu, \nu = 0, 1, 2, 3$)
\beq 
\begin{gathered}
(\ric_{\widehat g}(g))_{\mu\nu}  = - \ha g^{pq} \p_p \p_q u_{\mu\nu}  + g^{ab}g_{ps}\Gamma^p_{\mu b}\Gamma^s_{\nu a} + \ha (g_{\nu l} \Gamma^l_{ab} g^{aq}g^{bd}  \p_\mu u_{qd} + g_{\mu l} \Gamma^l_{ab} g^{aq} g^{bd}  \p_\nu u_{qd}).
\end{gathered}
\eeq 
The stress-energy tensor for the electromagnetic field is ($\alpha, \beta = 0, 1, 2, 3$)
\beq
\begin{gathered}
T_{em, \alpha\beta}= g^{\la\mu}F_{\alpha\mu} F_{\beta\la} - \frac{1}{4}g_{\alpha\beta}g^{\la\gamma}g^{\mu\delta}F_{\gamma\delta}F_{\la\mu}, \ \ F_{\alpha\beta} = \p_{\alpha}\phi_\beta - \p_{\beta} \phi_{\alpha}.
\end{gathered}
\eeq
Let $J= (J^0, \bar J)$ be the source where  $J^0$ is defined in \eqref{eqj0}. The reduced Einstein equation is of the form
\beq
- g^{pq} \p_p \p_q u_{\mu\nu} + 2P_{\mu\nu} - 2T_{em, \mu\nu} = (J^i\phi_i) g, \ \ \mu, \nu = 0, 1, 2, 3,
\eeq
where $P_{\mu\nu}$ is the semilinear term in $\ric_{\widehat g}(g)$ and $J^0$ is regarded as a nonlinear function of $g, J^i, i = 1, 2, 3$. Finally, the reduced Maxwell equations \eqref{maxeq} is
\beq
 -g^{\alpha\la} \p_\la\p_\alpha \phi_\beta   = g_{\beta\alpha}J^\alpha, \ \ \beta = 0, 1, 2, 3.
%g^{\alpha\la} \p_\la\p_\alpha \phi_\beta + \p_\beta(g^{\alpha\la})\p_{\la} \phi_{\alpha} - g^{\alpha\la}\Gamma^\mu_{\la\beta} (\p_{\alpha} \phi_\mu - \p_\mu \phi_{\alpha}) = -\mcj_\beta.
\eeq
In the above equations, $g^{\alpha\beta}$ can be computed as following
\beq
\begin{gathered}
g^{-1} = (h + u)^{-1} = (\Id + h^{-1}u + (h^{-1}u)^2 + (h^{-1}u)^3 + \cdots) h^{-1},
\end{gathered}
\eeq 
where the $2$-tensors are treated as matrices in local coordinates. The components are
\beqq\label{metexpan}
\begin{split}
g^{ab} &= h^{ab} + h^{aa'}u_{a'b'}h^{b'b} + h^{aa'}u_{a'c'} h^{c'c} u_{cb'}h^{b'b} + \cdots\\
& = h^{ab} + \sum_{a, b = 0}^3 h^{aa}h^{bb} u_{ab} + \sum_{a, b, c = 0}^3 h^{aa}h^{bb}h^{cc}u_{ac}u_{cb} + \cdots.
\end{split}
\eeqq
It is worth mentioning that the first line of the above formula is in Einstein summation but the second line is not and we made simplifications using properties of $h$. Using these formulas, we can find the nonlinear terms $\mcb$ in \eqref{pernon1} explicitly. Next, let $\bar J = \eps \bar \mcj$ and $\mcj^0$ be the linearization of $J^0$ defined in \eqref{eqmcj0}. The linearized equations \eqref{perlin} are simply
\beq
\begin{gathered}
 -h^{pq}  \p_p \p_q \dot g_{\mu\nu} = 0, \ \ \mu, \nu = 0, 1, 2, 3 \\[2pt]
 -h^{\alpha\la} \p_\la\p_\alpha \dot \phi_\beta = h_{\beta\alpha}\mcj^\alpha, \ \ \beta = 0, 1, 2, 3.
 \end{gathered}
\eeq 
We denote the causal inverse $Q_h = (-h^{pq}  \p_p \p_q)^{-1}$. The main result of this section is
\begin{theorem}\label{maingw}
Let $(M, h)$ be the Minkowski space-time with $M = \mbr^4$ and $M(T_0) = (-\infty,T_0)\times \mbr^3, T_0>0$. Let $Y$ be a 2-dimensional surface on $\{t = 0\}$ and $\mcy$ be the null hyper-surface from $Y$ i.e.\ $\mcy = \{\exp_h (t V) \in M(T_0):  t > 0, V\in L^+_YM\}$. Let $\bar \mcj^i \in I^{\mu + 1}(N^*Y), i = 1, 2, 3, \mu \leq -10$ be compactly supported. Consider the Einstein-Maxwell equations in wave gauge
\beqq\label{einmaxmin}
\begin{gathered}
\left\{\begin{array}{c}
\ein(g) =  T_{sour} \\[3pt]
\delta_g d\phi =  \eps \mcj^\flat\\
\div_g \mcj = 0
\end{array}\right. \text{ in } M(T_0),\\[2pt]
g = h, \ \  \phi = 0, \ \ \mcj^0 = 0,  \text{ in } M(T_0)\backslash J_g^+(\supp (\bar \mcj)).
\end{gathered}
\eeqq
For $\eps>0$ sufficiently small, the solution to \eqref{einmaxmin} on $M\backslash Y$  satisfies
\beq
\begin{gathered}
g = h + \eps^2 g_1 + o(\eps^2),\\
\phi = \eps \phi_1 + o(\eps^2),
\end{gathered}
\eeq
where the term in $o(\eps^2)$ is small in $H^4(M)$, such that on $M\backslash Y$,  
\begin{enumerate}
\item $g_1\in H^9(M), \phi_1\in H^8(M)$ with $\singsupp(g_1),  \singsupp(\phi_1)\subset \mcy$.
\item $g_1, \phi_1$ are non-vanishing if  $\bar\mcj$ is non-vanishing.
\end{enumerate}
\end{theorem}

\bpf
According to  the discussion in Section \ref{secfull}, for $\bar \mcj\in H^7(M(T_0))\subset C^4(M(T_0)) \subset E^4$ and $\eps$ small, we can find source $J_\eps = (J^0_\eps, \eps \bar \mcj)$ such that there exists a unique solution $\vec w \in E^4$ to \eqref{pernon1} in wave gauge with source term $J_\eps$. Let $\mcj^0 = \p_\eps J^0_\eps|_{\eps = 0}$ and we write $\vec\mcj = (\bold{0}, \eps \mcj)$ with $\bold{0}\in \bfb^{10}$. The solution to the linearized Einstein-Maxwell equations \eqref{perlin} is $\vec v = (\dot g, \dot \phi) = \bfq(\vec\mcj) \in I^{\mu-\ha}(\La_1\backslash Y; \bold{B}^{14}) \subset H^{8}(M; \bold{B}^{14})$, where  $\La_1$ is the flow out of $N^*Y\cap L^*M$ under $\La^h$. Notice that here $\vec v$ is a Lagrangian distribution but may not be conormal. The projection of $\La_1\backslash Y$ to $M$ is exactly $\mcy.$ Also, from \eqref{perlin}, we know that the  metric components $\dot g$ in $I^{\mu-\ha}(\La_1)$ of $\vec v$ are all zero. From the equations \eqref{pernon1} and \eqref{perlin}, we find that 
\beqq\label{expan2}
\begin{gathered}
\bold{P}(\vec w - \eps \vec v) + P_2(x, \vec w) + H_2(x, \vec w)  + o(\eps^2) = \eps^2 (\mcj^i \dot \phi_i) \vec h,
\end{gathered}
\eeqq
where $o(\eps^2)$ is in $H^4(M)$. We explain the terms. The term $\vec h = (h, \bold{0}), \bold{0}\in \bfb^4$ comes from the terms in the Einstein-Maxwell equations with $\mcj$ and we used the fact that $\dot g  = 0$. The term $P_2$ can be regarded as a  smooth $14\times 14$ matrix obtained from the expansion of $(h + u)^{-1}$:
\beq
P_2(x, \vec w) = (h u h)^{pq} \frac{\p^2 }{\p x^p \p x^q} \vec w,
\eeq
where $h$ and $u$ are treated as matrices.  The term $H_2$ is section valued in $\mathbf{B}^{14}$ and comes from the nonlinear terms of \eqref{pernon1}. More explicitly, the elements can be written as
\beq
\begin{gathered}
H_{2, \gamma} = \sum_{i, j = 1}^{14} \sum_{ \alpha, \beta = 1}^4 \sum_{a, b=0}^1 H^\gamma_{2, ij\alpha\beta ab} \p_\alpha^a w_i \p_\beta^b w_j,
\end{gathered}
\eeq
where the coefficients are all smooth and $\gamma = 1, 2, \cdots, 14.$ From \eqref{expan2}, we obtain 
\beq
\vec w = \eps \vec v - \bfq(P_2(x, \vec w) + H_2(x, \vec w)) + \eps^2 \bfq((\mcj^i \dot \phi_i) \vec h)+ o(\eps^2).
\eeq
Substitute $\vec w$ to the right hand side, we get 
\beqq\label{eqexpan1}
\vec w = \eps \vec v - \eps^2 \bfq(P_2(x, \vec v) + H_2(x, \vec v)) + \eps^2 \bfq((\mcj^i \dot \phi_i) \vec h) + o(\eps^2).
\eeqq
Notice that since $H^4(M)$ is an algebra and $\bfq$ is continuous from $H_{\comp}^4(M)$ to $H_{\loc}^{5}(M)$, the remainder term $o(\eps^2)$ is still in $H^4(M)$. Since $\dot g = 0$, we have $P_2(x, \vec v) = 0$. The terms in $H_2$ with metric components vanish as well. So it suffices to consider terms in $H_2$ which only have the electric potentials. We shall denote them by $\widehat H_2$. Observe that these terms can only come from the stress-energy tensor $T_{em}(\phi)$. Actually, we find using \eqref{stau} and \eqref{metexpan} that
\beqq\label{eqh2}
\begin{gathered}
\widehat H_{2, \alpha\beta}(x, \vec w) =  -2(h^{aa'}F_{\alpha a}F_{\beta a'} - \frac{1}{4} h_{\alpha\beta}h^{aa'}h^{bb'}F_{ab}F_{a'b'}), \ \ \alpha, \beta = 0, 1, 2, 3, \\ 
\widehat H_{2, \mu}(x, \vec w) = 0, \ \ \mu = 0, 1, 2, 3,
\end{gathered}
\eeqq
where the $F_\bullet$ are defined in terms of $\dot \phi$. Here $\widehat H_{2, \alpha\beta}$ are the terms of $\widehat H_2$ in the reduced Einstein equations and $\widehat H_{2, \mu}$ are the terms from the Maxwell equations. Also, we renumbered the section valued $\widehat H_2$. Now we can simplify \eqref{eqexpan1} to
\beq
\vec w = \eps \vec v - \eps^2\bfq (\widehat H_2(x, \dot \phi))+  \eps^2 \bfq((\mcj^i \dot \phi_i) \vec h) + o(\eps^2).
\eeq
In particular, the metric components are
\beq
u =   \eps^2g_1 + o(\eps^2), \ \ g_{1, \alpha\beta} = -Q_h(\widehat H_{2, \alpha\beta} + (\mcj^i \dot \phi_i)  h_{\alpha\beta}).
\eeq
Using the fact that $H^8(M)$ is an algebra and the continuity of $Q_h$, we see that $g_1\in H^9(M)$. The singular support properties follows from standard wave front analysis, see e.g.\ \cite[Section 1.3]{Du}.

Finally, we observe from \eqref{eqh2} that if $\alpha = \beta = 0$, we have that
\beq
\begin{split}
\widehat H_{2, 00} & = -2(h^{aa'}F_{0 a}F_{0 a'} - \frac{1}{4} h_{00}h^{aa'}h^{bb'}F_{ab}F_{a'b'}) 
= -2(\sum_{a = 0}^3  h^{aa}F_{0a}^2 + \sum_{a, b = 0}^3  \frac{1}{4}h^{aa}h^{bb}F_{ab}^2) \\
&= -2(\sum_{a = 1}^3 F_{0a}^2 - \ha \sum_{a = 1}^3F_{0a}^2 + \frac{1}{4} \sum_{a, b = 1}^3F^2_{ab}) 
 = -\sum_{a = 1}^3 F_{0a}^2  - \frac{1}{2} \sum_{a, b = 1}^3F^2_{ab}.
\end{split}
\eeq
Here we used the fact that $F_{ab}$ is antisymmetric and $F_{00}=0$. It is clear that $\widehat H_{2, 00}$ vanishes if and only if $F_{ab} = 0, a, b = 0, 1,2, 3$. It follows from the Maxwell equation \eqref{maxeq} and its linearization that $F = 0$ implies $\mcj = 0$. Thus if $\bar\mcj$ is non-vanishing, we conclude that $\dot \phi$ is non-vanishing from the linearized equations and hence $g_1$ is non-vanishing outside the support of $\mcj$. This ends the proof of the theorem.
\epf

%==================================%
\section{Interactions of distorted plane waves}\label{singu}
From this section, we return to the setting of a general globally hyperbolic vacuum spacetime $(M,\widehat g)$ instead of the Minkowski spacetime in Section 4. When four distorted plane waves meet at a point, new singularities could be produced. In \cite{KLU1},  the nature of these singularities  are analyzed using Gaussian beam solutions and stationary phase type arguments. Recently in \cite{LUW}, the authors studied such interactions carefully for scalar waves using paired Lagragnians and symbol calculus. We shall apply these results to the Einstein-Maxwell equations. 

Let $x^{(j)} \in V$ and $(x^{(j)}, \theta^{(j)})\in L^+M, j = 1, 2, 3, 4$ be such that 
\beq
\gamma_{x^{(j)}, \theta^{(j)}}([0, t_0])\subset V, \ \ x^{(j)}(t_0)\notin J_{\widehat g}^+(x^{(k)}(t_0)), \ \ j \neq k,
\eeq
which means that the points are causally independent. We define $K_j = K(x^{(j)}, \theta^{(j)}; t_0, s_0), j = 1, 2, 3, 4$ and $\La_j(x^{(j)}, \theta^{(j)}; t_0, s_0)$ similar to \eqref{ymani} and \eqref{lamani}.  Let ${\tau_j}, j = 1, 2, 3, 4$ be such that $\gamma_{x^{(j)}, \theta^{(j)}}({\tau_j})$ is the first conjugate point of $x^{(j)}$ along the geodesics and ${\tau_{\textrm{min}}} = \min_{j = 1, 2, 3, 4}({\tau_j})$. {\em In the rest of this section, we shall study the interactions only in the following set}  
\beq
\begin{gathered}
\mathcal{N}((\vec x, \vec \theta), t_0) = M(T_0)\backslash \bigcup_{j = 1}^4 J_{\widehat g}^+(\gamma_{x^{(j)}, \theta^{(j)}}({\tau_j})),  
\end{gathered}
\eeq
where $\vec x = (x^{(1)}, x^{(2)}, x^{(3)}, x^{(4)}), \vec \theta = (\theta^{(1)}, \theta^{(2)}, \theta^{(3)}, \theta^{(4)}),$ i.e.\ away from the causal future of the conjugate points.  In particular, $\La_j = N^*K_j$  in $\mathcal{N}((\vec x, \vec \theta), t_0). $

Recall that two submanifolds $X, Y$ of $M$ intersect transversally if 
\beq
T_q X + T_q Y = T_q M, \ \ \forall q\in X\cap Y.
\eeq
For the codimension $1$ submanifolds $K_i,  i = 1, 2, 3, 4$ we consider, $\La_i = N^*K_i \subset L^*M$ i.e.\ the co-vectors normal to $K_i$ are light-like. For the moment, we assume that they only intersect at $q_0$  transversally meaning 
\begin{enumerate}
\item $K_i, K_j$ intersect transversally at a codimension $2$ submanifold $K_{ij}, i< j$;
\item $K_i, K_j,  K_k$ intersect  at  a codimension $3$ submanifold $K_{ijk}, i<j<k$;
\item $K_i, i = 1, 2, 3, 4$ intersect at a point $q_0$. 
\end{enumerate}
In particular, (2) implies that $K_{ij}\cap K_k$ transversally and (3) implies that the four submanifolds intersect at a point $q_0$ and the normal co-vectors $\zeta_i$ to $K_i$ at $q_0$ are linearly independent so they span the cotangent space $T_{q_0}^*M$. We remark that for any $q\in M$, we can find $K_i$ intersect transversally at $q$.
For $i = 1, 2, 3, 4$, we shall denote 
\beq
\begin{gathered}
\La_i = N^*K_i; \ \ \La_{ij} = N^*K_{ij}, i< j; \ \ \La_{ijk} = N^*K_{ijk}, i<j<k; \ \ \La_{q_0} = T^*_{q_0}M\backslash 0.
\end{gathered}
\eeq
These are Lagrangian submanifolds of $T^*M$. We introduce the following notations
\beq
\begin{gathered}
\La^{(1)} = \bigcup_{i = 1}^4 \La_i; \ \ \La^{(3)} = \bigcup_{i, j, k = 1, i<j<k}^4\La_{ijk}.
%K^{(1)} = \cup_{i = 1}^4 K_i; \ \ K^{(3)} = \cup_{i, j, k = 1, i<j<k }^4 K_{ijk}.
\end{gathered}
\eeq
It is obvious that $\La^{(1), \widehat g} = \La^{(1)}$, but $\La^{(3), \widehat g}$ is not the same as $\La^{(3)}$. We denote $\Theta = \La^{(1)}\cup \La^{(3), \widehat g}$ which is of particular importance  below. 

We begin with  the interaction of scalar valued distorted plane waves. Let $v_i \in I^{\mu}(N^*K_i), i = 1,2 , 3, 4$  and $Q_{\widehat g} \in I^{-\frac{3}{2}, -\ha}(N^*\diag, \La_{\widehat g})$ be the causal inverse of $\square_{\widehat g}$ on $(M, \widehat g)$. For semilinear wave equations studied in \cite{LUW}, singularities of the following terms were analyzed
\beqq\label{fourtho}
\begin{gathered}
\mcy_1 = Q_{\widehat g}(cv_1v_2v_3v_4),\\[2pt] 
\mcy_2 = Q_{\widehat g}(av_1Q_{\widehat g}(bv_2v_3v_4)),\ \ \mcy_3 = Q_{\widehat g}(bv_1v_2Q_{\widehat g}(av_3v_4)),\\[2pt] 
\mcy_4 = Q_{\widehat g}(av_1Q_{\widehat g}(av_2Q_{\widehat g}(av_3v_4))), \ \ \mcy_5 = Q_{\widehat g}(aQ_{\widehat g}(av_1v_2)Q_{\widehat g}(av_3v_4)),
\end{gathered}
\eeqq
where $a, b, c$ are smooth functions on $M$. The terms $\mcy_i, i = 1, 2, 3, 4, 5$ involve multiplication of four conormal distributions whose singular support intersect at $q_0$. It is proved in Prop.\ 3.9 of \cite{LUW} that $\mcy_i$ have conormal singularities at $\La^{\widehat g}_{q_0}\backslash \Theta$ i.e.\ the flow out of $\La_{q_0}$ away from $\Theta$. We will see that such terms also appear in the asymptotic analysis of the Einstein-Maxwell equations, but they may contain derivatives. So we slightly generalize the result in \cite{LUW} to include cases when the order of $v_i$ are different. 
\begin{prop}\label{porder}
Let $v_i\in I^{\mu_i}(\La_i), i = 1, 2, 3, 4$ and $\widetilde \mu = \sum_{i = 1}^4 \mu_i$. Let $\mcv$ be a smooth vector field. We have the following conclusions
\begin{enumerate}\label{internew}
\item $Q_{\widehat g}(cv_1v_2v_3v_4)  \in I^{\widetilde \mu+\frac{3}{2}}(\La^{\widehat g}_{q_0}\backslash \Theta); $ 
 \item $Q_{\widehat g}(av_1\mcv Q_{\widehat g}(bv_2v_3v_4)), Q_{\widehat g}(bv_1v_2 \mcv Q_{\widehat g}(av_3v_4))  \in I^{\widetilde\mu+ \ha}(\La^{\widehat g}_{q_0}\backslash \Theta); $
 \item $Q_{\widehat g}(av_1 \mcv Q_{\widehat g}(av_2 \mcv Q_{\widehat g}(av_3v_4))), Q_{\widehat g}(a\mcv Q_{\widehat g}(av_1v_2)\mcv Q_{\widehat g}(av_3v_4))  \in I^{\widetilde \mu-\frac{1}{2}}(\La^{\widehat g}_{q_0}\backslash \Theta).$
\end{enumerate}
\end{prop}
\bpf
For a smooth vector field $\mcv$, $\mcv v_i \in I^{\mu+1}(N^*Y_i)$ if $v_i \in I^{\mu}(N^*Y_i), i = 1, 2, 3, 4$. This is addressed for $\mcv = \nabla_g$ in \cite[Lemma 4.1]{WZ} and the general case is the same. Also, that lemma tells that $\mcv Q_{\widehat g} \in I^{-\frac{3}{2}+ 1, -\ha}(N^*\diag, \La_{\widehat g})$ and the principal symbols can be found. Then the proof is that of \cite[Prop.\ 3.9]{LUW} by adjusting the orders.
\epf

We emphasis that the set $\Theta$ is the union of the wave front set of $v_i$ and $\La^{(3), \widehat g}$ is the wave front set of singularities generated by the triple wave interactions. So we only look at the new singularities in $\mcu^{(4)}$ produced by the four wave interactions.  
 
In Section 3.5 of \cite{LUW}, the principal symbols of the terms \eqref{fourtho} are found explicitly. For our purpose, we just need the symbol of $\mcy_3, \mcy_4$ and $\mcy_5$. Consider the symbols at $(q, \eta) \in \La^{\widehat g}_{q_0}\backslash \Theta$, which is joined with $(q_0, \zeta) \in \La_{q_0}$ by bi-characteristics. We can write $\zeta = \sum_{i = 1}^4\zeta_i$ where $\zeta_i\in N^*_{q_0}K_i$. Let $A_i$ be the principal symbols of $v_i$. Then we have
 \begin{multline}\label{eqy3y5}
 \begin{aligned}
 %&\sigma_{\La^{\widehat g}_{q_0}}(\mcy_1)(q, \eta) = (2\pi)^{-3}\sigma_{\La_{\widehat g}}(Q_{\widehat g})(q, \eta, q_0, \zeta)c(q_0)\prod_{i = 1}^4A_i(q_0, \zeta_i),\\
 &\sigma_{\La^{\widehat g}_{q_0}}(\mcy_3)(q, \eta) = (2\pi)^{-3}\sigma_{\La_{\widehat g}}(Q_{\widehat g})(q, \eta, q_0, \zeta)a(q_0) b(q_0) \frac{1}{|\zeta_3+\zeta_4|^2_{{\widehat g}^*(q_0)}}\cdot \prod_{i = 1}^4A_i(q_0, \zeta_i), \\
%&\sigma_{\La^{\widehat g}_{q_0}}(\mcy_2)(q, \eta) = (2\pi)^{-3}\sigma_{\La_{\widehat g}}(Q_{\widehat g})(q, \eta, q_0, \zeta)a(q_0) b(q_0) |\zeta_2+\zeta_3+\zeta_4|_{{\widehat g}^*(q_0)}^{-2} \cdot  \prod_{i = 1}^4A_i(q_0, \zeta_i),
&\sigma_{\La^{\widehat g}_{q_0}}(\mcy_4)(q, \eta) = (2\pi)^{-3}\sigma_{\La_{\widehat g}}(Q_{\widehat g})(q, \eta, q_0, \zeta)a^3(q_0)   \frac{1}{|\zeta_2+\zeta_3+\zeta_4|_{{\widehat g}^*(q_0)}^{2}}  \cdot \frac{1}{|\zeta_3+\zeta_4|_{{\widehat g}^*(q_0)}^{2}} \prod_{i = 1}^4A_i(q_0, \zeta_i),\\
&\sigma_{\La^{\widehat g}_{q_0}}(\mcy_5)(q, \eta) = (2\pi)^{-3}\sigma_{\La_{\widehat g}}(Q_{\widehat g})(q, \eta, q_0, \zeta)a^3(q_0) \frac{1}{|\zeta_3+\zeta_4|^2_{{\widehat g}^*(q_0)}}  \cdot  \frac{1}{|\zeta_1+\zeta_2|^2_{{\widehat g}^*(q_0)}} \prod_{i = 1}^4A_i(q_0, \zeta_i).
\end{aligned}
 \end{multline}
We remark that here we trivialized the density factors in the distributions in a local coordinate near $q_0$, so the symbols are functions. The generalization of the above formulas to the case with derivatives is quite straightforward because $\sigma(\mcv v_i) = \sigma(\mcv)\sigma(v_i)$, see also \cite[Lemma 4.1]{WZ}.\\

%%%%%%%%%%%
%\subsection{Singularities in the asymptotic expansion}\label{secasym}
Now we are ready to study the singularities for the Einstein-Maxwell equations. As in Section \ref{secmuls}, we assume that $\mcj^{(i)} \in I^{\mu+1}(N^*Y_i; \bfb^4), i = 1, 2, 3, 4$ are supported in $V$ and $\bar \mcj^{(i)}$ are compactly supported. In order that $\mcj^{(i)} \in E^{4}$, we shall take $\mu \leq -10$ so that $I^{\mu+1}(N^*Y_i)\subset H^{7}(M(T_0))\subset C^{4}(M(T_0))\subset E^{4}$. We denote $\vec \mcj^{(i)} = (\bold{0}, \mcj^{(i)})$ and let $\vec v^{(i)} = \bfq (\vec \mcj^{(i)}) \in I^{\mu-\ha}(N^*K_i; \bold{B}^{14})$ be distorted plane waves. Let $\eps_i, i = 1, 2, 3, 4$ be small parameters and $\vec \mcj = \sum_{i = 1}^4 \eps_i \vec \mcj^{(i)}.$ Then $\vec v = \sum_{i = 1}^4 \eps_i \vec v^{(i)} = (\dot g, \dot \phi)$ is the solution to the linearized equation \eqref{perlin} $\bold{P} \vec v = \vec \mcj$. We recall that $\dot g \in I^{\mu-\ha}(\La_i; \bold{B}^{10})$ is the %\HOX{Check the claim "$\vec w = (u, \phi)$ be the solution of the nonlinear equation \eqref{pernon1} with source $\vec \mcj$, where $u\in I^{\mu-\ha}(\La_i; \bold{B}^{10})$".} 
linearized metric component and $\dot \phi\in I^{\mu-\ha}(\La_i; \bold{B}^{4})$ is the linearized electromagnetic component. 

By the microlocal linearization condition, we let $J_\eps$ be  such that $\p_{\eps_i} J_\eps|_{\eps_i = 0} = \vec \mcj^{(i)}, i = 1, 2, 3, 4$. Let $\vec w = (u, \phi)$ be the solution of the nonlinear equation \eqref{pernon1} with source $J_\eps.$ From Prop.\ \ref{stabest}, we can write the asymptotic expansion of $\vec w$ as $\eps_i\rightarrow 0$
\beqq\label{uexpan}
\vec w = \vec v + \sum_{1\leq i< j \leq 4}\eps_i\eps_j \mcu^{(2)} + \sum_{1\leq i< j< k\leq 4}\eps_i\eps_j\eps_k \mcu^{(3)} + \eps_1\eps_2\eps_3\eps_4 \mcu^{(4)} + \mcr_\eps,
\eeqq
where $\mcr_\eps$ denotes the collection of terms in $H^4(M(T_0))$ and  $\bigcup_{i = 1}^4 O(\eps_i^2)$. In particular,
\beqq\label{eqmcu}
\mcu^{(4)} = \p_{\eps_1}\p_{\eps_2}\p_{\eps_3}\p_{\eps_4} \vec w|_{\{\eps_1=\eps_2=\eps_3=\eps_4=0\}}.
\eeqq
Our goal is to analyze the singularities in $\mcu^{(4)}$. It is not easy to find  the terms in \eqref{uexpan} explicitly because the nonlinear Einstein-Maxwell equations is complicated. We will simplify the computation by identifying the most singular terms in $\mcu^{(4)}$. 

 From equations \eqref{pernon1} and \eqref{perlin} which $\vec w, \vec v$ satisfy, we derive
\beqq\label{eqexpan}
\begin{gathered}
\bold{P}(\vec w - \vec  v) + \sum_{i = 2}^4 P_i(x, \vec w) + \sum_{i = 2}^4H_i(x, \vec w) + \mcr_\eps = 0.
\end{gathered}
\eeqq
We explain the terms in \eqref{eqexpan} in detail. For the moment, we shall ignore the asymptotic expansion terms  involving $\bar \mcj^{(i)}$ and we shall see later that they do not matter for our analysis.

 First of all, the $P_i, i = 2, 3, 4$ terms come from the quasilinear term $-g^{pq} \p_p\p_q \vec w_i , i = 1, 2, \cdots, 14$ in \eqref{pernon1}. In a given local coordinates, we can express two tensors as  $4\times 4$ matrices. It is easy to see that 
 \beqq
 \begin{split}
 g^{-1} &= (\widehat g + u)^{-1} = (\Id + \widehat g^{-1} u)^{-1} \widehat g^{-1}\\
 & = \widehat g^{-1} - \widehat g^{-1} u \widehat g^{-1} + (\widehat g^{-1} u)^2 \widehat g^{-1} + (\widehat g^{-1} u)^3 \widehat g^{-1} + \cdots.
 \end{split}
 \eeqq
Then the $P_i$ terms are
 \beq
 \begin{gathered}
 P_2(x, \vec w) = (\widehat g^{-1} u \widehat g^{-1})^{pq} \frac{\p^2 \vec w}{\p x^p \p x^q}, \ \ P_3(x, \vec w) = -(\widehat g^{-1} u \widehat g^{-1} u \widehat g^{-1})^{pq} \frac{\p^2 \vec w}{\p x^p \p x^q}, \\
 P_4(x, \vec w) = (\widehat g^{-1} u\widehat g^{-1} u \widehat g^{-1} u \widehat g^{-1})^{pq} \frac{\p^2 \vec w}{\p x^p \p x^q}
\end{gathered}
 \eeq
More precisely, we  write down the components of these terms as 
\beqq\label{eqformp}
\begin{split}
P_{2, i}(x, \vec w) &= \widehat g^{pa}u_{ab}\widehat g^{bq} \frac{\p^2 \vec w_{i}}{\p x^p \p x^q} \\
P_{3, i}(x, \vec w) & = -\widehat g^{pa} u_{ab} \widehat g^{bc} u_{cd} \widehat g^{dq} \frac{\p^2 \vec w_{i}}{\p x^p \p x^q} \\
P_{4, i}(x, \vec w) & = \widehat g^{pa} u_{ab} \widehat g^{bc} u_{cd} \widehat g^{de} u_{ef} \widehat g^{fq}\frac{\p^2 \vec w_{i}}{\p x^p \p x^q}, \ \  i = 1, 2, \cdots, 14.
\end{split}
\eeqq
Notice that each term has two derivatives and the coefficients of the derivatives are polynomials of the metric component. It is also convenient to regard them as multi-linear functions. For example, 
\beq
P_{2, i}(x, \vec w^{(1)}, \vec w^{(2)}) = \widehat g^{pa}u^{(1)}_{ab}\widehat g^{bq} \dfrac{\p^2 \vec w^{(2)}_{i}}{\p x^p \p x^q}, \ \ i = 1, 2, \cdots, 14.
\eeq

Next, the terms $H_i, i = 2, 3, 4$ in \eqref{eqexpan} come from the semilinear terms of \eqref{pernon1}. They are section valued in $\bfb^{14}$. Each component of $H_i$ is a sum of $i$-th order monomials of $\vec w_\bullet, \p \vec w_\bullet$ but at most quadratic in $\p \vec w_\bullet$, where $\bullet= 1, 2, \cdots, 14$ denote a generic index. This follows from the expression of the reduced Ricci tensor \eqref{ricre} and the fact that the Christoffel symbols only have first derivatives of the metric $g$. We can write the component of $H_i, i = 2, 3, 4$ as multilinear functions as
\beq
\begin{split}
H_{2, \theta}(x, \vec w^{(1)}, \vec w^{(2)}) &= \sum_{i, j = 1}^{14} \sum_{ \alpha, \beta = 1}^4 \sum_{a, b=0}^1 H^\theta_{2, ij\alpha\beta ab} \p_\alpha^a  \vec w^{(1)}_i \p_\beta^b  \vec w^{(2)}_j, \\
H_{3, \theta}(x, \vec w^{(1)}, \vec w^{(2)}, \vec w^{(3)}) &= \sum_{i, j,k = 1}^{14} \sum_{ \alpha, \beta, \gamma = 1}^4 \sum_{a, b, c =0, 1; a+b+c \leq 2} H^\theta_{3, ijk\alpha\beta\gamma abc} \p_\alpha^a \vec w^{(1)}_i \p_\beta^b \vec w^{(2)}_j \p_\theta^c \vec w^{(3)}_k,\\
H_{4, \theta}(x, \vec w^{(1)}, \vec w^{(2)}, \vec w^{(3)}, \vec w^{(4)}) &= \sum_{i, j,k,l  = 1}^{14}\sum_{\alpha, \beta, \gamma, \delta = 1}^4 \sum_{(a, b, c, d) \in \mca} H^\theta_{4, ijkl\alpha\beta\gamma\delta abcd} \p_\alpha^a \vec w^{(1)}_i \p_\beta^b \vec w^{(2)}_j \p_\gamma^c \vec w^{(3)}_k \p_\delta^d \vec w^{(4)}_l,
\end{split}
\eeq 
where the set $\mca = \{(a, b, c, d) : a, b, c, d =0, 1; a+b+c+d \leq 2\}$, $\theta = 1, 2, \cdots, 14$ and the coefficients are all smooth. We emphasis that the derivatives of $\vec w$ in these terms appear at most twice and such terms are especially important for the analysis below. We denote such terms i.e.\ terms in $H_i$ with two derivatives, by $\widehat H_i, i = 2, 3, 4.$   

For convenience, we denote
\beq
\begin{gathered}
G_i(x, \vec w) = P_i(x, \vec w) + H_{i}(x, \vec w), \ \ i = 2, 3, 4,
\end{gathered}
\eeq
then $G_i$ are polynomials in $\vec w, \p \vec w$ of order $i$. We also denote $\widehat G_i = P_i + \widehat H_i.$ Then the asymptotic expansion \eqref{eqexpan} can be written as
\beqq\label{eqite1}
\begin{gathered}
 \vec w =\vec  v - \bfq( G_2 + G_3 + G_4) + \mcr_\eps.
\end{gathered}
\eeqq
By the stability estimate, we have $\|\vec w\|_{E^{4}} \leq C\sum_{i = 1}^4 \eps_i$. Note $E^{4}\subset H^4(M)$ is an algebra. We have 
\beq
\| \vec w_i \vec w_j\|_{H^{4}(M)} \leq C(\sum_{i = 1}^4 \eps_i)^2.
\eeq
Therefore, by the continuity of $\bfq$, we obtained the first term in \eqref{uexpan}. Also, we notice that all the terms in $\mcr_\eps$ are in $H^{4}(M)$. To get other terms, we will iterate the formula \eqref{eqite1} i.e.\ plug \eqref{eqite1} to the right hand side of \eqref{eqite1}. This will generate many terms. However, we only need terms of the order $\eps_1\eps_2\eps_3\eps_4$ which can only be obtained from the multiplication of four terms of $\vec v$ and $\p \vec v$.

\begin{prop}\label{intert}
Consider the fourth order interaction term $\mcu^{(4)}$ defined in \eqref{eqmcu}.
\begin{enumerate}
\item If $\bigcap_{j = 1}^4\gamma_{x_j, \theta_j}(t) = q_0$, for $s_0$ sufficiently small so that $K_i$ only intersect at $q_0$,  then on $\mathcal{N}((\vec x, \vec \theta), t_0)$, we can write $\mcu^{(4)} = \bfq(\mch  + \widehat \mch)$ such that 
 \beq
 \bfq(\mch) \in I^{4\mu + \frac{3}{2}}(\La^{\widehat g}_{q_0}\backslash \Theta; \bold{B}^{14}) \text{ and } \bfq(\widehat \mch) \in I^{4\mu+ \ha}(\La^{\widehat g}_{q_0}\backslash \Theta; \bold{B}^{14}).
 \eeq 
The term $\mch$ is given by
\begin{multline}\label{symbh1}
\begin{gathered}
\mch =  \sum_{(i, j, k, l) \in \sigma(4)} \bigg( \widehat H_3(x, \vec v^{(i)}, \vec v^{(j)}, \bfq (\widehat H_2(x, \vec v^{(k)}, \vec v^{(l)}))) \\
+ \widehat H_3(x, \vec v^{(i)}, \bfq (\widehat H_2(x, \vec v^{(j)}, \vec v^{(k)})), \vec v^{(j)}) + \widehat H_3(x,  \bfq (\widehat H_2(x, \vec v^{(i)}, \vec v^{(j)})), \vec v^{(k)}, \vec v^{(l)}) \bigg)\\
 -\sum_{(i, j, k, l) \in \sigma(4)} \widehat G_2(x, \bfq(\widehat H_2(x, \vec v^{(i)}, \vec v^{(j)})), \bfq(\widehat H_2(x, \vec v^{(k)}, \vec v^{(l)}))) \\
 -\sum_{(i, j, k, l) \in \sigma(4)} \bigg( \widehat H_2(x, \vec v^{(i)}, \bfq(P_2(x, \bfq(\widehat H_2(x, \vec v^{(j)}, \vec v^{(k)})), \vec v^{(l)}))) \\
  + \widehat H_2(x, \bfq(P_2(x, \bfq(\widehat H_2(x, \vec v^{(i)}, \vec v^{(j)})), \vec v^{(k)})), \vec v^{(l)})\bigg)
\end{gathered}
\end{multline}
with $\sigma(4)$ denoting the set of permutations of $(1, 2, 3, 4)$.

\item If $\bigcap_{j = 1}^4\gamma_{x_j, \theta_j}(t) = \emptyset$,  then $\mcu^{(4)}$  is smooth on $\mathcal{N}((\vec x, \vec \theta), t_0) $ microlocally away from $\Theta$.
\end{enumerate}
\end{prop} 
\bpf

(1) We divide the proof into two steps. 

\textbf{Step 1:} We first determine the most singular terms in $\mcu^{(4)}$. We claim that $\mcu^{(4)} = \bfq(\mch + \widehat \mch)$ such that
 \beq
 \bfq(\mch) \in I^{4\mu + \frac{3}{2}}(\La^{\widehat g}_{q_0}\backslash \Theta, \bold{B}^{14}), \ \  \bfq(\widehat \mch) \in I^{4\mu+ \ha}(\La^{\widehat g}_{q_0}\backslash \Theta, \bold{B}^{14}).
 \eeq
Here $\mch = \sum_{i = 1}^3 \mch_{i}$ in which 
\begin{equation}\label{eqmch}
\begin{gathered}
\mch_{1} = -\sum_{(i, j, k, l) \in \sigma(4)} \widehat G_4(x, \vec v^{(i)}, \vec v^{(j)}, \vec v^{(k)}, \vec v^{(l)})\\ 
\mch_{2} = \sum_{(i, j, k, l) \in \sigma(4)} \bigg( \widehat G_3(x, \vec v^{(i)}, \vec v^{(j)}, \bfq (\widehat G_2(x, \vec v^{(k)}, \vec v^{(l)}))) + \widehat G_3(x, \vec v^{(i)}, \bfq (\widehat G_2(x, \vec v^{(j)}, \vec v^{(k)})), \vec v^{(j)})\\
+ \widehat G_3(x,  \bfq (\widehat G_2(x, \vec v^{(i)}, \vec v^{(j)})), \vec v^{(k)}, \vec v^{(l)})\bigg) \\
  +  \sum_{(i, j, k, l) \in \sigma(4)} \bigg(\widehat G_2(x, \bfq(\widehat G_3(x, \vec v^{(i)}, \vec v^{(j)}, \vec v^{(k)})), \vec v^{(l)}) + \widehat G_2(x, \vec v^{(i)}, \bfq(\widehat G_3(x, \vec v^{(j)}, \vec v^{(k)}, \vec v^{(l)}))) \bigg) ,
\end{gathered}
\end{equation}
and
\begin{equation}
\begin{gathered}\label{eqmch2}
 \mch_{3}  = -  \sum_{(i, j, k, l) \in \sigma(4)} \widehat G_2(x, \bfq(\widehat G_2(x, \vec v^{(i)}, \vec v^{(j)})), \bfq(\widehat G_2(x, \vec v^{(k)}, \vec v^{(l)})))\\
 -\sum_{(i, j, k, l) \in \sigma(4)} \bigg(\widehat G_2(x, \vec v^{(i)}, \bfq(\widehat G_2(x, \vec v^{(j)}, \bfq(\widehat G_2(x, \vec v^{(k)}, \vec v^{(l)}))))) \\
+ \widehat G_2(x, \vec v^{(i)}, \bfq(\widehat G_2(x, \bfq(\widehat G_2(x, \vec v^{(j)}, \vec v^{(k)})), \vec v^{(l)}))) \\
  + \widehat G_2(x, \bfq(\widehat G_2(x, \vec v^{(i)}, \bfq(\widehat G_2(x, \vec v^{(j)}, \vec v^{(k)})))), \vec v^{(l)})\\
  + \widehat G_2(x, \bfq(\widehat G_2(x, \bfq(\widehat G_2(x, \vec v^{(i)}, \vec v^{(j)})), \vec v^{(k)})), \vec v^{(l)})\bigg)
\end{gathered}
\end{equation}

%\beqq\label{eqmch}
%\begin{split}
%\mch_{1} &= \sum_{(i, j, k, l) \in \sigma(4)} \widehat H_4(x, \vec v^{(i)}, \vec v^{(j)}, \vec v^{(k)}, \vec v^{(l)}) + P_2(x, \vec v^{(i)}, \vec v^{(j)})\widehat H_2(x, \vec v^{(k)}, \vec v^{(l)}),\\
%\mch_{2} &= \sum_{(i, j, k, l) \in \sigma(4)} [\widehat H_3(x, \vec v^{(i)}, \vec v^{(j)}, \bfq(\widehat H_2(x, \vec v^{(k)}, \vec v^{(l)}))) + \widehat H_2(x, \vec v^{(i)}, \bfq(\widehat H_3(x, \vec v^{(j)}, \vec v^{(k)}, \vec v^{(l)})))],  \\
%\mch_{3} &= \sum_{(i, j, k, l) \in \sigma(4)} \widehat H_2(x, \bfq(\widehat H_2(x, \vec v^{(i)}, \vec v^{(j)})), \bfq(\widehat H_2(x, \vec v^{(k)}, \vec v^{(l)}))) \\
%&+ \sum_{(i, j, k, l) \in \sigma(4)} \widehat H_2(x, \vec v^{(i)}, \bfq(\widehat H_2(x, \vec v^{(j)}, \bfq(\widehat H_2(x, \vec v^{(k)}, \vec v^{(l)})))))
%\end{split}
%\eeqq

We use the formula \eqref{eqite1} to iterate to get the asymptotic terms. First of all, we put \eqref{eqite1} to the right hand side of \eqref{eqite1} to get one term from $G_4(x, \vec w)$: 
\beq
\begin{split}
&G_4(x, \vec w) = G_4(x, \vec v) + \mcr_\eps = \sum_{(i, j, k, l) \in \sigma(4)} G_4(x, \vec v^{(i)}, \vec v^{(j)}, \vec v^{(k)}, \vec v^{(l)}) + \mcr_\eps.
\end{split}
\eeq
The summation terms consist of two types of terms because $G_4 = P_4 + H_4$. The terms from $P_4$ can be found in \eqref{eqformp} and they all have two derivatives. Using Prop.\ \ref{porder}, we conclude that $\bfq (P_4(x, \vec v^{(i)}, \vec v^{(j)}, \vec v^{(k)}, \vec v^{(l)})) \in I^{4\mu + \frac{3}{2}}(\La^{\widehat g}_{q_0}\backslash \Theta)$. The terms from $H_4$ are of the form
\beq
A_{ijkl\alpha\beta m n}\vec v^{(a)}_i \vec v^{(b)}_j \p^m_\alpha \vec v^{(c)}_k \p_\beta^n \vec v^{(d)}_l, \ \ m, n\leq 1; \ \ i, j, k, l = 1, \cdots, 14; \ \ \alpha, \beta = 1, 2, 3, 4,
\eeq
and $a, b, c, d$ are permutations of $1,2,3,4$ (Note this is not in Einstein summation.) We can apply Prop.\ \ref{porder} to conclude that when $m, n = 1$, the term after applying $\bfq$ is in $I^{4\mu + \frac{3}{2}}(\La^{\widehat g}_{q_0}\backslash \Theta)$ and otherwise in $I^{4\mu + \ha}(\La^{\widehat g}_{q_0}\backslash \Theta)$. When $m = n = 1$, the terms only come from $\widehat H_4$. Thus we obtain the leading term $\mch_{1}$.

Next consider the other terms in the asymptotic expansion. To get order $\eps_1\eps_2\eps_3\eps_4$ terms, we need to iterate twice or three times using \eqref{eqite1}. From the term $G_3(x, \vec w)$, we get 
\beqq\label{eqh3}
\begin{split}
 G_3(x, \vec w) & = -\sum_{(i, j, k, l) \in \sigma(4)} \bigg( G_3(x, \vec v^{(i)}, \vec v^{(j)}, \bfq (G_2(x, \vec v^{(k)}, \vec v^{(l)}))) \\
+& G_3(x, \vec v^{(i)}, \bfq (G_2(x, \vec v^{(j)}, \vec v^{(k)})), \vec v^{(j)}) + G_3(x,  \bfq (G_2(x, \vec v^{(i)}, \vec v^{(j)})), \vec v^{(k)}, \vec v^{(l)})\bigg) + \mcr_\eps.
\end{split}
\eeqq
%\beqq
%\begin{split}
%&= \sum_{i, j,k = 1}^{14}\sum_{\alpha, \beta, \gamma = 1}^4 \sum_{a, b, c =0}^1 H_{3, ijk\alpha\beta\gamma abc} \p_\alpha^a \vec v_i \p_\beta^b \vec v_j \cdot \\
%&\p_\gamma^c \bfq_k(G_2(x, \vec v)) + \mcr_\eps, \ \ a+b+c \leq 2.
%\end{split}
%\eeqq
From Prop.\ \ref{porder}, we know that the term after applying $\bfq$ is in $I^{4\mu+\frac{3}{2}}(\La^{\widehat g}_{q_0}\backslash \Theta)$ if the $G_i, i = 2, 3$ terms involved have two derivatives i.e. they are $\widehat G_i, i = 2, 3.$ Otherwise, the terms are in $I^{4\mu +\ha}(\La^{\widehat g}_{q_0}\backslash \Theta)$ which are less singular. So we get one piece in $\mch_{2}$.  

Finally, from the term $G_2(x, \vec w),$ we get from the iteration using \eqref{eqite1} that 
\begin{equation}\label{eqg2}
\begin{split}
 &G_2(x, \vec w)  \\[2pt] 
 = & G_2(x, \vec v - \bfq(G_2(x, \vec w) + G_3(x, \vec w))) + \mcr_\eps \\[2pt] 
= & G_2(x, \bfq(G_2(x, \vec v)), \bfq(G_2(x, \vec v))) -  G_2(x, \vec v, \bfq(G_2(x, \vec w) + G_3(x, \vec v))) \\
&  - G_2(x, \bfq(G_2(x, \vec w) + G_3(x, \vec v)), \vec v)  + \mcr_\eps\\[2pt] 
 =  & G_2(x, \bfq(G_2(x, \vec v)), \bfq(G_2(x, \vec v))) -  G_2(x, \vec v, \bfq(G_3(x, \vec v)))   - G_2(x, \bfq(G_3(x, \vec v)), \vec v)  \\
 &  -  G_2(x, \vec v, \bfq(G_2(x, \vec v -  \bfq(G_2(x, \vec v)))))  - G_2(x, \bfq(G_2(x, \vec v - \bfq(G_2(x, \vec v))), \vec v) + \mcr_\eps\\[2pt] 
  = &   G_2(x, \bfq(G_2(x, \vec v)), \bfq(G_2(x, \vec v))) -  G_2(x, \vec v, \bfq(G_3(x, \vec v)))   - G_2(x, \bfq(G_3(x, \vec v)), \vec v)  \\
 &  +  G_2(x, \vec v, \bfq(G_2(x, \vec v, \bfq(G_2(x, \vec v))))) + G_2(x, \bfq(G_2(x, \vec v, \bfq(G_2(x, \vec v)))), \vec v)  \\
 & +  G_2(x, \bfq(G_2(x, \vec v, \bfq(G_2(x, \vec v))), \vec v) + G_2(x, \bfq(G_2(x, \bfq(G_2(x, \vec v)), \vec v), \vec v) + \mcr_\eps. 
\end{split}
\end{equation}
Using $\vec v = \sum_{i = 1}^4 \vec v^{(i)}$, we find that 
\beq
\begin{gathered}
 G_2(x, \vec v, \bfq(G_3(x, \vec v)))  + G_2(x, \bfq(G_3(x, \vec v)), \vec v)\\
 =  \sum_{(i ,j, k, l) \in \sigma(4)} \bigg( G_2(x, \vec v^{(i)}, \bfq(G_3(x, \vec v^{(j)}, \vec v^{(k)}, \vec v^{(l)})))   + G_2(x, \bfq(G_3(x, \vec v^{(i)}, \vec v^{(j)}, \vec v^{(k)})), \vec v^{(l)}) \bigg)
\end{gathered}
\eeq
Applying Prop.\ \ref{porder}, we see that the leading order term is achieved when the $G_i, i = 2, 3$ are $\widehat G_i$ and the leading terms are in $I^{4\mu + \frac{3}{2}}(\La^{\widehat g}_{q_0}\backslash \Theta)$. So we get the other piece of $\mch_2$. Using the same argument, we can obtain the terms in $\mch_3$ from the rest of terms in \eqref{eqg2}. The details are omitted here. 

Finally, let's consider the asymptotic expansion terms involving $\bar \mcj^{(i)}$ which we have dropped in \eqref{eqexpan}. Notice that the our source $J_\eps = (J^0_\eps, \bar J_\eps)$ where $\bar J_\eps^{a} = \sum_{i = 1}^4\eps_i \bar\mcj^{(i), a}, a = 1, 2, 3$ and according to \eqref{eqmcj0} we have
\beq 
\begin{split}
J^{0}_\eps %&= -(-\det  g)^{-\ha} \sum_{a= 1}^3 \int_0^t \p_a((-\det  g)^\ha   \bar J_\eps^a )ds 
=  -(-\det  g)^{-\ha} \sum_{a= 1}^3 \int_0^t \p_a((-\det  g)^\ha   \sum_{i = 1}^4\eps_i \bar\mcj^{(i), a} )ds 
\end{split}
\eeq 
So for example the order $\eps_i\eps_j, i\neq j$ terms involving $\bar\mcj^{(i)}$ in \eqref{eqexpan}  is a summation of the product of $\vec w_a, a = 1, 2, \cdots, 14$ and $\bar\mcj^{(i)}_b, b = 0, 1, 2, 3$ and possibly their derivatives in $y$ and integrals in $t$.  But $\WF(\bar\mcj^{(i)})$ is close to $\La_i\cap N^*Y_i$. By a wave front analysis and the argument in Lemma \ref{lmconormal}, we conclude that the wave front of these terms are included in $\La^{(1)}$. For the order $\eps_i\eps_j\eps_k, i<j < k$ terms involving $\bar\mcj^{(i)}$, a similar analysis tells that the wave front set of these terms are contained in $\La^{(1)}$ as well, and finally the wave front set of order $\eps_1\eps_2\eps_3\eps_4$ terms involving $\bar\mcj^{(i)}$ are contained in $\Theta$.  This finishes the proof of the claim.

\textbf{Step 2:} 
We recall the notation that $\vec v^{(i)} = (\dot g^{(i)}, \dot \phi^{(i)})$ is the linearized wave with source $\vec \mcj^{(i)}$, and $\vec v = (\dot g, \dot \phi) = \sum_{i = 1}^4  \eps_i \vec v^{(i)}$.  We already observed that the metric components $\dot g$ are all zero because of the choice of the source and the fact that the linearized equation \eqref{perlin} are decoupled. We use these to further identify the most singular terms in $\mch.$

(i) We see that 
\beq
H_{4, \alpha\beta}(x, \vec v, \vec v) = 0 \text{ and } P_{4, \alpha\beta}(x, \vec v, \vec v) = 0, \ \ \alpha, \beta = 0, 1, 2, 3
\eeq
because they all have the metric component. This implies that $\mch_1$ vanishes. 

(ii) Notice that $P_3(x, \vec v^{(i)}, \vec v^{(j)}, \vec v^{(k)})$ is at least quadratic in the metric components  and $P_2(x, \vec v^{(i)}, \vec v^{(j)})$ and $\hat H_3(x, \vec v^{(i)}, \vec v^{(j)}, \vec v^{(k)})$ are at least linear in the metric components. As the metric components of the linearized waves $\vec v^{(i)}$ vanish, the non-vanishing term in $\mch_2$ is 
\begin{equation}\label{eqmch}
\begin{gathered}
\mch_2 =  \sum_{(i, j, k, l) \in \sigma(4)} \bigg( \widehat H_3(x, \vec v^{(i)}, \vec v^{(j)}, \bfq (\widehat H_2(x, \vec v^{(k)}, \vec v^{(l)}))) + \widehat H_3(x, \vec v^{(i)}, \bfq (\widehat H_2(x, \vec v^{(j)}, \vec v^{(k)})), \vec v^{(j)}) \\  + \widehat H_3(x,  \bfq (\widehat H_2(x, \vec v^{(i)}, \vec v^{(j)})), \vec v^{(k)}, \vec v^{(l)}) \bigg).
\end{gathered}
\end{equation}

(iii) The components in $\widehat H_2(x, \vec v, \vec v)$ involving $\dot g$ components vanish. For analyzing terms in $\mch_3$, we identify the $\widehat H_2$ terms which have $\phi$ components. In the reduced Einstein equations, such $\widehat H_{2}$ terms are quadratic in $\phi$  and they come from the stress-energy term $T_{em}(\phi)$. For the Maxwell equations, $\widehat G_2$ terms are only linear in $\phi$ and this is when $\widehat G_2=P_2$. As a result, the electromagnetic potential component of $\bfq (\widehat G_2(x, \vec v^{(i)}, \vec v^{(j)})), i, j = 1, 2, 3, 4, i\neq j$ are zero but the metric components could be non-zero when we have $\bfq (\widehat H_2(x, \vec v^{(i)}, \vec v^{(j)}))$. Then only the electromagnetic components of $\widehat G_2(x, \vec v^{(k)}, \bfq (\widehat H_2(x, \vec v^{(i)}, \vec v^{(j)}))$ could be non-zero when $\widehat G_2$ is $P_2$. Thus, we can simplify $\mch_3$ to
\beqq\label{eqmch}
\begin{gathered}
\mch_3 =  -\sum_{(i, j, k, l) \in \sigma(4)} \widehat G_2(x, \bfq(\widehat H_2(x, \vec v^{(i)}, \vec v^{(j)})), \bfq(\widehat H_2(x, \vec v^{(k)}, \vec v^{(l)}))) \\
-\sum_{(i, j, k, l) \in \sigma(4)} \bigg( \widehat H_2(x, \vec v^{(i)}, \bfq(P_2(x, \bfq(\widehat H_2(x, \vec v^{(j)}, \vec v^{(k)})), \vec v^{(l)}))) \\
  + \widehat H_2(x, \bfq(P_2(x, \bfq(\widehat H_2(x, \vec v^{(i)}, \vec v^{(j)})), \vec v^{(k)})), \vec v^{(l)})\bigg),
\end{gathered}
\eeqq
in which all $\hat H_2$ are quadratic in $\phi$ components. This completes the proof of part (1).

%\textbf{Step 3:} Notice that every term in $\mcu^{(4)}$ is a linear combination of the model terms $\mcy_3$ in Section \ref{intera} so that the principal symbols of $\mcu^{(4)}$ are fourth order homogeneous polynomials of the principal symbols of $\vec v^{(i)}$. 

(2) The proof is the the same as that of Prop.\ 4.1 of \cite{LUW} for the scalar case. If $K_i, i = 1, 2, 3, 4$ do not intersect, the singularities of $\mcu^{(4)}$ are at most conic (when three of $K_i$ intersect) and the wave front set is contained in $\Theta$.  
\epf

Finally, we show that by choosing distorted plane waves, the newly generated singularities of $\mcu^{(4)}$ are not always vanishing. 
\begin{prop}\label{symb} 
Suppose that geodesics $\gamma_{x^{(j)}, \theta^{(j)}}(\R_+), j = 1, 2, 3, 4$ intersect at $\gamma_{x^{(j)}, \theta^{(j)}}(t_j) = q_0$ and  $\dot\gamma_{x^{(j)}, \theta^{(j)}}(t_j)$, $j=1,2,3,4$, are linearly independent. Let $\vec v^{(i)} = (\dot g^{(i)}, \dot \phi^{(i)})$ be distorted plane waves propagating near geodesics  $\gamma_{x^{(j)}, \theta^{(j)}}(\R_+)$, described in the beginning of Section 5.2, so that $\vec v^{(i)}\in I^{\mu - \ha}(N^*K_i; \bold{B}^{14})$ in $ M(T_0)\backslash Y_i.$ For the fourth order interaction term $\mcu^{(4)}$ produced by the waves $\vec v^{(i)}$, we denote by $\mcu^{(4), met}$ the metric component and by $\mcu^{(4), em}$ the electromagnetic potential component. Let $(q, \eta)\in \La_{q_0}^{\widehat g}\backslash \Theta$ which is joined to $(q_0, \zeta)\in \La_{q_0}$ by bicharacteristics, and $ \zeta^{(i)} \in N^*_{q_0}K_i$ 
be such that $\zeta = \sum_{i = 1}^4 \zeta^{(i)}$.

We have the following conclusions for the principal symbol of $\mcu^{(4)}$ in $I^{4\mu + \frac 32}(\La_{q_0}^{\widehat g}\backslash \Theta; \bold{B}^{14})$.
\begin{enumerate}
\item [(i)] $\sigma(\mcu^{(4), em})(q, \eta) = 0$.
\item  [(ii)] $\sigma(\mcu^{(4), met})(q, \eta) = \sigma(\bfq_{\widehat g})(q, \eta, q_0, \zeta)
 \mcp(\zeta^{(1)}, \zeta^{(2)}, \zeta^{(3)}, \zeta^{(4)}, A^{(1)}, A^{(2)}, A^{(3)}, A^{(4)})$, where 
$ A^{(i)}=\sigma(\dot \phi^{(i)})(q_0,\zeta^{(i)})$ are the principal symbols of $\dot \phi^{(i)}$ at $(q_0,\zeta^{(i)})\in N^*K_i$,
$\mcp$ is a $4\times 4$ matrix whose component $\mcp^{\alpha\beta}, \alpha, \beta = 0, 1, 2, 3$ are homogeneous polynomial of degree $4$ in $A^{(i)}\in (\mbr^4)^4$, $i = 1, 2, 3, 4$ whose coefficients are
real-analytic functions defined on 
% and $\mcp^{\alpha\beta}(\cdot)$ are non-vanishing on any {\ctext open set of \HOX{component $\zeta$ is removed from definition of $\mcx$.}
\beq
\begin{gathered}
\mcx(\zeta) = \{(\zeta^{(1)}, \zeta^{(2)}, \zeta^{(3)}, \zeta^{(4)}) \in (L_{q_0}^{*}M)^4:\ \sum_{i = 1}^4\zeta^{(i)}=\zeta  \}.
 \end{gathered}
\eeq
 
 \item  [(iii)] Assume that $q$ has a neighborhood $B$  such that the intersection of $B$ and
 the light cone $\mathcal L^+_{q_0}$ is a smooth 3-dimensional surface $S=B\cap  \mathcal L^+_{q_0}$.
Let $\text{Sym}^2(T_{q}M\otimes T_{q}M)$ denote the symmetric two tensors at $q_0$ and consider its
 subspace $\text{Sym}^2(T_{q}S\otimes T_{q}S)$.  Then any open subset  $W\subset \mcx(\zeta)$ 
 contains $\vec \zeta=(\zeta^{(1)}, \zeta^{(2)}, \zeta^{(3)}, \zeta^{(4)})$ such that $\zeta^{(1)}, \zeta^{(2)}, \zeta^{(3)}, \zeta^{(4)}$ are linearly independent 
 and there are $ A^{(i)}=\sigma(\dot \phi^{(i)})(q_0,\zeta^{(i)}), i = 1, 2, 3, 4$ satisfying the  linearized gauge condition $\widehat g^{\alpha\la}\zeta^{(i)}_{\la}A^{(i)}_\alpha =0 $ and $V\in Sym^2(T_qS\otimes T_qS)$ 
 such that  the fourth interaction term $\mcu^{(4)}$,  corresponding to the interaction of waves $\vec v^{(i)}$,  satisfies $\langle \sigma(\mcu^{(4), met})(q, \eta), V\rangle\not =0$.
\end{enumerate}
\end{prop}

We remark that parts (ii)-(iii) of the proposition can be interpreted as following: 
Consider any $(x^{(i)}, \theta^{(i)}), i = 1, 2, 3, 4$ so that the geodesic $\gamma_{x^{(i)}, \theta^{(i)}}$ intersect at $q_0$. For any $(q,\eta)\in \La_{q_0}^{\widehat g}\backslash \Theta$, one can find $(\widetilde x^{(i)}, \widetilde \theta^{(i)})$ in an arbitrarily small neighborhood of $(x^{(i)}, \theta^{(i)})$ and distorted plane waves $\vec v^{(i)}$ associated with $\gamma_{\widetilde x^{(i)}, \widetilde \theta^{(i)}}$ such that the interaction term  $\mcu^{(4)}$ has non-vanishing singularities at $(q,\eta)$ i.e.\ $\mcu^{(4)}$ is not smooth at $q.$ 
Moreover, some of these singularities have such polarization that they can be observed by taking inner products of the field $\mcu^{(4)}$ and
tensor products of vector fields that are tangent to the light cone $ \mathcal L^+_{q_0}$.

\bpf[Proof of Prop.\ \ref{symb}] 
\textbf{Claim (i).} 
 From Prop.\ \ref{intert}, we can write the symbol of $\mcu^{(4)}$  at $(q,\eta)$ as 
\beqq\label{symint}
\sigma(\mcu^{(4)})(q, \eta) = \sigma(\bold{Q})(q, \eta, q_0, \zeta) \sigma(\mch)(q_0, \zeta),
\eeqq
where $\mch$ is given in Prop.\ \ref{intert}. 
From the proof of  Prop.\ \ref{intert}, we've seen that it suffices to consider in $\mch$ the terms $\widehat H_2$ that are quadratic in $\phi$ the electromagnetic components, and such terms only come from the reduced Einstein equation. In particular, the electromagnetic component  $\widehat H_{2, \mu}, \mu = 0, 1, 2, 3$ are all zero. However, from the Maxwell equation \eqref{pernon1}, we observe that the nonlinear terms $\widehat P_{2, \mu}, \widehat H_{3, \mu}, \mu = 0, 1, 2, 3$ are always linear in $\phi$. Therefore, from the expression of $\mch$ above, we see that $\mch_{\mu} = 0, \mu = 0, 1, 2, 3$. This implies that the principal symbol $\sigma(\mcu^{(4), em})(q, \eta) = 0$.

\textbf{Claim (ii).} 
We denote below $\vec\zeta=(\zeta^{(1)}, \zeta^{(2)}, \zeta^{(3)}, \zeta^{(4)})$ and $\vec A=( A^{(1)}, A^{(2)}, A^{(3)}, A^{(4)})$. As we already discussed in part (i), in $\mch$ it suffices to consider $\widehat H_2$ which are quadratic in $\phi$. The terms in \eqref{symbh1} are similar to $\mcy_3, \mcy_4, \mcy_5$ in  discussion after Prop.\ \ref{porder} for which we have the principal symbol, see \eqref{eqy3y5}. Using these formulas for \eqref{symbh1}, we see that the principal symbol $\sigma(\mch)(q_0, \zeta)= \mcp(\vec \zeta,\vec A)$ is a fourth order polynomial of $\vec A^{(i)}_a, i = 1, 2, 3, 4, a = 0, 1, 2, 3$. Therefore, they are real analytic functions of $\vec\zeta$  and $\vec A$, see more details in the proof of (iii). This proves claim (ii).

\textbf{Claim (iii).}
Now we compute the metric component of $\sigma(\mch)(q_0, \zeta)$ in \eqref{symint} in terms of $\vec\zeta$ and $\vec A$ to show they are non-zero functions. Without loss of generality, we can use  local coordinate near $q_0$ such that $\widehat g$ is the standard Minkowski metric $h$ in $\R^4$ at $q_0$. Then we can use the Einstein-Maxwell equations in Minkowski space-time found in Section \ref{gwave}. We remark that the derivatives of the metric at $q_0$ in the local coordinates do not contribute in the principal symbols and thus we can do the considerations below in the case when the space $(M,\widehat g)$ is the Minkowski space.

We assume that $\zeta$ is chosen to be $\zeta  = (9/4, -7/4, 1, 1)$. Consider the following set of vectors 
\beqq\label{vectors}
\begin{gathered}
\zeta^{(1)}=(1,0,1,0),\ \ \zeta^{(2)}=(1,0,0,1), \ \ 
\zeta^{(3)}=  (-3/4, -3/4, 0,0),\ \ \zeta^{(4)} =  (1,-1,0,0).
\end{gathered}\hspace{-6mm}
\eeqq
These are light-like vectors and $\sum_{i = 1}^4\zeta^{(i)} = (9/4, -7/4, 1, 1)=\zeta$. Let $\imath$ be the imaginary unit i.e.\ $\imath^2 = -1$. We shall take 
\beqq\label{first example}
\begin{gathered}
A^{(1)}= -(0,0,0,\imath),\ \ A^{(2)}= -(0,0,\imath,0), \ \ 
A^{(3)}=  -(0, 0, 0,\imath),\ \ A^{(4)} =  -(0,0,0,\imath)
\end{gathered}
\eeqq
as the symbols of $\dot \phi^{(i)}$ at $\zeta^{(i)}$.
{We will also denote the symbols given in (\ref{first example}) by ${}^1A^{(i)}=A^{(i)}$ and ${}^1\vec A=({}^1A^{(i)})_{i=1}^4$.}
 With these choices, we see that the symbols $A^{(i)}$ satisfy the microlocal linearized gauge condition i.e.\
\beqq\label{recall  microlocal linearized gauge condition}
h^{\alpha\la} \zeta^{(i)}_\la A^{(i)}_\alpha = 0, \ \ i = 1, 2, 3, 4.
\eeqq 
{Next we compute  the principal symbol of $F^{(i)}$ at $(q_0, \zeta^{(i)})$ (in $I^{\mu+1}$) using $\sigma(F^{(i)}_{\alpha\beta}) = \imath\zeta^{(i)}_\alpha A^{(i)}_\beta - \imath\zeta^{(i)}_\beta A^{(i)}_{\alpha}$. Due to the choices of the vectors, many components of $F$ vanish.} Actually, the non-zero ones are
\beq
\begin{gathered}
F^{(1)}_{03} = 1, \ \ F^{(1)}_{23} = 1; \ \ F^{(2)}_{02} = 1, \ \ F^{(2)}_{32} = 1; \\
F^{(3)}_{03} = -3/4, \ \ F^{(3)}_{13} = -3/4; \ \ F^{(4)}_{03} = 1, \ \ F^{(4)}_{13} = -1,
\end{gathered}
\eeq
and their (anti-)symmetric terms. By straightforward computations carried out in Appendix \ref{secapp}, we find that 
\beqq\label{matrix with five non-zero elements}
\mcp = \sigma(\mch)(q_0, \zeta) =  c_\pi
  \begin{pmatrix}
    * & * & -5.3 & -0.7 \\
    * & * & 2.1 & 4.9 \\
-5.3 & 2.1 & * & -4.5\\
-0.7 & 4.9 & -4.5& * 
  \end{pmatrix}
\eeqq
where $*$ stands for the elements which we didn\rq{}t compute and $c_\pi = (2\pi)^{-3}$ is a non-zero constant. 
This proves that $\mcp$ hence the symbol of $\mcu^{(4), met}$ is non-vanishing for the choice of $\zeta^{(i)}$ and $A^{(i)}={}^1A^{(i)}$. 

We need to construct four other examples of choices of  $A^{(i)}$ such that the symbol of $\mch^{\alpha\beta}$ is non-vanishing.
%It remains to show part (3). 
%For convenience, we {\ctext denote $\vec\zeta=(\zeta^{(1)}, \zeta^{(2)}, \zeta^{(3)}, \zeta^{(4)})$ and $\vec A=( A^{(1)}, A^{(2)}, A^{(3)}, A^{(4)})$.}
%$\vec A = (A^{(1)}, A^{(2)}, A^{(3)}, A^{(4)})$. 
%According to the result in Step 2, we found one set of vectors which we now denote by ${}^1\vec A$ such that the symbol $\sigma(\mch({}^1\vec A))$ has five non-vanishing components. 
%the dependency of $\mch$ on $\vec A$. 
For simplicity, we define an operator $T$ from $\text{Sym}^2(T_{q_0}M\otimes T_{q_0}M)\simeq \mbr^{10}$  to $\mbr^5$ as follows:
\beq
T(B) = (B_{02}, B_{03}, B_{12}, B_{13}, B_{23}),
\eeq
where $B\in\text{Sym}^2(T_{q_0}M\otimes T_{q_0}M)$. In this notation, we have that 
\beq
T\mcp(\vec \zeta,{}^1\vec A) = c_\pi (-5.3, -0.7, 2.1, 4.9, -4.5).
\eeq
In Appendix \ref{secapp}, we construct four set of vectors ${}^a \vec A, a = 2, 3, 4, 5$ that satisfy the microlocal linearized gauge conditions ({recall  microlocal linearized gauge condition}) and that 
  $T\mcp(\vec \zeta,{}^a\vec A)$, $ a = 1, 2, 3, 4, 5$ are non-zero and linearly independent in $\mbr^5$. 
% This is done in Appendix \ref{secapp2}. 
Actually, the proof is based on the same kind of computation as for ${}^1 \vec A$. 
 We emphasize that in these examples we keep the vectors $\zeta^{(i)}$ fixed and we construct different choices of $\vec A$.
%In particular, this implies that there are ${}^a \vec A$  such that $\sigma(\mch({}^a \vec A)),$ $ a = 1, 2, 3, 4, 5$ are linearly independent in $\text{Sym}^2(T_{q_0}M\otimes T_{q_0}M)$.

Next we consider the dependency of  $\sigma(\mch)$, or more precisely, of $\mcp(\vec \zeta,\vec A)$ on the variables
$\vec \zeta$  and $\vec A$. We will consider the case when $\zeta^{(i)}=\alpha_i \xi^{(i)}$, with $\alpha_i\in \R$ and $\xi^{(i)}\in L_{q_0}^{*, +}M$, and 
{denote  $\vec\xi=(\xi^{(1)}, \xi^{(2)}, \xi^{(3)}, \xi^{(4)})\in (L_{q_0}^{*, +}M)^4$. 
Also, let
\beq
\begin{gathered}
\mathcal Y  = \{(\xi^{(1)}, \xi^{(2)}, \xi^{(3)}, \xi^{(4)}, A^{(1)}, A^{(2)}, A^{(3)}, A^{(4)}) \in (L_{q_0}^{*, +}M)^4\times (\mbr^4)^4:\\
 \text{$A^{(i)}$ satisfy the linearized gauge condition $\widehat g^{\alpha\la}\xi_{\la}^{(i)}A^{(i)}_\alpha =0$. }  \}.
 \end{gathered}
\eeq
When $\xi^{(1)}, \xi^{(2)}, \xi^{(3)}, \xi^{(4)}$ are linearly independent,
let $\alpha_i(\zeta,\vec\xi)\in \mbr$ be such that $\alpha_i=\alpha_i(\zeta,\vec\xi)$ satisfy
 $\zeta = \sum_{i = 1}^4\alpha_i\xi^{(i)}$. Then $\alpha_i(\zeta,\vec\xi)$  is a quotient of real-analytic functions, that is, $\alpha_i(\zeta,\vec\xi)=p_i(\zeta,\vec\xi)/q_i(\zeta,\vec\xi)$  where
$p_i(\zeta,\vec\xi)$ and $q_i(\zeta,\vec\xi)$  are real-analytic functions of $(\zeta,\vec\xi)\in  (L_{q_0}^{*, +}M)^5$ and $q_i(\zeta,\vec\xi)$   is not identically vanishing.  
Next, let us consider $\zeta= (9/4, -7/4, 1, 1)$. Here vectors $\zeta^{(i)}$ are given in (\ref{vectors}) and we use $\xi^{(i)}=\tilde \zeta^{(i)} = \zeta^{(i)}/\alpha_i$  where $(\alpha_1,\alpha_2,\alpha_3,\alpha_4)=(1,1,-3/4,1)$. 
Let
$$ S(\vec\xi,\vec A)=\mcp(\alpha_1(\zeta,\vec\xi)\xi^{(1)},
\alpha_2(\zeta,\vec\xi)\xi^{(2)},\alpha_3(\zeta,\vec\xi)\xi^{(3)},\alpha_4(\zeta,\vec\xi)\xi^{(4)},\vec A),\quad \hbox{for }(\vec\xi,\vec A)\in \mathcal Y.
%\sigma(\mch(\vec\xi,\vec A))(q_0,\zeta_0).
$$ 
Note we have chosen $\zeta$ to a fixed light-like vector, function $S(\vec\xi,\vec A)$,  that gives the symbol
$\sigma(\mch)(q_0, \zeta)$   of the fourth interaction source corresponding to the waves $\vec v^{(i)}$ having directions $\xi^{(i)}$  and symbols $A^{(i)}$,
%
%and then $$
%\sigma(\mch)(q_0, \zeta)=S(\vec\xi,\vec A)=\mcp(\alpha_1(\zeta,\vec\xi)\xi_1,
%\alpha_2(\zeta,\vec\xi)\xi_2,\alpha_3(\zeta,\vec\xi)\xi_3,\alpha_4(\zeta,\vec\xi)\xi_4,\vec A)
%$$
can be considered as a quotient of real analytic functions of $\vec\xi$  and $\vec A$.

}

Now we let $e_a \in \text{Sym}^2(T_{q_0}M\otimes T_{q_0}M), a = 1, 2, 3, 4, 5$ be constructed as follows.
\beq
\begin{gathered}
e_{1, 02} = e_{1, 20} = 1, \ \  e_{2, 03} = e_{2, 30} = 1, \\
e_{3, 12} = e_{3, 21} = 1, \ \  e_{4, 13} = e_{4, 31} = 1,\ \ e_{5, 23} = e_{5, 32} = 1
\end{gathered}
\eeq
and all the rest of the components are zero.  Let $Z = \text{span}\{e_a, a = 1, 2, 3, 4, 5\}$ be a $5$ dimensional subspace of $\text{Sym}^2(T_{q_0}M\otimes T_{q_0}M)$. Then, $T$  is a projection onto the space $Z$.

{Consider the map
\beq
D(\vec\xi,({}^a\vec A)_{a=1}^5) =%=(\vec\xi,({}^a\vec A)_{a=1}^5)\mapsto
\det\bigg(\big (\bra S (\vec\xi,{}^a\vec A),e_b\cet\big)_{a,b=1}^5\bigg).
\eeq
By substituting in $D(\vec\xi,({}^a\vec A)_{a=1}^5)$  the  vectors $\vec \xi$ and the symbols
given in formula (\ref{first example}) and  the other four sets, for which $TS (\vec\xi,{}^a\vec A)$, $a=1,2,3,4,5,$ are linearly independent, we see 
that the function $D:\mathcal Y \to \mathbb C$ obtains at some point $(\vec\xi,({}^a\vec A)_{a=1}^5)\in \mathcal Y$
a finite and non-zero value, and thus it is not identically vanishing. Since
$\mathcal Y$ is a real-analytic manifold and $D:\mathcal Y\to \mathbb C\cup\{\infty\}$ is a 
quotient of real-analytic functions that is not identically vanishing, we see that $D:\mathcal Y \to \mathbb C$ does not vanish or is infinity in any open subset
of $\mathcal Y$. Thus we see that in any open set $W\subset (L_{q_0}^{*, +}M)^4$ there are $\vec\xi=(\xi^{(1)}, \xi^{(2)}, \xi^{(3)}, \xi^{(4)}) \in W$  and $(\vec\xi ,({}^a\vec A)_{a=1}^5)\in \mathcal Y$ 
such that vectors $\xi^{(1)}, \xi^{(2)}, \xi^{(3)}, \xi^{(4)}$  are linearly independent and matrices $S (\vec\xi,{}^a\vec A)$, $a=1,2,\cdots,5$ are linearly independent.
Let $ \zeta^{(i)}=\alpha_i(\zeta, \vec\xi)\xi^{(i)}$, $i=1,2,3,4$ and ${}^a\vec v^{(i)} = ({}^a\dot g^{(i)}, {}^a\dot\phi^{(i)})\in I^{\mu-\ha}(N^*K_i)$, $a=1,2,3,4,5$ be distorted plane waves described in the beginning of this section, such that ${}^a\dot g^{(i)}=0$  and $\sigma({}^a\dot \phi^{(i)})(q_0,\zeta_i)={}^aA^{(i)}$. 
%For
%$c=(c_a)_{a=1}^5\in \mathbb C^5$, we define 
%\beqq\label{c combinations}
%\vec v^{(i)}(x)=\sum_{a=1}^5 c_a\,\cdotp {}^a\vec v^{(i)}(x),\quad %\phi^{(i)}(x)=\sum_{a=1}^5 c_a\,\cdotp  {}^a\vec \phi^{(i)}(x),\quad
% A^{(i)}_c=\sum_{a=1}^5 c_a\,\cdotp  {}^a A^{(i)}.
%\eeqq
When ${}^a\mcu^{(4)}$ is the wave produced by the linearized waves ${}^a\vec v^{(i)}(x)$ corresponding to the vectors $ \vec\zeta$ and matrixes ${}^aA^{(i)}$, % given in (\ref{c combinations}), 
we see that
\beqq\label{inner products}
\langle \sigma({}^a\mcu^{(4), met})(q, \eta), V\rangle  = \langle S(\vec\xi,{}^a\vec A),\mcr^T V\rangle,\quad \mcr = \sigma(Q_{\widehat g})(q, \eta, q_0, \zeta)
\eeqq
where $V\in \text{Sym}^2(T_{q}M\otimes T_{q}M)$.
Since the linear map $\mcr$ is bijective,
and the space $\text{Sym}^2(T_{q}S\otimes T_{q}S)$ has co-dimension 4 in $\text{Sym}^2(T_{q}M\otimes T_{q}M)$,
also the space $\mcr^T(\text{Sym}^2(T_{q}S\otimes T_{q}S))$ has co-dimension 4 in $\text{Sym}^2(T_{q_0}M\otimes T_{q_0}M)$.
If the dualities (\ref{inner products}) would vanish for all $a=1,2,3,4,5$, 
and for all $\mcr^T V\in \mcr^T(\text{Sym}^2(T_{q}S\otimes T_{q}S))$, then
all the vectors 
$S(\vec\xi,{}^a\vec A)$, $a=1,2,3,4,5$ would be in the four dimensional space $(\mcr^T(\text{Sym}^2(T_{q}S\otimes T_{q}S)))^\perp.$
However, since the vectors  $S(\vec\xi,{}^a\vec A)$,  $a=1,2,3,4,5$ are linearly independent, this is not possible.

% the space $\text{span}\{ S(\vec\xi,{}^a\vec A)\ c=(c_a)_{a=1}^5\in \mathbb C^5\}$  has dimension 5,
%we see that there exists $c=(c_a)_{a=1}^5\in \mathbb C^5$ such that for the linearized waves (\ref{c combinations}
The above implies that for some $V\in \text{Sym}^2(T_{q}S\otimes T_{q}S)$ and for
 some $a\in \{1,2,3,4,5\}$, the linearized waves $\vec v^{(i)}(x)$, 
propagating near geodesics $\gamma_{q_0,\zeta_i^{\sharp}}$  and having  symbols $A^{(i)}={}^aA^{(i)}$ at $(q_0,\zeta^{(i)})$, %
are such that the wave $\mcu^{(4)}$ produced by the fourth order interaction  
satisfies $\langle \sigma(\mcu^{(4), met})(q, \eta), V\rangle \not=0$. This proves the claim (iii).}
%
%that for any $V\in Z$ if $\langle \sigma(\mch({}^a\vec A)), V\rangle  = 0$ for $a = 1, 2, 3, 4, 5$, then $V = 0$. In other words, for any $V\in Z, V\neq 0$, we can find $\vec A$ so that $\langle \sigma(\mch(\vec A)), V\rangle  \neq 0$. Because $\mcr \doteq \sigma(Q_g)(q, \eta, q_0, \zeta)$ is an invertible matrix, we conclude that 
%\beq
%\langle \sigma(\mcu^{(4), met})(q, \eta), V\rangle  = \langle \sigma(\mch)(q_0, \zeta),\mcr^T V\rangle  \neq 0
%\eeq
%for any $V\in (\mcr^{T})^{-1} (Z)$ which is a $5$ dimensional subspace of $\text{Sym}^2(T_{q}M\otimes T_{q}M)$. This finishes the proof of part (3) hence completes the proof of the proposition.
 \epf

%==================================%
\section{Determination of the space-time}\label{secinv}
We finish the proof of Theorem \ref{main1}. Since the rest of the proof is essentially the same as that of \cite[Theorem 1.1]{KLU1} and the argument is quite involved, we shall only repeat the key argument and refer the interested readers to \cite{KLU1} for more details. Basically, we need to consider two things. First, we analyzed the singularities in wave gauge and we need to show such singularities can be observed in the data set {in Fermi coordinates, that is roughly speaking, by freely falling observers that move in the perturbed spacetime}. This follows from \cite[Section 4]{KLU1}. Second, from the observed singularities, one can determine the earliest light observation set and the space-time structure. {This is contained in \cite[Section 4]{KLU}, see also \cite[Section 3.5 and 5]{KLU1}.} 

We first show that an analogue of \cite[Lemma 4.2]{KLU} holds for our problem, which is the key lemma for changing observations to Fermi coordinate. We remark that Prop.\ \ref{intert} and Prop.\ \ref{symb} are the analogous of Prop.\  3.3 and Prop.\ 3.4 of \cite{KLU1}. Following \cite[Section 4]{KLU1}, we say that the {\em interaction condition (I)} is satisfied for $y\in V$ with light-like vectors $(\vec x, \vec \xi) = ((x^{(i)}, \xi^{(i)}))_{i = 1}^4$ and $t_0 > 0$, if \\
\indent \textbf{(I) }  There exist $q_0 \in \bigcap_{i = 1}^4 \gamma_{(x^{(i)}(t_0), \xi^{(i)}(t_0))}((0, \tau_i))$, $\zeta \in L_{q_0}^+M$ and $t\geq 0$ such that $y = \gamma_{q_0, \zeta}(t)$. Here we recall that $\tau_i$ are defined in the beginning of Section \ref{singu} i.e.\ $\gamma_{(x^{(i)}(t_0), \xi^{(i)}(t_0))}(\tau_i)$ is the first conjugate point. 
 
Also, we say that $y\in V$ satisfies the singularity {\em detection condition (D)} with light-like directions $(\vec x, \vec \xi)$ and $t_0, \widehat s >0$ if \\
\indent \textbf{(D) } For any $s, s_0\in (0, \widehat s)$ and $i = 1, 2, 3, 4$, there are $(\widetilde x^{(i)}, \widetilde \xi^{(i)})$ in the $s$-neighborhood of $(x^{(i)},  \xi^{(i)})$, $(2s)$-neighborhood $B^{(i)}$ of $x^{(i)}$ (these neighborhoods are defined using the Riemannian metric $\widehat g^+$) satisfying 
\beq
J^+_{\widehat g}(B^{(j)}) \cap J^-_{\widehat g}(B^{(j)})\subset V_{\widehat g} \text{ and } J^+_{\widehat g}(B^{(j)})\cap J^+_{\widehat g}(B^{(k)})=\emptyset \text{ for all $j\not =k$.}
\eeq
Also, consider sources $\mcj^{(i)}$ constructed as in the beginning of Section \ref{singu}. Let $(g_\eps, \phi_\eps)$ be the solutions to \eqref{einmax} with sources $J_\eps$ such that $\p_{\eps_i} J_\eps|_{\eps_i = 0} = \mcj^{(i)}$ and $\overline J_\eps^{(i)}$ are supported in $B^{(i)}$. {Let $\Phi_{g_\eps}$ be the normal coordinates with respect to $g_\eps$ centered at the point $y_\eps=\Psi_\eps^{-1}(\Psi_0(y))$, where $\Psi_\eps$ are the Fermi coordinates, see (\ref{Fermi coordinates}).}
 Then $\p_{\eps_1}\p_{\eps_2}\p_{\eps_3}\p_{\eps_4}\Phi_{g_\eps}^*g_\eps|_{\eps_1 = \eps_2 = \eps_3 = \eps_4 = 0}$ is not smooth {near $\Phi_0(y)$}. We remark that due to this condition we can determine the wave gauge coordinates in the sets $J^+_{\widehat g}(B^{(j)})\cap  J^-_{\widehat g}(B^{(j)})$. \\
 
Using \cite[Lemma 4.1]{KLU1} and Prop.\ \ref{intert} and \ref{symb}, we can adapt the proof of \cite[Lemma 4.2]{KLU1} to our case and conclude that using the data set $\mcd(\delta)$, we can determine whether the condition $(D)$ is valid for the given point $y\in V$ or not. We remark that here the key to adapt the proof of \cite[Lemma 4.2]{KLU1} is Prop.\ \ref{symb} part (3), which is required by \cite[Lemma 4.1]{KLU1}. {Then, \cite[Lemma 4.2]{KLU1}, parts (1)-(2),
tell that the interaction condition (I) and the singularity detection condition (D) are equivalent when $y$  is in not the future of any cut point of geodesics
$\gamma_{(x^{(i)}(t_0), \xi^{(i)}(t_0))}$.}

Next, the light observation set of $q\in M$ in $V$ is defined as $\mcp_V(q) = \mcl_q^+\cap V$. The earliest light observation set is defined as 
\beq
\mce_V(q) = \{x\in \mcp_V(q): \hbox{there is no $y\in \mcp_V(q)$ and future-pointing time-like path}\\
 \hbox{$\alpha:[0, 1] \rightarrow V$ such that $\alpha(0) = y$ and }\alpha(1) = x \}\subset V,
\eeq
 see \cite[Def.\ 1.1]{KLU1}. For $W\subset M$ open, the collection of the earliest light observation sets with source points in $W$ is $\mce_V(W) = \{\mce_V(q): q\in W\}.$ In particular, if $q_0$ is the interaction point as in Prop.\ \ref{intert}, then $\mce_V(q_0) \subset \mathcal{N}((\vec x, \vec \theta), t_0). $ 
From Section 2.2.1 of \cite{KLU1}, we know that $\mce_V(q_0)$ contains a $3$-dimensional submanifold hence is not empty. 

The argument in \cite[Section 5]{KLU1} can be applied to show that from the data set $\mcd(\delta)$, we can determine the set $\mce_{V, \widehat g}(I_{\widehat g}(p_-, p_+))$. This means that we can produce ``artificial point source" at any point $q$ of $I_{\widehat g}(p_-,p_+)$ and determine the intersection of the observation set $V_{\widehat g}$ and the light cone emanating from the point $q.$ 
Now for the set up of Theorem \ref{main1}, we can {apply \cite[Theorem 1.2]{KLU}  to conclude that $\widehat g^{(1)}$ is conformal to $\widehat g^{(2)}$. Moreover, since $\widehat g^{(1)}, \widehat g^{(2)}$ are Ricci flat, we can further apply \cite[Corollary 1.3 ]{KLU}} to conclude that the conformal diffeomorphism is an isometry i.e.\ there exists a diffeomorphism $\Psi: I_{\widehat g^{(1)}}(p_-, p_+) \rightarrow I_{\widehat g^{(2)}}(p_-, p_+)$ such that $\Psi^*\widehat g^{(2)} = \widehat g^{(1)}$. This finishes the proof of Theorem \ref{main1}.

\appendix
%%%%%%%%%%%%%%%%%%%%%%%%%%%%%%%%%%%%%%%%
\section{Computation of principal symbols of the interaction term}\label{secapp}
We complete the calculation required in the proof of Prop.\ \ref{symb}. Our goal is to compute the principal symbols of $\mch$ given in Prop.\ \ref{intert}. For convenience, we split $\mch$ into three parts
\beqq\label{symint1}
\begin{gathered}
\mch_1 =  \sum_{(i, j, k, l) \in \sigma(4)} \bigg( \widehat H_3(x, \vec v^{(i)}, \vec v^{(j)}, \bfq (\widehat H_2(x, \vec v^{(k)}, \vec v^{(l)}))) \\
+ \widehat H_3(x, \vec v^{(i)}, \bfq (\widehat H_2(x, \vec v^{(j)}, \vec v^{(k)})), \vec v^{(j)}) + \widehat H_3(x,  \bfq (\widehat H_2(x, \vec v^{(i)}, \vec v^{(j)})), \vec v^{(k)}, \vec v^{(l)}) \bigg)\\
 \mch_2 =  -\sum_{(i, j, k, l) \in \sigma(4)} \widehat G_2(x, \bfq(\widehat H_2(x, \vec v^{(i)}, \vec v^{(j)})), \bfq(\widehat H_2(x, \vec v^{(k)}, \vec v^{(l)}))) \\
\mch_3 =  -\sum_{(i, j, k, l) \in \sigma(4)} \bigg( \widehat H_2(x, \vec v^{(i)}, \bfq(P_2(x, \bfq(\widehat H_2(x, \vec v^{(j)}, \vec v^{(k)})), \vec v^{(l)}))) \\
  + \widehat H_2(x, \bfq(P_2(x, \bfq(\widehat H_2(x, \vec v^{(i)}, \vec v^{(j)})), \vec v^{(k)})), \vec v^{(l)})\bigg)
\end{gathered}
\eeqq
For the computation, we recall that we use local coordinate near $q_0$ so that $\widehat g$ is the standard Minkowski metric at $q_0$, and we found the Einstein-Maxwell equations in Minkowski space-time in Section \ref{gwave}. Let $(u, \phi) = (g - \widehat g, \phi)$ and  $F = d\phi$. For our purpose, we need the nonlinear terms $\widehat H_2, \widehat H_3$ which are quadratic in $\phi$, which can only come from the stress-energy tensor. For $\alpha, \beta = 0, 1, 2, 3,$ we have
\beq
\begin{split}
(\widehat H_2)_{\alpha\beta} &= 2(-h^{\la\mu}F_{\alpha\la} F_{\beta\mu} + \frac{1}{4}h_{\alpha\beta}h^{\la\gamma}h^{\mu\delta}F_{\gamma\delta}F_{\la\mu}), \\ 
%= -h^{\la\la'}F_{\alpha\la} F_{\beta\la'} +   \frac{1}{4}h_{\alpha\beta}h^{\la\la'}h^{\mu\mu'}F^2_{\la\mu'},
&=  2\big(\sum_{\la = 0}^3 -h^{\la\la}F_{\alpha\la} F_{\beta\la} + \sum_{\la, \mu = 0}^3 \frac{1}{4}h_{\alpha\beta}h^{\la\la}h^{\mu\mu}F^2_{\la\mu}\big)
\end{split}
\eeq
and for $\alpha\neq \beta$, we have 
\beq
\begin{split}
(\widehat H_3)_{\alpha\beta} &= 2(-h^{aa'}h^{bb'}u_{ab}F_{\alpha a'} F_{\beta b'} + \frac{1}{4}u_{\alpha\beta}h^{aa'}h^{bb'}F_{ab}F_{a'b'})\\
&= \sum_{a, b = 0}^3 2(-h^{aa}h^{bb}u_{ab}F_{\alpha a} F_{\beta b} + \frac{1}{4}u_{\alpha\beta}h^{aa}h^{bb}F^2_{ab}).
\end{split}
\eeq
Also, we need
\beq
P_2(x, u, v) = (h^{-1} u h^{-1})^{pq}\frac{\p^2 }{\p x^p \p x^q} u.
\eeq
Sometimes, it is convenient to write $\widehat H_2$ in matrix form. Let $H = (h_{ij})$ and $F^{(i)} = d\dot \phi^{(i)}, i = 1, 2, 3, 4.$ Notice that $F^{(i)}$ are anti-symmetric so $F^{(i)}_{\beta\mu} = -F^{(i)}_{\mu \beta}, \mu, \beta = 0, 1, 2, 3$. Thus we obtain that 
\beq
\begin{gathered}
\widehat H_2(x, v^{(i)}, v^{(j)}) =  2 F^{(i)} H F^{(j)} + \ha H \text{Tr}(F^{(i)}F^{(j)}), \text{  where }\\
\text{Tr}(F^{(i)}F^{(j)}) = \sum_{\la, \mu = 0}^3 h^{\la\la}h^{\mu\mu} F^{(i)}_{\la\mu} F^{(j)}_{\la\mu}, \ \ 1 \leq i < j \leq 4.
\end{gathered}
\eeq
Also, we notice that $\widehat H_2(x, v^{(i)}, v^{(j)}) = \widehat H^T_2(x, v^{(j)}, v^{(i)})$, where $T$ denotes the transpose of a matrix. In view of Prop.\ \ref{intert}, we can write the causal inverse $\bfq = Q_{\widehat g} \bold{Id}$ at $q_0$ below.  
The computation of the symbol is straightforward but very lengthy, so we compute each term in a separate subsection.   

\subsection{Computation of $\sigma(\mch_1)$}
In this subsection, we write $\mch_1 = \mch$. For $\alpha \neq \beta$, we have  
\beqq\label{eqh3}
\begin{split}
\mch_{\alpha\beta} = &\sum_{(i, j, k, l)\in \sigma(4)}\widehat H_{3, \alpha\beta}(\bfq \widehat H_2(\vec v^{(i)}, \vec v^{(j)}), \vec v^{(k)}, \vec v^{(l)}) \\
= &4\sum_{(i, j, k, l)\in \sigma(4)}-h^{aa'}h^{bb'}Q_{\widehat g}(-h^{\la\mu}F^{(i)}_{a\la} F^{(j)}_{b\mu} + \frac{1}{4}h_{ab}h^{\la\gamma}h^{\mu\delta}F^{(i)}_{\gamma\delta}F^{(j)}_{\la\mu})F^{(k)}_{\alpha a'} F^{(l)}_{\beta b'} \\
&+ 4\sum_{(i, j, k, l)\in \sigma(4)} \frac{1}{4}Q_{\widehat g}(-h^{\la\mu}F^{(i)}_{\alpha\la} F^{(j)}_{\beta\mu})h^{aa'}h^{bb'}F^{(k)}_{ab}F^{(l)}_{a'b'}\\
%= &\sum_{(i, j, k, l) \in \sigma(4)}  -h^{aa}h^{bb}Q_{\widehat g}[-h^{\la\la}F^{(i)}_{a\la}F^{(j)}_{b\la}]F^{(k)}_{\alpha a} F^{(l)}_{\beta b}  - \sum_{(i, j, k, l) \in \sigma(4)} h^{aa}h^{bb} Q_{\widehat g}[\frac{1}{4} h_{ab} h^{\la\la}h^{\mu\mu}F^{(i)}_{\la\mu}F_{\la\mu}^{(j)}]F^{(k)}_{\alpha a} F^{(l)}_{\beta b} \\
%&- \sum_{(i, j, k, l) \in \sigma(4)}  \frac{1}{4} Q_{\widehat g}(h^{\la\la} F^{(i)}_{\alpha\la}F^{(j)}_{\beta\la})h^{aa}h^{bb}F^{(k)}_{ab}F^{(l)}_{ab} \\
= &  4\big(\mci_{1, \alpha\beta} + \mci_{2, \alpha\beta} + \mci_{3, \alpha\beta}\big), 
\end{split}
\eeqq
where $F^{(i)} = d\dot \phi^{(i)} \text{ and }$ 
\beq 
\begin{split}
\mci_{1, \alpha\beta} = &\sum_{a, b, \la = 0}^3 \sum_{(i, j, k, l) \in \sigma(4)}   h^{aa}h^{bb}h^{\la\la}Q_{\widehat g}[F^{(i)}_{a\la}F^{(j)}_{b\la}]F^{(k)}_{\alpha a} F^{(l)}_{\beta b} \\
 \mci_{2, \alpha\beta} = & -\sum_{a, \la, \mu = 0}^3  \sum_{(i, j, k, l) \in \sigma(4)} \frac{1}{4}  h^{\la\la}h^{\mu\mu}h^{aa} Q_{\widehat g}[ F^{(i)}_{\la\mu}F_{\la\mu}^{(j)}]F^{(k)}_{\alpha a} F^{(l)}_{\beta a} \\
\mci_{3, \alpha\beta} = &- \sum_{a, b, \la = 0}^3 \sum_{(i, j, k, l) \in \sigma(4)}  \frac{1}{4}h^{\la\la} h^{aa}h^{bb} Q_{\widehat g}(F^{(i)}_{\alpha\la}F^{(j)}_{\beta\la})F^{(k)}_{ab}F^{(l)}_{ab}.
\end{split}
\eeq
To see the structure of the terms, we find that 
\beqq\label{compfor}
\begin{gathered}
\mci_1 = \sum_{(i, j, k, l) \in \sigma(4)}  F^{(k)}HQ_{\widehat g}(F^{(i)}HF^{(j)})HF^{(l)}, \\
\mci_2 =   \frac{1}{4} \sum_{(i, j, k, l) \in \sigma(4)} Q_{\widehat g}(\text{Tr}(F^{(i)}F^{(j)}))F^{(k)}HF^{(l)}, \\
\mci_3 =   \frac{1}{4} \sum_{(i, j, k, l) \in \sigma(4)}  Q_{\widehat g}(F^{(i)}HF^{(j)})\text{Tr}(F^{(k)}F^{(l)}).
\end{gathered}
\eeqq
In these formulas, $Q_{\widehat g}$ applies to each element of the matrix. Also, the sign of $\mci_2, \mci_3$ are positive because we used the anti-symmetry of $F$. We emphasize that we can only use these formula to express the off-diagonal elements.  
{We know (see e.g.\ Lemma 4.1 \cite{WZ}) that if $\dot \phi^{(i)}\in I^\mu(\La)$, then the principal symbol  of $F^{(i)}_{\alpha\beta}$ (in $I^{\mu+1}(\La)$) are given by 
\beq
\sigma(F^{(i)}_{\alpha\beta})(x, \zeta) = \imath\zeta_\alpha \sigma(\dot \phi^{(i)}_\beta)(x, \zeta) - \imath\zeta_\beta \sigma(\dot \phi^{(i)}_{\alpha})(x, \zeta), \ \ (x, \zeta) \in \La.
\eeq
So at the intersection point $q_0$, we get $\sigma(F^{(i)}_{\alpha\beta})(q_0, \zeta^{(i)}) = \imath\zeta^{(i)}_\alpha A^{(i)}_\beta  - \imath\zeta^{(i)}_\beta A^{(i)}_{\alpha}$. \\
%By substituting these formulas in (\ref{eqh3}) and (\ref{compfor}), and using the expressions of the principal symbols in Section \ref{intera}, we see   that
% $\sigma(\mch)(q_0, \zeta)= \mcp(\vec \zeta,\vec A)$ is a real analytic function of $\vec\zeta$  and $\vec A$. This proves claim (ii)
%} 

%\subsubsection{Computation for the first set of vectors}
We start with the set of vectors considered in  Prop.\ \ref{symb}
\beq
\begin{gathered}
\zeta^{(1)}=(1,0,1,0),\ \ \zeta^{(2)}=(1,0,0,1), \ \
\zeta^{(3)}=  (-3/4, -3/4, 0,0),\ \ \zeta^{(4)} =  (1,-1,0,0)
\end{gathered}
\eeq
so that the sum $\zeta = (9/4, -3/4, 1, 1)$ is light-like. The vectors ${}^1\vec A$ are taken as 
\beq
\begin{gathered}
{}^1A^{(1)}= -(0,0,0,\imath),\ \ {}^1A^{(2)}= -(0,0,\imath,0), \ \
{}^1A^{(3)}=  -(0, 0, 0,\imath),\ \ {}^1A^{(4)} =  -(0,0,0,\imath).
\end{gathered}
\eeq
We calculate the principal symbol $\sigma(\mch({}^1\vec A))$   in detail for this set of vectors. 
We start with some preparations. For the set of vectors $\zeta^{(i)}, i = 1, 2, 3, 4$, we first compute  $\mcg(i, j) = |\zeta^{(i)} + \zeta^{(j)}|_h^2 = 2h(\zeta^{(i)}, \zeta^{(j)})$:  
\beq
\begin{gathered}
\mcg(1, 2) = -2, \ \ \mcg(1, 3) = 3/2, \ \ \mcg(1, 4) =  -2, \\
\mcg(2, 3) = 3/2, \ \ \mcg(2, 4) = -2, \ \ \mcg(3, 4) = 3.
\end{gathered}
\eeq
Next, we compute the trace terms $\mcw(i, j) = \text{Tr}(F^{(i)}F^{(j)})$: 
\beq
\begin{gathered}
\mcw(1, 2) = -2, \ \ \mcw(1, 3) = 3/2, \ \ \mcw(1, 4) =  -2, \\
\mcw(2, 3) = 0, \ \ \mcw(2, 4) = 0, \ \ \mcw(3, 4) = 3.
\end{gathered}
\eeq

Now we compute the terms $\mch_{\alpha\beta}$ using the formulas. Although there are many summation terms due to the permutations, eventually the non-trivial terms are few thanks to the sparseness of $F^\bullet$. We remark that one can also compute using the matrix representation \eqref{compfor}. We split the rest into five parts, each dealing with one term.\\

 %%%%%%%%%%%%%%%%%%%%%%%
\textbf{(a) $\mch_{02}$:} We begin with
\beq
\begin{split}
 \mci_{1, 02} = &\sum_{(i, j, k, l) \in \sigma(4)}  \sum_{a, b, \la = 0}^3 h^{aa}h^{bb}h^{\la\la}Q_{\widehat g}[F^{(i)}_{a\la}F^{(j)}_{b\la}]F^{(k)}_{0 a} F^{(l)}_{2 b}  \\
 = & \sum_{(i, j, k, l) \in \sigma(4)} \sum_{\la = 0}^3 h^{\la\la}Q_{\widehat g}[F^{(i)}_{2\la}F^{(j)}_{3\la}]F^{(k)}_{0 2} F^{(l)}_{23}+ \sum_{(i, j, k, l) \in \sigma(4)}  \sum_{\la, b = 0}^3 h^{bb}h^{\la\la}Q_{\widehat g}[F^{(i)}_{3\la}F^{(j)}_{b\la}]F^{(k)}_{03} F^{(l)}_{2 b}\\
 = & \sum_{(i, j, k, l) \in \sigma(4)} \sum_{\la = 0}^3  -h^{\la\la}Q_{\widehat g}[F^{(i)}_{3\la}F^{(j)}_{0\la}]F^{(k)}_{03} F^{(l)}_{2 0} + \sum_{(i, j, k, l) \in \sigma(4)} \sum_{\la = 0}^3  h^{\la\la}Q_{\widehat g}[F^{(i)}_{3\la}F^{(j)}_{3\la}]F^{(k)}_{03} F^{(l)}_{2 3}\\
  = &  \sum_{(i, j, k, l) \in \sigma(4)}  -Q_{\widehat g}[F^{(i)}_{30}F^{(j)}_{30}]F^{(k)}_{03} F^{(l)}_{2 3} +  \sum_{(i, j, k, l) \in \sigma(4)}   Q_{\widehat g}[F^{(i)}_{31}F^{(j)}_{31}]F^{(k)}_{03} F^{(l)}_{2 3}\\
   = &  -2Q_{\widehat g}[F^{(1)}_{30}F^{(3)}_{30}]F^{(4)}_{03} F^{(2)}_{2 3}  -2Q_{\widehat g}[F^{(1)}_{30}F^{(4)}_{30}]F^{(3)}_{03} F^{(2)}_{2 3} -2Q_{\widehat g}[F^{(3)}_{30}F^{(4)}_{30}]F^{(1)}_{03} F^{(2)}_{2 3}\\
   &  - 2Q_{\widehat g}[F^{(3)}_{31}F^{(4)}_{31}]F^{(1)}_{03} F^{(2)}_{2 3}.
\end{split}
\eeq
The symbols can be calculated using the expressions in Section \ref{singu}.  For simplicity, we shall ignore the $c_\pi = (2\pi)^{-3}$ factor  as well as the multiple of $4$ in the computations. We add them in the final answer. For example, we shall compute using
\beq
\sigma(Q_{\widehat g}[F^{(i)}_{\bullet}F^{(j)}_{\bullet}]F^{(k)}_{\bullet} F^{(l)}_{\bullet}) =  \frac{1}{|\zeta^{(i)} + \zeta^{(j)}|^2_h}\sigma(F^{(i)}_{\bullet})\sigma(F^{(j)}_{\bullet})\sigma(F^{(k)}_{\bullet})\sigma(F^{(l)}_{\bullet}),
\eeq
where $\bullet$ stands for a generic index. Hereafter the principal symbols $\sigma(F^{(i)})$ are always evaluated at $(q_0, \zeta^{(i)})$ which shall be omitted in the notations. Now we use the values of $F^\bullet$ and $\mcg(\cdot, \cdot)$ to get
\beq
\begin{split}
\sigma(\mci_{1, 02}) &= -2 (2/3) (-1) (3/4) 1 (-1) -2 (-1/2)(-1)(-1)(-3/4) (-1) -2(1/3)(3/4)(-1)1(-1) \\
& + 2 (1/3)(3/4)1\cdot 1 (-1) = -5/4.
\end{split}
\eeq

Next, consider 
\beq
\begin{split}
 \mci_{2, 02} =    - \sum_{(i, j, k, l) \in \sigma(4)} \sum_{\la, \mu = 0}^3\frac{1}{4}  h^{\la\la}h^{\mu\mu} Q_{\widehat g}[ F^{(i)}_{\la\mu}F_{\la\mu}^{(j)}]F^{(k)}_{0 3} F^{(l)}_{2 3} 
\end{split}
\eeq
When we compute the symbol, we will have $\mcw(i, j)$ which is zero for $(i, j) = (2, 3)$ and $(2, 4)$. Now we determine the non-trivial terms for all permutations of $(i, j, k, l)$. First, $l$ has to be $1$ or $2$. If $l = 1$, then $k$ must be $3$ or $4$ so that at least one of $i, j$ is $2$. But then the corresponding $\mcw$ term become zero. So it suffices to take $l = 2$. Then we have
\beq
\begin{split}
 \mci_{2, 02} =   &  \sum_{\la, \mu = 0}^3 \bigg( -  \frac{1}{2}  h^{\la\la}h^{\mu\mu} Q_{\widehat g}[ F^{(3)}_{\la\mu}F_{\la\mu}^{(4)}]F^{(1)}_{0 3} F^{(2)}_{2 3} 
 -  \frac{1}{2}    h^{\la\la}h^{\mu\mu} Q_{\widehat g}[ F^{(1)}_{\la\mu}F_{\la\mu}^{(4)}]F^{(3)}_{0 3} F^{(2)}_{2 3}\\
 & -  \frac{1}{2}    h^{\la\la}h^{\mu\mu} Q_{\widehat g}[ F^{(1)}_{\la\mu}F_{\la\mu}^{(3)}]F^{(4)}_{0 3} F^{(2)}_{2 3}\bigg).
\end{split}
\eeq
Then we compute the symbol and find that 
\beq
\sigma(\mci_{2, 02}) = -\ha [\frac{3}{3} 1 (-1) + \frac{-2}{-1} (-3/4) (-2) + \frac{3/2}{3/4} 1 (-1)] = 5/8.
\eeq
Finally, consider 
\beq
\begin{split}
\mci_{3, 02} = &- \sum_{(i, j, k, l) \in \sigma(4)}   \sum_{a, b = 0}^3 \frac{1}{4} h^{aa}h^{bb} Q_{\widehat g}(F^{(i)}_{03}F^{(j)}_{23})F^{(k)}_{ab}F^{(l)}_{ab}
\end{split}
\eeq
This is similar to $\mci_{2, 02}$. When $j = 1$, the symbol of the terms vanishes. For $j = 2,$ we obtain  
\beq
\begin{split}
\mci_{3, 02} = & \sum_{a, b = 0}^3\bigg( - \frac{1}{2}  h^{aa}h^{bb} Q_{\widehat g}(F^{(1)}_{03}F^{(2)}_{23})F^{(3)}_{ab}F^{(4)}_{ab}- \frac{1}{2}  h^{aa}h^{bb} Q_{\widehat g}(F^{(3)}_{03}F^{(2)}_{23})F^{(1)}_{ab}F^{(4)}_{ab} \\
&  - \frac{1}{2}   h^{aa}h^{bb} Q_{\widehat g}(F^{(4)}_{03}F^{(2)}_{23})F^{(1)}_{ab}F^{(3)}_{ab}\bigg)
\end{split}
\eeq
Then we can find the symbol is $\sigma(\mci_{3, 02}) = -5/8$. Finally, we have $\sigma(\mch_{02}) = 4c_\pi(-5/4).$\\

%%%%%%%%%%%%%%%%%%%%%%
\textbf{(b) $\mch_{12}$:} We start with 
\beq
\begin{split}
\mci_{1, 12} = &\sum_{(i, j, k, l) \in \sigma(4)} \sum_{b, \la = 0}^3 h^{bb}h^{\la\la}Q_{\widehat g}[F^{(i)}_{3\la}F^{(j)}_{b\la}]F^{(k)}_{1 3} F^{(l)}_{2 b}  \\
 =  &\sum_{(i, j, k, l) \in \sigma(4)}  \sum_{\la = 0}^3-h^{\la\la}Q_{\widehat g}[F^{(i)}_{3\la}F^{(j)}_{0\la}]F^{(k)}_{1 3} F^{(2)}_{2 0}  +   \sum_{(i, j, k, l) \in \sigma(4)} \sum_{\la = 0}^3 h^{\la\la}Q_{\widehat g}[F^{(i)}_{3\la}F^{(j)}_{3\la}]F^{(k)}_{1 3} F^{(l)}_{2 3} \\
 =  &  \sum_{(i, j, k, l) \in \sigma(4)} -Q_{\widehat g}[F^{(i)}_{30}F^{(j)}_{30}]F^{(k)}_{1 3} F^{(l)}_{2 3}.
\end{split}
\eeq
We know that $l$ must be $1, 2$ and $k$ must be $3, 4$. So we just need to consider these permutations.  
\beq
\begin{split}
\sigma(\mci_{1, 12}) &= -2\sigma(Q_{\widehat g}[F^{(1)}_{30}F^{(4)}_{30}]F^{(3)}_{1 3} F^{(2)}_{2 3})-2\sigma(Q_{\widehat g}[F^{(1)}_{30}F^{(3)}_{30}]F^{(4)}_{1 3} F^{(2)}_{2 3}) \\
%&-2\sigma(Q_{\widehat g}[F^{(2)}_{30}F^{(4)}_{30}]F^{(3)}_{1 3} F^{(1)}_{2 3})-2\sigma(Q_{\widehat g}[F^{(2)}_{30}F^{(3)}_{30}]F^{(4)}_{1 3} F^{(1)}_{2 3})\\
& = -2 (-1/2) (-1) (-1) (-3/4) (-1) - 2 (2/3) (-1) (3/4) (-1) (-1) = 7/4.
\end{split}
\eeq
Next, we consider 
\beq
\begin{split}
\mci_{2, 12} = & - \sum_{(i, j, k, l) \in \sigma(4)} \sum_{\la, \mu = 0}^3\frac{1}{4}  h^{\la\la}h^{\mu\mu} Q_{\widehat g}[ F^{(i)}_{\la\mu}F_{\la\mu}^{(j)}]F^{(k)}_{1 3} F^{(l)}_{2 3}.
\end{split}
\eeq
We know that $l$ must be $2$ because if $l = 1$, we always have $i, j$ are $2, 3$ or $2, 4$ and the corresponding $\mcw$ vanish. So we have
\beq
\begin{split}
\sigma(\mci_{2, 12}) &= \sum_{\la, \mu = 0}^3\bigg( - \frac{1}{2}   \sigma( h^{\la\la}h^{\mu\mu} Q_{\widehat g}[ F^{(1)}_{\la\mu}F_{\la\mu}^{(4)}]F^{(3)}_{1 3} F^{(2)}_{2 3})- \frac{1}{2}  \sigma( h^{\la\la}h^{\mu\mu} Q_{\widehat g}[ F^{(1)}_{\la\mu}F_{\la\mu}^{(3)}]F^{(4)}_{1 3} F^{(2)}_{2 3})\bigg) \\
 &= -\ha \cdot \frac{-2}{-2} (-3/4) (-1) - \ha \cdot \frac{3/2}{3/2} (-1) (-1) = -7/2.
\end{split}
\eeq
Finally, we find that 
\[
\mci_{3, 12} =  - \sum_{(i, j, k, l) \in \sigma(4)}  \sum_{a, b = 0}^3 \frac{1}{4} h^{aa}h^{bb} Q_{\widehat g}(F^{(i)}_{1 3}F^{(j)}_{2 3})F^{(k)}_{ab}F^{(l)}_{ab}.
\]
By similar consideration, we get
\beq
\begin{split}
\sigma(\mci_{3, 12}) &=  \sum_{a, b = 0}^3\bigg( - \frac{1}{2}\sigma(  h^{aa}h^{bb} Q_{\widehat g}(F^{(3)}_{1 3}F^{(2)}_{2 3})F^{(1)}_{ab}F^{(4)}_{ab})- \frac{1}{2} \sigma(  h^{aa}h^{bb} Q_{\widehat g}(F^{(4)}_{1 3}F^{(2)}_{2 3})F^{(1)}_{ab}F^{(3)}_{ab})\bigg)\\
 &= -\ha \cdot \frac{-2}{3/2} (-3/4) (-1) - \ha \cdot \frac{3/2}{-2} (-1)(-1) = 7/2
\end{split}
\eeq
Summing up the symbols, we get $\sigma(\mch_{12}) = 4c_\pi(7/4).$\\

%%%%%%%%%%%%%%%%%%%%%%%
\textbf{(c) $\mch_{03}$:} We start with 
\beq
\begin{split}
 \mci_{2, 03} &=    - \sum_{(i, j, k, l) \in \sigma(4)}  \sum_{\la, \mu, a = 0}^3\frac{1}{4}  h^{\la\la}h^{\mu\mu}h^{aa} Q_{\widehat g}[ F^{(i)}_{\la\mu}F_{\la\mu}^{(j)}]F^{(k)}_{0 a} F^{(l)}_{3 a}  \\
 & =  - \sum_{(i, j, k, l) \in \sigma(4)}  \sum_{\la, \mu = 0}^3 \frac{1}{4}  h^{\la\la}h^{\mu\mu}  Q_{\widehat g}[ F^{(i)}_{\la\mu}F_{\la\mu}^{(j)}]F^{(k)}_{0 2} F^{(l)}_{3 2}  =   \sum_{\la, \mu = 0}^3 -\frac{1}{4}  2 h^{\la\la}h^{\mu\mu}  Q_{\widehat g}[ F^{(3)}_{\la\mu}F_{\la\mu}^{(4)}]F^{(2)}_{0 2} F^{(1)}_{3 2}.
\end{split}
\eeq
So the symbol is $\sigma(\mci_{2, 03}) = -\ha \cdot \frac{3}{3} 1 (-1) =   1/2.$ 
Similarly, we have 
\[
\mci_{3, 03} =   \sum_{\la, \mu = 0}^3-\frac{1}{4}  2 h^{\la\la}h^{\mu\mu}   F^{(3)}_{\la\mu}F_{\la\mu}^{(4)} Q_{\widehat g}[F^{(2)}_{0 2} F^{(1)}_{3 2}],
\]
 and the symbol is $\sigma(\mci_{3, 03}) = -\ha \cdot \frac{3}{-2} 1 (-1) = -3/4.$  Finally we calculate
\beq
\begin{split}
 &\mci_{1, 03} = \sum_{(i, j, k, l) \in \sigma(4)}   \sum_{a, b, \la = 0}^3 h^{aa}h^{bb}h^{\la\la}Q_{\widehat g}[F^{(i)}_{a\la}F^{(j)}_{b\la}]F^{(k)}_{0a} F^{(l)}_{3b}\\
  = & \sum_{(i, j, k, l) \in \sigma(4)} \sum_{b, \la = 0}^3  h^{bb}h^{\la\la}Q_{\widehat g}[F^{(i)}_{2\la}F^{(j)}_{b\la}]F^{(k)}_{02} F^{(l)}_{3b}  +  \sum_{(i, j, k, l) \in \sigma(4)} \sum_{b, \la = 0}^3   h^{bb}h^{\la\la}Q_{\widehat g}[F^{(i)}_{3\la}F^{(j)}_{b\la}]F^{(k)}_{03} F^{(l)}_{3b}   \\
   = & \sum_{(i, j, k, l) \in \sigma(4)}  \sum_{b = 0}^3 h^{bb} Q_{\widehat g}[F^{(i)}_{23}F^{(j)}_{b3}]F^{(k)}_{02} F^{(l)}_{3b}   -  \sum_{(i, j, k, l) \in \sigma(4)} \sum_{\la = 0}^3  h^{\la\la}Q_{\widehat g}[F^{(i)}_{3\la}F^{(j)}_{0\la}]F^{(k)}_{03} F^{(l)}_{30} \\
   & +   \sum_{(i, j, k, l) \in \sigma(4)}  \sum_{\la = 0}^3   h^{\la\la}Q_{\widehat g}[F^{(i)}_{3\la}F^{(j)}_{1\la}]F^{(k)}_{03} F^{(l)}_{31} +   \sum_{(i, j, k, l) \in \sigma(4)}\sum_{\la = 0}^3    h^{\la\la}Q_{\widehat g}[F^{(i)}_{3\la}F^{(j)}_{2\la}]F^{(k)}_{03} F^{(l)}_{32} \\
   %=  & -\sum_{(i, j, k, l) \in \sigma(4)}    Q_{\widehat g}[F^{(i)}_{23}F^{(j)}_{03}]F^{(k)}_{02} F^{(l)}_{30} +  \sum_{(i, j, k, l) \in \sigma(4)}   Q_{\widehat g}[F^{(i)}_{23}F^{(j)}_{13}]F^{(k)}_{02} F^{(l)}_{31}\\
   %&-  \sum_{(i, j, k, l) \in \sigma(4)}  Q_{\widehat g}[F^{(i)}_{32}F^{(j)}_{02}]F^{(k)}_{03} F^{(l)}_{30}   -  \sum_{(i, j, k, l) \in \sigma(4)}  Q_{\widehat g}[F^{(i)}_{30}F^{(j)}_{20}]F^{(k)}_{03} F^{(l)}_{32}\\
%    =& - Q_{\widehat g}[F^{(1)}_{23}F^{(3)}_{03}]F^{(2)}_{02} F^{(4)}_{30} - Q_{\widehat g}[F^{(1)}_{23}F^{(4)}_{03}]F^{(2)}_{02} F^{(3)}_{30} + Q_{\widehat g}[F^{(1)}_{23}F^{(3)}_{13}]F^{(2)}_{02} F^{(4)}_{31} \\
%    &+ Q_{\widehat g}[F^{(1)}_{23}F^{(4)}_{13}]F^{(2)}_{02} F^{(3)}_{31}  - Q_{\widehat g}[F^{(1)}_{32}F^{(2)}_{02}]F^{(3)}_{03} F^{(4)}_{30}- Q_{\widehat g}[F^{(1)}_{32}F^{(2)}_{02}]F^{(4)}_{03} F^{(3)}_{30} \\
%    &-Q_{\widehat g}[F^{(3)}_{30}F^{2)}_{20}]F^{(4)}_{03} F^{(1)}_{32}-Q_{\widehat g}[F^{(4)}_{30}F^{2)}_{20}]F^{(3)}_{03} F^{(1)}_{32}
=& -  Q_{\widehat g}[F^{(1)}_{23}F^{(3)}_{03}]F^{(2)}_{02} F^{(4)}_{30} - Q_{\widehat g}[F^{(1)}_{23}F^{(4)}_{03}]F^{(2)}_{02} F^{(3)}_{30} + Q_{\widehat g}[F^{(1)}_{23}F^{(3)}_{13}]F^{(2)}_{02} F^{(4)}_{31} \\
   + &Q_{\widehat g}[F^{(1)}_{23}F^{(4)}_{13}]F^{(2)}_{02} F^{(3)}_{31}  - 2Q_{\widehat g}[F^{(1)}_{32}F^{(2)}_{02}]F^{(3)}_{03} F^{(4)}_{30}  -Q_{\widehat g}[F^{(3)}_{30}F^{(2)}_{20}]F^{(4)}_{03} F^{(1)}_{32} -Q_{\widehat g}[F^{(4)}_{30}F^{(2)}_{20}]F^{(3)}_{03} F^{(1)}_{32}
\end{split}
\eeq
We compute the symbol as
\beq
\begin{split}
\sigma(\mci_{1, 03}) &= -(2/3)1(-3/4)1(-1) - (-\ha)1\cdot 1\cdot 1(3/4) + (2/3) 1 (-3/4) 1 \cdot 1\\
&+ (-\ha) 1 (-1) 1 (3/4) - 2(-\ha) (-1) 1 (-3/4) (-1) - (2/3) (3/4) (-1) 1\cdot (-1)\\
& - (-\ha)(-1)(-1)(-3/4)(-1) = -9/8.
\end{split}
\eeq
Summing up the symbols, we get $\sigma(\mch_{03}) = 4c_\pi(-11/8).$\\

%%%%%%%%%%%%%%%%%%%%%
\textbf{(d) $\mch_{13}$:} We compute
\beq
\begin{split}
 &\mci_{1, 13} = \sum_{(i, j, k, l) \in \sigma(4)}  \sum_{b, \la = 0}^3   h^{bb}h^{\la\la}Q_{\widehat g}[F^{(i)}_{3\la}F^{(j)}_{b\la}]F^{(k)}_{1 3} F^{(l)}_{3 b}  \\
  =  &\sum_{(i, j, k, l) \in \sigma(4)} \sum_{\la = 0}^3   \bigg( - h^{\la\la}Q_{\widehat g}[F^{(i)}_{3\la}F^{(j)}_{0\la}]F^{(k)}_{1 3} F^{(l)}_{3 0} +  h^{\la\la}Q_{\widehat g}[F^{(i)}_{3\la}F^{(j)}_{1\la}]F^{(k)}_{1 3} F^{(l)}_{3 1}     +   h^{\la\la}Q_{\widehat g}[F^{(i)}_{3\la}F^{(j)}_{2\la}]F^{(k)}_{1 3} F^{(l)}_{3 2}\bigg)\\
     =  &-\sum_{(i, j, k, l) \in \sigma(4)}   Q_{\widehat g}[F^{(i)}_{32}F^{(j)}_{02}]F^{(k)}_{1 3} F^{(l)}_{3 0}  -  \sum_{(i, j, k, l) \in \sigma(4)}   Q_{\widehat g}[F^{(i)}_{3 0}F^{(j)}_{20}]F^{(k)}_{1 3} F^{(l)}_{3 2}\\
      = & -   Q_{\widehat g}[F^{(1)}_{32}F_{02}^{(2)}]F^{(3)}_{1 3} F^{(4)}_{3 0} -   Q_{\widehat g}[F^{(3)}_{3 0}F^{(2)}_{20}]F^{(4)}_{1 3} F^{(1)}_{3 2} -   Q_{\widehat g}[F^{(1)}_{32}F_{02}^{(2)}]F^{(4)}_{1 3} F^{(3)}_{3 0}  -   Q_{\widehat g}[F^{(4)}_{3 0}F^{(2)}_{20}]F^{(3)}_{1 3} F^{(1)}_{3 2}.
\end{split}
\eeq
We can compute the symbol as
\beq
\begin{split}
\sigma(\mci_{1, 13}) & = -(-1/2) (-1)1(-3/4)(-1) -(2/3) (3/4) (-1) (-1) (-1) \\
&- (-1/2) (-1) 1 (-1) (3/4) - (-1/2) (-1) (-1) (-3/4) (-1)  = 7/8.
\end{split}
\eeq

For $\mci_{2, 13}$, we notice that $\mci_{2, 13} =    - \sum_{(i, j, k, l) \in \sigma(4)} \sum_{\la, \mu = 0}^3  \frac{1}{4} h^{\mu\mu}h^{\la\la} Q_{\widehat g}[ F^{(i)}_{\la\mu}F_{\la\mu}^{(j)}]F^{(k)}_{1 3} F^{(l)}_{3 3} =0.$ Similarly, $\mci_{3, 13} = 0$. Thus $\sigma(\mch_{13}) = 4c_\pi \sigma(\mci_{1, 13}) = 4c_\pi (7/8).$\\

%%%%%%%%%%%%%%%%%%%
\textbf{(e) $\mch_{23}$:} We begin with
\beq
\begin{split}
 \mci_{1, 23}& = \sum_{(i, j, k, l) \in \sigma(4)}\sum_{a, b, \la = 0}^3    h^{aa}h^{bb}h^{\la\la}Q_{\widehat g}[F^{(i)}_{a\la}F^{(j)}_{b\la}]F^{(k)}_{2 a} F^{(l)}_{3 b}  \\
 & =  \sum_{(i, j, k, l) \in \sigma(4)} \sum_{b = 0}^3  -h^{bb} Q_{\widehat g}[F^{(i)}_{03}F^{(j)}_{b3}]F^{(k)}_{2 0} F^{(l)}_{3 b} +  \sum_{(i, j, k, l) \in \sigma(4)} \sum_{b, \la = 0}^3   h^{bb}h^{\la\la}Q_{\widehat g}[F^{(i)}_{ 3\la}F^{(j)}_{b\la}]F^{(k)}_{2 3} F^{(l)}_{3 b}.
\end{split}
\eeq
There are many terms in the summations so we deal the two summations separately.
\beq
\begin{split}
 \mca_1\doteq &\sum_{(i, j, k, l) \in \sigma(4)}  \sum_{b = 0}^3  -h^{bb} Q_{\widehat g}[F^{(i)}_{03}F^{(j)}_{b3}]F^{(k)}_{2 0} F^{(l)}_{3 b} \\
 = &  \sum_{(i, j, k, l) \in \sigma(4)} \bigg(  Q_{\widehat g}[F^{(i)}_{03}F^{(j)}_{03}]F^{(k)}_{2 0} F^{(l)}_{30} - Q_{\widehat g}[F^{(i)}_{03}F^{(j)}_{13}]F^{(k)}_{2 0} F^{(l)}_{3 1} - Q_{\widehat g}[F^{(i)}_{03}F^{(j)}_{23}]F^{(k)}_{2 0}  F^{(l)}_{3 2}\bigg)\\
  = &  2\bigg(  Q_{\widehat g}[F^{(3)}_{03}F^{(4)}_{03}]F^{(2)}_{2 0} F^{(1)}_{30} + Q_{\widehat g}[F^{(1)}_{03}F^{(4)}_{03}]F^{(2)}_{2 0} F^{(3)}_{30} + Q_{\widehat g}[F^{(1)}_{03}F^{(3)}_{03}]F^{(2)}_{2 0} F^{(4)}_{30} \bigg)\\ 
  & +\bigg(- Q_{\widehat g}[F^{(1)}_{03}F^{(3)}_{13}]F^{(2)}_{2 0} F^{(4)}_{3 1} - Q_{\widehat g}[F^{(1)}_{03}F^{(4)}_{13}]F^{(2)}_{2 0} F^{(3)}_{3 1}\bigg). 
\end{split}
\eeq
Then we find 
\beq
\begin{split}
\sigma(\mca_1) &= 2\bigg( (1/3) (-3/4)1 (-1)(-1) + (-\ha) 1 \cdot 1(-1) (3/4) +   (2/3)1 (-3/4) (-1) (-1)  \bigg)\\
& + \bigg(- (2/3) 1 (-3/4) (-1) 1 - (-\ha ) 1 (-1) (-1) (3/4)\bigg) = -7/8.
\end{split}
\eeq

For the other summation, we have
\beq
\begin{split}
\mca_2 & \doteq \sum_{(i, j, k, l) \in \sigma(4)} \sum_{b, \la = 0}^3   h^{bb}h^{\la\la}Q_{\widehat g}[F^{(i)}_{ 3\la}F^{(j)}_{b\la}]F^{(k)}_{2 3} F^{(l)}_{3 b}\\
 = & \sum_{(i, j, k, l) \in \sigma(4)} \sum_{\la = 0}^3 \bigg(  -h^{\la\la}Q_{\widehat g}[F^{(i)}_{ 3\la}F^{(j)}_{0\la}]F^{(k)}_{2 3} F^{(l)}_{3 0}  + h^{\la\la}Q_{\widehat g}[F^{(i)}_{ 3\la}F^{(j)}_{1\la}]F^{(k)}_{2 3} F^{(l)}_{3 1} +  h^{\la\la}Q_{\widehat g}[F^{(i)}_{ 3\la}F^{(j)}_{2\la}]F^{(k)}_{2 3} F^{(l)}_{3 2} \bigg)\\
=& \sum_{(i, j, k, l) \in \sigma(4)}\bigg(  - Q_{\widehat g}[F^{(i)}_{ 32}F^{(j)}_{02}]F^{(k)}_{2 3} F^{(l)}_{3 0}  - Q_{\widehat g}[F^{(i)}_{ 30}F^{(j)}_{20}]F^{(k)}_{2 3} F^{(l)}_{3 2}  \bigg) = 0.
\end{split}
\eeq
Therefore, $\sigma(\mci_{1, 23}) = 1/8.$ 
Next, we compute
\beq
\begin{split}
 &\mci_{2, 23}  =     \sum_{(i, j, k, l) \in \sigma(4)} \sum_{\la, \mu = 0}^3  \frac{1}{4}  h^{\la\la}h^{\mu\mu}  Q_{\widehat g}[ F^{(i)}_{\la\mu}F_{\la\mu}^{(j)}]F^{(k)}_{2 0} F^{(l)}_{3 0} \\
  &=   \sum_{\la, \mu = 0}^3 \bigg(\frac{1}{2}  h^{\la\la}h^{\mu\mu}  Q_{\widehat g}[ F^{(3)}_{\la\mu}F_{\la\mu}^{(4)}]F^{(2)}_{2 0} F^{(1)}_{3 0} +  \frac{1}{2}  h^{\la\la}h^{\mu\mu}  Q_{\widehat g}[ F^{(1)}_{\la\mu}F_{\la\mu}^{(4)}]F^{(2)}_{2 0} F^{(3)}_{3 0} + \frac{1}{2}  h^{\la\la}h^{\mu\mu}  Q_{\widehat g}[ F^{(1)}_{\la\mu}F_{\la\mu}^{(3)}]F^{(2)}_{2 0} F^{(4)}_{3 0}\bigg).
\end{split}
\eeq
Therefore, the symbol is
\beq
\begin{split}
\sigma(\mci_{2, 23}) = \ha \cdot \frac{3}{3} (-1) (-1) + \ha\cdot \frac{-2}{-2} (-1)(3/4) + \ha \cdot\frac{3/2}{3/2} (-1)(-1) = 5/8.
\end{split}
\eeq

Finally, we consider
\beq
\begin{split}
&\mci_{3, 23} =   \sum_{(i, j, k, l) \in \sigma(4)} \sum_{a, b = 0}^3  \frac{1}{4}  h^{aa}h^{bb} Q_{\widehat g}(F^{(i)}_{20}F^{(j)}_{30})F^{(k)}_{ab}F^{(l)}_{ab}\\
&=  \sum_{a, b  = 0}^3 \bigg( \frac{1}{2}  h^{aa}h^{bb}   F^{(3)}_{ab}F_{ab}^{(4)} Q_{\widehat g}[F^{(2)}_{2 0} F^{(1)}_{3 0}] +  \frac{1}{2}  h^{aa}h^{bb}    F^{(1)}_{ab}F_{ab}^{(4)} Q_{\widehat g}[F^{(2)}_{2 0} F^{(3)}_{3 0}] + \frac{1}{2}  h^{aa}h^{bb}   F^{(1)}_{ab}F_{ab}^{(3)}Q_{\widehat g}[F^{(2)}_{2 0} F^{(4)}_{3 0}]\bigg).
\end{split}
\eeq
The symbol can be found as
\beq
\begin{split}
\sigma(\mci_{3, 23}) =  \ha \cdot \frac{3}{-2} (-1) (-1) + \ha \cdot \frac{-2}{3/2} (-1) (3/4)+ \ha \cdot \frac{3/2}{-2}(-1) (-1) = -5/8.
\end{split}
\eeq
Summing up the symbols we get $\sigma(\mch_{23}) = 4c_\pi(-7/8).$ \\

To conclude, we obtained 
\beq
 \sigma(\mch_1({}^1\vec A))(q_0, \zeta) = c_\pi
  \begin{pmatrix}
    * & * & -5 & -5.5 \\
    * & * & 7 & 3.5 \\
-5 & 7 & * & -3.5\\
-5.5 & 3.5 & -3.5& * 
  \end{pmatrix}
\eeq

%%%%%%%%%%%%%%%
%\subsection{Computation for another four set of vectors}\label{secapp2}
For the proof of Prop.\ \ref{symb} part (3), we need to find another four set of vectors ${}^a\vec A, a = 2, 3, 4, 5$ so that the principal symbols  $\sigma(\mch_1({}^a\vec A)), a = 1, 2, 3, 4, 5$ are linearly independent.  We will fix the choices of $\zeta^{(i)}, i = 1, 2, 3, 4$  and vary the choices of $\vec A$. This is reasonable because each $A^{(i)}, i = 1, 2, 3, 4$ varies in a $3$ dimensional vector space (we lose one dimension due to the gauge condition). In the following, we shall  omit the details of the computation as they are the same as or ${}^1\vec A$. The result can be checked quite straightforwardly with some tedious computations. Also, one can check the result via symbolic computations using the formulas \eqref{compfor}. 

We take ${}^2\vec A$ as
\beq
\begin{gathered}
{}^2A^{(1)}= -(0,\imath,0, 0),\ \ {}^2A^{(2)}= -(0,0,\imath,0), \ \ 
{}^2A^{(3)}=  -(0, 0, 0,\imath),\ \ {}^2A^{(4)} =  -(0,0,0,\imath).
\end{gathered}
\eeq
Notice that compare to ${}^1\vec A$, we just changed $A^{(1)}$. We find 
\beq
\sigma(\mch_1({}^2\vec A)) = c_\pi
  \begin{pmatrix}
    * & * & 3.5 & 3.5 \\
    * & * & -3.5 & -5.5 \\
3.5 & -3.5 & * & -7\\
3.5 & -5.5 & -7& * 
  \end{pmatrix}
\eeq

We take ${}^3\vec A$ as
\beq
\begin{gathered}
{}^3A^{(1)}= -(0,0,0,\imath),\ \ {}^3A^{(2)}= -(0,\imath,0,0), \ \ 
{}^3A^{(3)}=  -(0, 0, 0,\imath),\ \ {}^3A^{(4)} =  -(0,0,0,\imath).
\end{gathered}
\eeq
Notice that compare to ${}^1\vec A$, we just changed $A^{(2)}$. We find 
\beq
\sigma(\mch_1({}^3\vec A)) = c_\pi
  \begin{pmatrix}
    * & * & 7 & 7 \\
    * & * & -5 & -7 \\
7 & -5 & * & 7\\
7 & -7 & 7& * 
  \end{pmatrix}
\eeq

We take ${}^4\vec A$ as
\beq
\begin{gathered}
{}^4A^{(1)}= -(0,0,0,\imath),\ \ {}^4A^{(2)}= -(0,0,\imath,0), \ \ 
{}^4A^{(3)}=  -(0, 0, \imath,0),\ \ {}^4A^{(4)} =  -(0,0,0,\imath).
\end{gathered}
\eeq
Here we just changed $A^{(3)}$ compared to ${}^1\vec A$. We find 
\beq
\sigma(\mch_1({}^4\vec A)) = c_\pi
  \begin{pmatrix}
    * & * & 5 & -5.5 \\
    * & * & 0 & -3.5 \\
5& 0 & * & 1\\
-5.5 & -3.5 & 1& * 
  \end{pmatrix}
\eeq

We take ${}^5\vec A$ as
\beq
\begin{gathered}
{}^5A^{(1)}= -(0,0,0,\imath),\ \ {}^5A^{(2)}= -(0,0,\imath,0),  \ \ 
{}^5A^{(3)}=  -(0, 0, 0,\imath),\ \ {}^5A^{(4)} =  -(0,0,\imath,0).
\end{gathered}
\eeq
Here we changed $A^{(4)}$ compared to ${}^1\vec A$. We find 
\beq
\sigma(\mch_1({}^5\vec A)) =  c_\pi
  \begin{pmatrix}
    * & * & -5.5 & 5 \\
    * & * & -3.5 & 0 \\
-5.5 & -3.5 & * & 1\\
5 & 0 & 1& * 
  \end{pmatrix}
\eeq

This completes the calculation for $\sigma(\mch_1)$ needed in the proof of Prop.\ \ref{symb}.

\subsection{Computation of $\sigma(\mch_2)$}
The computation of this term is relatively simple. We shall use the matrix form. 
For convenience, we let 
\beq
\mcw^{(ij)} = \bfq\bigg(\widehat H_2(x, \vec v^{(i)}, \vec v^{(j)}) +  \widehat H_2(x, \vec v^{(j)}, \vec v^{(i)})\bigg), \ \  1 \leq i < j \leq 4.
\eeq
Then we can write 
\beqq\label{symsim}
\begin{gathered}
\sigma(\mch_2)(q_0, \zeta) = - \sum  \sigma \bigg(\widehat P_2(x,  \mcw^{(ij)},  \mcw^{(kl)}) \bigg)(q_0, \zeta).
\end{gathered} 
\eeqq
and the summation is over the pairs $(i, j) = (1, 2), (1, 3), (1, 4), (2, 3), (2, 4), (3, 4). $ So we reduce to summation of six terms. 
Using the expressions of principal symbols in Section \ref{singu}, we find that
\beq
\begin{gathered}
-\sigma \bigg(\widehat P_2(x,  \mcw^{(12)},  \mcw^{(34)}) \bigg)(q_0, \zeta)   =  (H \sigma(\mcw^{(12)}) H)^{pq} (\zeta^{(3)}_p + \zeta^{(4)}_p) (\zeta^{(3)}_q + \zeta^{(4)}_q) \sigma(\mcw^{(34)}) \\
 = \frac{c_\pi}{2h(\zeta^{(1)}, \zeta^{(2)})}\bigg(H \big( 2 \sigma(F^{(1)}) H \sigma(F^{(2)}) + 2 \sigma(F^{(2)}) H \sigma(F^{(1)}) + H \text{Tr}(\sigma(F^{(1)})\sigma(F^{(2)}))\big) H\bigg)^{pq} \\
 \cdot (\zeta^{(3)}_p + \zeta^{(4)}_p) (\zeta^{(3)}_q + \zeta^{(4)}_q)  \\
 \cdot \frac{1}{2h(\zeta^{(3)}, \zeta^{(4)})} \big( 2 \sigma(F^{(3)}) H \sigma(F^{(4)}) + 2 \sigma(F^{(4)}) H \sigma(F^{(3)}) + H \text{Tr}(\sigma(F^{(3)})\sigma(F^{(4)}))\big).
\end{gathered}
\eeq
The other terms in \eqref{symsim} can be written in a similar way and are omitted here. So the computation is reduced to several matrix multiplications. Again, we do the calculation for the first set of vectors. For our choices of ${}^1\vec A$, we have
\beq
\begin{gathered}
\sigma(F^{(1)})%(q_0, \zeta^{(1)}) 
= \begin{pmatrix}
    0 & 0 & 0 & 1 \\
    0 & 0 & 0 & 0 \\
    0 & 0 & 0 & 1 \\
    -1 & 0 & -1 & 0 \\
\end{pmatrix}, \ \ 
\sigma(F^{(2)})%(q_0, \zeta^{(2)})  
= \begin{pmatrix}
    0 & 0 & 1 & 0 \\
    0 & 0 & 0 & 0 \\
    -1 & 0 & 0 & -1 \\
    0 & 0 & 1 & 0 \\
\end{pmatrix}, \\ 
\sigma(F^{(3)})%(q_0, \zeta^{(3)}) 
=  \begin{pmatrix}
    0 & 0 & 0 & -\frac 34 \\
    0 & 0 & 0 &  -\frac 34 \\
    0 & 0 & 0 &  0\\
    \frac 34 & \frac 34& 0 & 0 \\
\end{pmatrix}, \ \ 
\sigma(F^{(4)})%(q_0, \zeta^{(4)})  
=  \begin{pmatrix}
    0 & 0 & 0 & 1\\
    0 & 0 & 0 & -1 \\
    0 & 0 & 0 & 0 \\
    -1 & 1 & 0 & 0 \\
\end{pmatrix}.
\end{gathered}
\eeq
By straightforward computations, we get that 
\beq
\begin{gathered}
-\sigma \bigg(\widehat P_2(x,  \mcw^{(12)},  \mcw^{(34)}) \bigg)
 = c_\pi \begin{pmatrix}
    0 &  0 & 0 &  0\\
   0  & 0 & 0 &  0\\
  0  & 0&  3&  0\\
   0  & 0 & 0 & -3
   \end{pmatrix}, \ \   -\sigma \bigg(\widehat P_2(x,  \mcw^{(34)},  \mcw^{(12)}) \bigg)(q_0, \zeta)
 = 0
\end{gathered}
\eeq
Similarly, we have 
\beq
\begin{gathered}
-\sigma \bigg(\widehat P_2(x,  \mcw^{(13)},  \mcw^{(24)}) \bigg) 
 = c_\pi \begin{pmatrix}
    0 &  0 & -8 &  0\\
   0  & 0 & 8 &  0\\
  -8  & 8&  0& -8\\
   0  & 0 & -8 & 0
   \end{pmatrix}, \\
   -\sigma \bigg(\widehat P_2(x,  \mcw^{(24)},  \mcw^{(13)}) \bigg) 
 = c_\pi \begin{pmatrix}
    -1 &  -1 & -1 &  0\\
   -1  & -1 & -1 &  0\\
  -1  & -1&  -1& 0\\
   0  & 0 & 0 & 1
   \end{pmatrix},
   \end{gathered}
\eeq
\beq
\begin{gathered}
   -\sigma \bigg(\widehat P_2(x,  \mcw^{(14)},  \mcw^{(23)}) \bigg) 
 = c_\pi \begin{pmatrix}
    0 &  0 &  0.75&  0\\
   0  & 0 & 0.75 &  0\\
 0.75  & 0.75&  0& 0.75\\
   0  & 0 & 0.75 & 0
   \end{pmatrix}, \\
      -\sigma \bigg(\widehat P_2(x,  \mcw^{(23)},  \mcw^{(14)}) \bigg) 
 =  c_\pi \begin{pmatrix}
    6 &  -6 & 6 &  0\\
   -6 & 6 & -6 &  0\\
  6  & -6&  6& 0\\
   0  & 0 & 0 & -6
   \end{pmatrix}, 
   \end{gathered}
\eeq
So finally, for this set of vectors  ${}^1\vec A$, we denote the term $\mch_2$ by $\mch_2({}^1\vec A)$ and we find that 
\beq
\begin{gathered}
 \sigma(\mch_2({}^1 \vec A))(q_0, \zeta)
 = c_\pi \begin{pmatrix}
    5 &  -7 & -2.25 &  0\\
   -7 & 5 & 1.75 &  0\\
  -2.25  & 1.75&  8&  -7.25\\
   0  & 0 & -7.25 & -8
   \end{pmatrix}.
\end{gathered}
\eeq
For the other sets of vectors, we get
\beq
\begin{gathered}
    \sigma(\mch_2({}^2 \vec A))(q_0, \zeta)
 = c_\pi \begin{pmatrix}
    0 &  0 & -7 &  7\\
   0  & 0 & 5 &  -5\\
 -7 & 5&  -0.875& 0\\
   7  & -5 & 0 & 0.875
   \end{pmatrix}\\
    \sigma(\mch_2({}^3 \vec A))(q_0, \zeta)
 = c_\pi \begin{pmatrix}
    0 &  0 & -8.75 &  8.75\\
   0  & 0 & 7.25 &  -7.25\\
 -8.75 & 7.25&  -18.375& 0\\
   8.75  & -7.25 & 0 & 18.375
   \end{pmatrix}\\
    \sigma(\mch_2({}^4 \vec A))(q_0, \zeta)
 = c_\pi \begin{pmatrix}
 9.25 &  -8.75 & 4 &  2.25\\
   -8.75  & 5.25 & -2 &  0.25\\
 4  & -2&  10.75& -6\\
   -2.25  & 0.25 & -6 & -6.75
   \end{pmatrix}\\
   \sigma(\mch_2({}^5 \vec A))(q_0, \zeta)
 = c_\pi \begin{pmatrix}
9.25 &  -8.75 & 2.25 & 4\\
    -8.75  & 5.25 & 0.25 & -2\\
   2.25  & 0.25 &  -6.75 &  -6\\
 4  & -2 &-6 & 10.75
   \end{pmatrix}
\end{gathered}
\eeq

\subsection{Computation of $\sigma(\mch_3)$} 
We simplify the symbols as following. We know that 
\beq
\begin{gathered}
-\sigma(\widehat H_2(x, \vec v^{(i)}, \bfq(P_2(x, \bfq(\widehat H_2(x, \vec v^{(j)}, \vec v^{(k)})), \vec v^{(l)}))))\\
= \frac{c_\pi}{|\zeta^{(j)} + \zeta^{(k)}|_h^2 |\zeta^{(j)} + \zeta^{(k)} + \zeta^{(l)}|_h^2} \bigg(H \sigma(\widehat H_2(\vec v^{(j)}, \vec v^{(k)}))H\bigg)^{pq} \zeta^{(l)}_p\zeta^{(l)}_q \cdot \sigma(\widehat H_2(\vec v^{(i)}, \mcv^{(l)})),
\end{gathered}
\eeq
where
\beq
\begin{gathered}
 \sigma(\widehat H_2(v^{(i)}, \mcv^{(l)})) \doteq  2 \sigma(F^{(i)}) H \sigma(\mcf^{(l)}) + \ha H \text{Tr}(\sigma(F^{(i)})\sigma(\mcf^{(l)})), \\
\sigma(\mcf_{\alpha\beta}^{(l)}) =  \imath \tilde \zeta^{(l)}_\alpha  \sigma(\dot \phi^{(l)}_\beta) -  \imath \tilde \zeta^{(l)}_\beta \sigma(\dot \phi^{(l)}_{\alpha}),\ \ \tilde \zeta^{(l)} = \zeta^{(j)} + \zeta^{(k)} + \zeta^{(l)}.
\end{gathered}
\eeq
Notice that the principal symbols of $\mcf^{(l)}$ is evaluated at $(q_0, \tilde \zeta^{(l)})$. The other piece is similar
\beq
\begin{gathered}
-\sigma(\widehat H_2(x, \bfq(P_2(x, \bfq(\widehat H_2(x, \vec v^{(j)}, \vec v^{(k)})), \vec v^{(l)})) , \vec v^{(i)}))\\
= \frac{c_\pi}{|\zeta^{(j)} + \zeta^{(k)}|_h^2 |\zeta^{(j)} + \zeta^{(k)} + \zeta^{(l)}|_h^2} \bigg(H \sigma(\widehat H_2(\vec v^{(j)}, \vec v^{(k)}))H\bigg)^{pq} \zeta^{(l)}_p\zeta^{(l)}_q \cdot \sigma(\widehat H_2(\mcv^{(l)}, \vec v^{(i)})).
\end{gathered}
\eeq
To simplify the computation a bit, we set 
\beq
\begin{split}
\sigma(\mco^{il}) &= \sigma(\widehat H_2(v^{(i)}, \mcv^{(l)}) + \widehat H_2(\mcv^{(i)}, v^{(l)}))\\
 &= 2 \sigma(F^{(i)}) H \sigma(\mcf^{(l)}) +2 \sigma(\mcf^{(i)}) H \sigma(F^{(l)}) +  H \text{Tr}(\sigma(F^{(i)})\sigma(\mcf^{(l)})),
\end{split}
\eeq
so that 
\beq
\begin{gathered}
\sigma(\mch_3) =   \sum  \frac{c_\pi}{|\zeta^{(j)} + \zeta^{(k)}|_h^2 |\zeta^{(j)} + \zeta^{(k)} + \zeta^{(l)}|_h^2} \bigg(H \sigma(\mcw^{jk})H\bigg)^{pq} \zeta^{(l)}_p\zeta^{(l)}_q \cdot \sigma(\mco^{il}),
\end{gathered}
\eeq
where the summation is over $l, i = 1, 2, 3, 4, l\neq i$ and $(j, k)$ are determined by the choice of $l, i$. The computation is straightforward and one can use many terms $\mcw^\bullet$ already computed in $\mch_2$. However, we need to compute twelve terms. The details will be omitted. For the first choice of the vectors ${}^1\vec A$, we get 
\beq
\begin{gathered}
 \sigma(\mch_3({}^1 \vec A))(q_0, \zeta)
 = c_\pi \begin{pmatrix}
    -4 &  7 & 1.95 &  4.8\\
   7  & -14 & -6.65 &  1.4\\
  1.95  &-6.65&  1&  6.25\\
   4.8 & 1.4& 6.25 & 9
   \end{pmatrix}.
\end{gathered}
\eeq
For the other sets of vectors, we get
\beq
\begin{gathered}
    \sigma(\mch_3({}^2 \vec A))(q_0, \zeta)
 = c_\pi \begin{pmatrix}
    -5.6 &  9.8 & 3.15 &  -3.5\\
   9.8  & -5.6 & -1.75 &  6.5\\
 3.15 & -1.75&  5.775& -3.5\\
   -3.5  & 6.5 & -3.5 & -5.775
   \end{pmatrix}\\
    \sigma(\mch_3({}^3 \vec A))(q_0, \zeta)
 = c_\pi \begin{pmatrix}
    -11.9 & 9.5 & 0 &  -10.15\\
   9.5  & -11.9 & 5.3 &  6.25\\
 0 & 5.3&  12.775& 0\\
   -10.15  & 6.25 & 0 & -12.775
   \end{pmatrix}\\
    \sigma(\mch_3({}^4 \vec A))(q_0, \zeta)
 = c_\pi \begin{pmatrix}
 -6.25 &  8.75 & 0 &  -1.75\\
   8.75  &-6.25 & 1.6 &  3.65\\
 0 & 1.6&  -1.75& 4\\
   -1.75 & 3.65 & 4 & 1.75
   \end{pmatrix}\\
   \sigma(\mch_3({}^5 \vec A))(q_0, \zeta)
 = c_\pi \begin{pmatrix}
 -6.25 &  8.75 & -1.75 &  0\\
   8.75  &-6.25 & 3.65 &  1.6\\
 -1.75 & 3.65&  1.75& 4\\
   0 & 1.6 & 4 & -1.75
   \end{pmatrix}
   \end{gathered}
\eeq

Finally, we sum up our computation results to get $\sigma(\mch({}^a \vec A)) = \sum_{i = 1}^3 \sigma(\mch_i({}^a \vec A)), a = 1, 2, 3, 4, 5.$ Recall that it suffices to show that the vectors $\mct^a \doteq T(\sigma(\mch({}^a \vec A))), a = 1, 2,3,4, 5$ are linearly independent where $T$ is defined in the proof of Prop.\ \ref{symb}. So we have 
\beq
\begin{split}
T(\sigma(\mch({}^1 \vec A))) &= c_\pi (-5.3, -0.7, 2.1, 4.9, -4.5),\\
T(\sigma(\mch({}^2 \vec A))) &= c_\pi (-0.35, 7, -0.25, -4, -10.5), \\
T(\sigma(\mch({}^3 \vec A))) &= c_\pi (-1.75, 5.6, 7.55, -8, 7), \\
T(\sigma(\mch({}^4 \vec A))) &= c_\pi (9, -5, -0.4, 0.4, -1), \\
T(\sigma(\mch({}^5 \vec A))) &= c_\pi (-5, 9, 0.4, -0.4, -1).
\end{split}
\eeq
We check that the rank of these five vectors is $5$. 
This completes the proof of the claim in Prop.\ \ref{symb}.

%===============================REFERENCE==========================================%

\end{document}